\RequirePackage{fix-cm}
\documentclass[twocolumn]{svjour3}          
\smartqed  
\usepackage{graphicx}
\usepackage{amsmath}
\usepackage{amssymb} 
\usepackage{empheq}
\usepackage{xcolor}
\usepackage[centerlast,small,bf]{caption}
\usepackage{float}
\usepackage{subfig}
\usepackage[T1]{fontenc}
\usepackage{appendix}
\usepackage{supertabular}

\newcommand{\DDt}{\frac{\partial}{\partial t}}

\newcommand{\myvec}[1]{\mathbf{#1}}

\newcommand*\widefbox[1]{\fbox{\hspace{0.2em}#1\hspace{0.2em}}}

\newcommand\dashedrule{\mbox{%
  \rule[0.5ex]{1mm}{.75pt}\hspace{1mm}\rule[0.5ex]{1mm}{.75pt}\hspace{1mm}\rule[0.5ex]{1mm}{.75pt}\hspace{1mm}\rule[0.5ex]{1mm}{.75pt}}}

\definecolor{darkolivegreen}{rgb}{0.33, 0.42, 0.18}

\usepackage{tabularx}
\newcolumntype{L}[1]{>{\raggedright\arraybackslash}p{#1}}
\newcolumntype{C}[1]{>{\centering\arraybackslash}p{#1}}
\newcolumntype{R}[1]{>{\raggedleft\arraybackslash}p{#1}}
\journalname{Computational Geosciences}
\begin{document}

\title{Non-Isothermal, Multi-phase, Multi-component Flows through Deformable Methane Hydrate Reservoirs
\thanks{We gratefully acknowledge the support for the first author by the German Research Foundation (DFG), through project no. WO 671/11-1}
}


\author{Shubhangi Gupta	\and
        Rainer Helmig	\and 
        Barbara Wohlmuth 
}


\institute{S. Gupta \at
              Chair for Numerical Mathematics, \\ 
              Technical University Munich, \\
              Boltzmannstra\ss e 3, 85748 Garching bei M\"unchen, Germany\\
              Tel.: +49-89-28918439\\
              \email{gupta@ma.tum.de}           
           \and
           R. Helmig \at
              Dept. of Hydromechanics and Modelling of Hydrosystems \\
	      University of Stuttgart, \\
	      Pfaffenwaldring 61, 
	      70569 Stuttgart,
	      Germany 
           \and
           B. Wohlmuth \at
              Chair for Numerical Mathematics, \\ 
              Technical University Munich, \\
              Boltzmannstra\ss e 3, 85748 Garching bei M\"unchen Germany\\
}

\date{Received: date / Accepted: date}

\maketitle

\begin{abstract}
We present a hydro-geomechanical model for subsurface methane hydrate systems. 
Our model considers kinetic hydrate phase change and non-isothermal, multi-phase, multi-component flow in elastically deforming soils. 
The model accounts for the effects of hydrate phase change and pore pressure changes on the mechanical properties of the soil.
It also accounts for the effect of soil deformation on the fluid-solid interaction properties relevant to reaction and transport processes (e.g., permeability, capillary pressure, reaction surface area).
We discuss a 'cause-effect' based decoupling strategy for the model and 
present our numerical discretization and solution scheme. 
We then proceed to identify the important model components and couplings which are most vital for a hydro-geomechanical hydrate simulator, namely, 
1) dissociation kinetics, 
2) hydrate phase change coupled with non-isothermal two phase two component flow,
3) two phase flow coupled with linear elasticity (poroelasticity coupling), and finally 
4) hydrate phase change coupled with poroelasticity (kinetics-poroelasticity coupling).
To show the versatility of our hydrate model, we numerically simulate test problems where, for each problem, we methodically isolate one out of the four aforementioned model components or couplings.
A special emphasis is laid on the kinetics-poroelasticity coupling for which 
we present a test problem where an axially loaded hydrate bearing sand sample experiences a spontaneous shift in the hydrate stability curve causing the hydrate to melt.
For this problem we present an analytical solution for pore-pressure, which we subsequently use to test the accuracy of the numerical scheme. 
Finally, we present a more complex $3D$ example where all the major model components are put together to give an idea of the model capabilities.
The setting is based on a subsurface hydrate reservoir which is destabilized through depressurization using a low pressure gas well.
In this example, we simulate the melting of hydrate, methane gas generation, 
and the resulting ground subsidence and stress build-up in the vicinity of the well.

\keywords{Methane hydrate reservoir \and hydro-geomechanical model \and kinetics-poroelasticity coupling}
\subclass{MSC 74F10 \and MSC 74F25 \and MSC 76S05 \and MSC 76V05}

\end{abstract}

\section{Introduction}\label{sec:introduction}
\paragraph{}
Methane hydrates are formed when water molecules form a cage-like structure and trap a large number of methane molecules within, forming a crystalline solid similar to ice \cite{Sloan1998}. 
Methane hydrates are thermodynamically stable under conditions of low temperature and high pressure and occur naturally in permafrost regions or below ocean/sea floors \cite{DaweThomas2007}.
If warmed or depressurized, methane hydrates destabilize and dissociate into water and methane gas.  
Methane hydrates are a very dense source of methane gas. One cubic meter of methane hydrate stores approximately 164 standard cubic meters of methane gas.
Also, the energy content of methane occurring in hydrate form is immense, possibly even exceeding the combined energy content of all other conventional fossil fuels \cite{Milkov2004}. 

\paragraph{}
The research in methane hydrates stems from a 'three-fold' motivation driven by concerns of methane gas extraction and production feasibility, global warming potential and climate change concerns,  
and inherent geo-hazards of mining/drilling induced destabilization of methane hydrates in subsurface reservoirs.
Phenomenological modeling and numerical simulation of these systems is thus vital 
1) for optimizing recovery techniques for extracting methane from hydrate bearing sediments,
2) for conducting studies and making predictions for mitigating bore-hole, local and regional slope stability hazards, 
3) for sequestering carbon-dioxide in gas hydrate, 
4) for possible application in natural gas storage and transport, and 
5) for evaluating the role of gas hydrate in global carbon cycle.

\paragraph{}
Methane hydrate formations are a fairly complex subsurface system characterized by a large number of highly interdependent physical phenomena. 
The typical physical processes occuring in a stimulated hydrate reservoir include 
1) hydrate phase change, 
2) non-isothermal multi-phase, multi-component flow, 
3) geomechanical deformation of the hydrate bearing sediment, and 
4) change in the hydraulic as well as the mechanical properties of the hydrate bearing sediment.
Thus, any detailed study of these reservoirs and their possible applications in energy, environment, and quantification of geo-hazards requires the development of 
1) \textit{multi-physics mathematical model} that accounts for the afore-mentioned processes and captures their interdependencies, and 
2) robust, and efficient \textit{numerical tools} capable of handling multi-physics models and performing large scale simulations.

\paragraph{}
Several mathematical models have been proposed (e.g. Tsypkin \cite{Tsypkin1991}, Ahmadi et al. \cite{Ahmadi2004}, Yousif et al. \cite{Yousif1991}, 
Sun and Mohanty \cite{SunMohanty2005}, Liu and Flemmings \cite{LiuFlemming2007}, Moridis \cite{Moridis2003,Moridis2007})
and different numerical simulators have been developed
(e.g. MH21-HYDRES \cite{MH21-HYDRES}, STOMP-HYD \cite{STOMP-HYD}, UMSICHT HyRes \cite{UMSICHT-HyRes}, TOUGH-HYDRATE \cite{TOUGH-HYDRATEmanual}) for simulating hydrate reservoirs.
These models and simulators consider mainly hydrate phase change and fluid flow while neglecting the geomechanical effects. 

\paragraph{}
It has been widely recognized in the hydrate community that the destabilization of hydrates can cause significant consolidation and ground deformation, 
and in extreme cases it can even trigger landslides \cite{SultanCochonat2004}.
Thus, a lot of experimental work has been done in characterizing the mechanical properties and deformation behaviour of methane hydrates and hydrate-bearing sediments 
(e.g. Hyodo et al. \cite{Hyodo2005,Hyodo2014}, Lee et al. \cite{Lee2010}).
Several mathematical models have been proposed to extend the above mentioned hydrate-reservoir model concepts to include geo-mechanics 
(e.g. Rutqvist and Moridis \cite{RutqvistMoridis2008}, Kimoto et al. \cite{Kimoto2010}).
Rutqvist and Moridis \cite{RutqvistMoridis2008} and Klar et al. \cite{KlarSogaNg2010} have coupled TOUGH-HYDRATE with the commercial geomechanical code FLAC$^{3D}$ \cite{Itasca2009} 
to investigate the hydro-geomechanical behaviour of hydrate reservoirs.
Kimoto et al. \cite{Kimoto2010} have developed their own chemo-thermo-hydro-mechanical simulator. Their simulator uses an elasto viscoplastic model to simulate deformation. 

\paragraph{}
In our work, we focus on building a consistent mathematical and numerical framework for hydrate systems from \textit{'ground zero'} up using a multi-physics approach.
The primary objective is to capture the dynamic coupling between the transport and mechanical processes observed at macroscopic scales.
In this paper, we first present our mathematical model. 
The application in focus is limited to \textit{gas production} from thermally stimulated or de-pressurized reservoirs. 
We break the model down into its functional \textit{building blocks} and try to identify the mechanism of the information exchange between them.
This is important to establish a consistent feedback loop between the processes. It also helps us to clearly identify the various couplings present in the multi-physics model.
We then present the numerical solution scheme followed by test problems where we methodically isolate each of the couplings identified during the model break-down step.
Through these test problems we 1) verifying each of the model components making up the hydrate simulator, and 2) show the versatility of the model in the variety of hydrate reservoir related problems it can handle.
Finally, we present a $3D$ example problem where a typical subsurface hydrate reservoir is destabilized by depressurization.

\section*{Nomenclature}
    \begin{supertabular}{L{0.65cm} L{0.65cm} L{6.5cm}}
      $\chi^{\kappa}_{\alpha}$ 	& \multicolumn{2}{L{7.2cm}}{ Mole fraction of component $\kappa=CH_4,H_2O$ in phase $\alpha=g,w$ } \\
      $\dot g^{CH_4}$ 		& \multicolumn{2}{L{7.2cm}}{ $CH_4$ generation rate } \\ 
      $\dot g^{H_2O}$ 		& \multicolumn{2}{L{7.2cm}}{ $H_2O$ generation rate } \\ 
      $-\dot g^{Hyd}$, & $\ -\dot g^{h}$ & Hydrate consumption rate 			 \\
      $\dot Q_h$		& \multicolumn{2}{L{7.2cm}}{ Heat of hydrate phase change			} \\
      $\dot q_{m_\alpha}^{\kappa}$& \multicolumn{2}{L{7.2cm}}{ Volumetric injection rate for phase $\alpha=g,w$ } \\
      $\epsilon$ 		& \multicolumn{2}{L{7.2cm}}{ Volumetric or isotropic strain 			} \\
      $\mathbf{g}$ 		& \multicolumn{2}{L{7.2cm}}{ Gravity vector 					} \\ 
      $\mathbf{J}_{\alpha}^{\kappa}$ & \multicolumn{2}{L{7.2cm}}{ Diffusion flux of component $\kappa=CH_4,H_2O$ in phase $\alpha=g,w$ } \\
      $\mathbf{u}$ 		& \multicolumn{2}{L{7.2cm}}{ Sediment displacement vector 			} \\
      $\mathbf{v}_s$ 		& \multicolumn{2}{L{7.2cm}}{ Sediment displacement velocity 			} \\
      $\mathbf{v}_{\beta,t}$ 	& \multicolumn{2}{L{7.2cm}}{ Total phase velocity 				} \\ 
      $\mathbf{v}_{\beta}$ 	& \multicolumn{2}{L{7.2cm}}{ Phase velocity relative to the sediment 		} \\
      $\mu_{\alpha}$		& \multicolumn{2}{L{7.2cm}}{ Dynamic viscocity of phase $\alpha=g,w$ 		} \\
      $\nu_{sh}$ 		& \multicolumn{2}{L{7.2cm}}{ Poisson ratio for composite solid 			} \\
      $\phi$ 			& \multicolumn{2}{L{7.2cm}}{ Actual or total porosity 				} \\ 
      $\phi_{eff}$		& \multicolumn{2}{L{7.2cm}}{ Apparent or ffective porosity 			} \\
      $\rho_{\gamma}$ 		& \multicolumn{2}{L{7.2cm}}{ Density of phase $\gamma=g,w,h,s$ 			} \\ 
      $\rho_{sh}$ 		& \multicolumn{2}{L{7.2cm}}{ Density of composite solid 			} \\
      $\sigma \ \ ,$ & $ \sigma' $ & Isotropic total and effective stresses 					  \\
      $\tau$ 			& \multicolumn{2}{L{7.2cm}}{ Tortuosity 					} \\
      $\tilde \epsilon$ 	& \multicolumn{2}{L{7.2cm}}{ strain tensor 					} \\
      $\tilde \sigma \ \ $, & $\tilde \sigma'$ & Total and effective stress tensors 		 		  \\
      $A_{rs}$ 			& \multicolumn{2}{L{7.2cm}}{ Specific reaction area  				} \\
      $A_{s}$ 			& \multicolumn{2}{L{7.2cm}}{ Specific surface area of hydrate-bearing sediment 	} \\
      $B_{\gamma}$		& \multicolumn{2}{L{7.2cm}}{ Bulk modulus of phases $\gamma=g,w,h,s$ 		} \\
      $B_{m}$ 			& \multicolumn{2}{L{7.2cm}}{ Bulk modulus of bulk REV  				} \\
      $B_{sh}$ 			& \multicolumn{2}{L{7.2cm}}{ Bulk modulus of composite matrix			} \\
      $Cp_{\alpha}$ 		& \multicolumn{2}{L{7.2cm}}{ Specific heat capacity at constant pressure of phase $\alpha=g,w$ } \\  
      $Cv_{\gamma}$		& \multicolumn{2}{L{7.2cm}}{ Specific heat capacity at constant volume of phase $\gamma=g,w,h,s$ } \\
      $D^{\alpha} $ 		& \multicolumn{2}{L{7.2cm}}{ Binary diffusion constant in phase $\alpha=g,w$ 	} \\
      $E_h \ \ $,  & $E_s$ 	& Young's modulus for hydrate and soil						  \\
      $E_{sh}$ 			& \multicolumn{2}{L{7.2cm}}{ Young's modulus for composite solid 		} \\
      $f_g$ 			& \multicolumn{2}{L{7.2cm}}{ Gas phase fugacity 				} \\
      $G_{sh} \ \ , $ & $\lambda_{sh}$ 	& Lame's parameters					 		  \\
      $h_{\alpha} $ 		& \multicolumn{2}{L{7.2cm}}{ Flow enthalpy of phase $\alpha=g,w$ 		} \\
      $H$	 		& \multicolumn{2}{L{7.2cm}}{ Henry's constant for methane 			} \\
      $K$ 			& \multicolumn{2}{L{7.2cm}}{ Intrinsic permeability of hydrate bearing sediment } \\
      $k^c_{\gamma}$ 		& \multicolumn{2}{L{7.2cm}}{ Thermal conductivity of phase $\gamma=g,w,h,s$	} \\ 
      $k^c_{eff}$ 		& \multicolumn{2}{L{7.2cm}}{ Lumped thermal conductivity of the REV 		} \\ 
      $k_{d}$ 			& \multicolumn{2}{L{7.2cm}}{ Hydrate dissociation rate constant 		} \\ 
      $k_{r\alpha}$ 		& \multicolumn{2}{L{7.2cm}}{ Relative permeability of phase $\alpha=g,w$ 	} \\
      $k_{reac}$ 		& \multicolumn{2}{L{7.2cm}}{ Rate constant for kinetic phase change of hydrate 	} \\
      $M_{\kappa}$		& \multicolumn{2}{L{7.2cm}}{ Molar mass of component $\kappa=H_2O,CH_4,Hyd$ 	} \\
      $N_{Hyd}$, & $N_h$ 		& Hydration number 					 		  \\
      $P^{sat}_{H_2O}$ 		& \multicolumn{2}{L{7.2cm}}{ Saturated water vapor pressure 			} \\
      $P_{\alpha}$		& \multicolumn{2}{L{7.2cm}}{ Pressure of phase $\alpha=g,w$ 			} \\
      $P_c$			& \multicolumn{2}{L{7.2cm}}{ Capillary pressure 					} \\
      $P_{e}\ \ $, & $P_{eqb}$	&  Equilibrium pressure for hydrate phase 		 \\
      $P_{eff}$			& \multicolumn{2}{L{7.2cm}}{ Effective fluid pressure 				} \\
      $R_u$			& \multicolumn{2}{L{7.2cm}}{ Universal gas constant				} \\
      $S_{\alpha r}$ 		& \multicolumn{2}{L{7.2cm}}{ Effective aqueous phase saturation			} \\ 
      $S_{\beta}$		& \multicolumn{2}{L{7.2cm}}{ Saturation of phase $\beta=g,w,h$ 			} \\       
      $T$ 			& \multicolumn{2}{L{7.2cm}}{ Temperature 						} \\ 
      $u_{\gamma}$ 		& \multicolumn{2}{L{7.2cm}}{ Internal energy of phase $\gamma=g,w,h,s$ 	} \\ 
    \end{supertabular}

\section{Mathematical model} \label{sec:mathematicalModel}

\begin{figure*}
  \centering
    \includegraphics[scale=0.45]{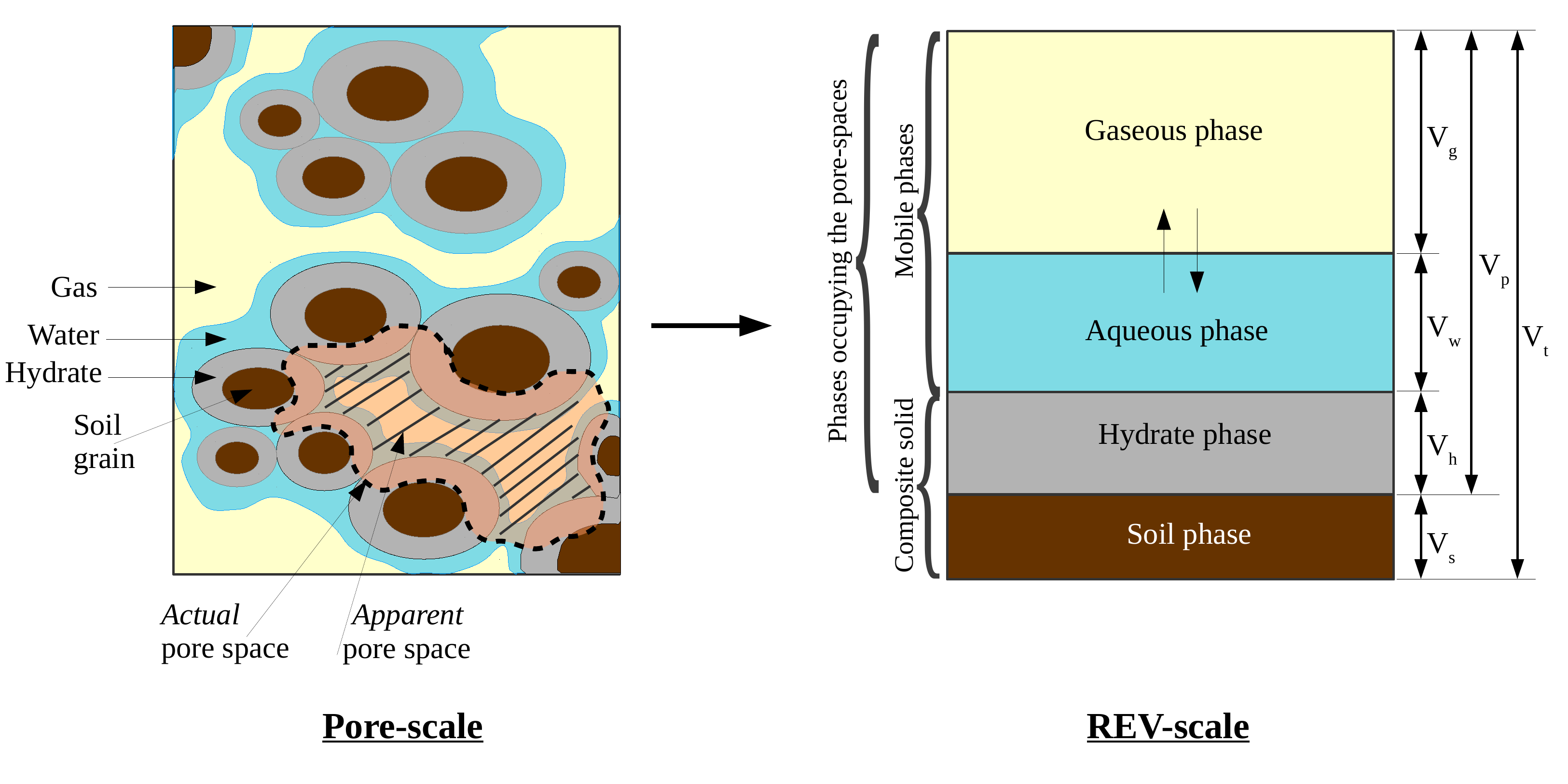}
  \caption{Pore-scale to REV-scale}  \label{fig:REVDescription}
\end{figure*}

\paragraph{} We consider three \textit{molecular components}: $CH_4$, $H_2O$ and $CH_4.\left(H_2O\right)_{N_h}$ (Hydrate), 
which are present in three distinct \textit{phases}: gaseous, aqueous, and solid.
The gaseous phase comprises of molecular methane and molecular water in vapour form. 
The aqueous phase comprises of molecular water and dissolved molecular methane.
The solid phase comprises of pure methane hydrate and soil grains.
The soil grains are assumed to form a material continuum which provides the skeletal structure to the porous medium. We shall refer to this as \textit{solid-matrix}.
The aqueous, gaseous, and hydrate phases exist in the void spaces of this solid matrix (See Fig. \ref{fig:REVDescription}).
At this stage the adsorption of methane gas on the surfaces of the solid matrix and the hydrate is not considered.

\paragraph{} At the pore-scale, we make a distinction between the \textit{actual pore-space}, which is the void space outside the soil-grains, 
and the \textit{apparent pore-space} which is that part of the void space that is not occupied by hydrate and is thus available for flow of water and gas.
At the REV-scale, this translates to \textit{actual} or \textit{total porosity} $\ \phi=\dfrac{V_p}{V_t}$, 
and \textit{apparent} or \textit{effective porosity} $\ \phi_{eff}=\dfrac{V_p-V_h}{V_t}$ (See Fig. \ref{fig:REVDescription}).
This distinction is important as it gives us a conceptual advantage in isolating the effects of deformation and hydrate melting on the hydraulic properties of the porous medium. 
(See Section \ref{sec:hydraulicProperties})

\paragraph{} To describe the hydraulic properties in Section \ref{sec:hydraulicProperties} and the effect of hydrate melting on these properties, 
we make an additional assumption that at the pore-scale the hydrate coats the soil grain perfectly, and the water phase forms a film over the hydrate. 

\paragraph{} In the subsequent discussion, the phases occupying the pore space (gaseous, aqueous, hydrate) will be denoted by '$\beta $' $ = g, w, h $ respectively, 
the mobile phases (gaseous and aqueous) will be denoted by '$\alpha$' $ = g, w $, 
and the molecular components will be denoted by '$\kappa$' $ = CH_4, H_2O, Hyd $. 
The soil matrix will be designated with the subscript '$s$'. 
The soil+hydrate composite matrix will be designated with the subscript '$sh$'.
'$\gamma$' would be used to denote all phases, i.e., $\gamma = g,w,h,$ and $s$.

\paragraph{\textbf{Assumptions} } The definition of the mathematical model is based on the following set of assumptions.
\begin{itemize}
 \item Gas hydrate reservoir is of SI type, i.e., the gas caged in the hydrate is purely methane.
 \item Ice formation is neglected.
 \item Aqueous phase is free of any salinity.
 \item Local thermal equilibrium (LTE) condition prevails in a representative control volume (REV) so that the solid and the neighboring fluid are at the same temperature.
 \item Gaseous and aqueous phases are treated as ideal mixtures.
 It is further assumed that the chemical species in the fluid phases attain chemical equilibrium instantaneously.
 It is important to emphasize that the chemical equilibrium exists only between the fluid phases, not between the fluids and hydrate. 
 The hydrate phase change is a relatively slower process compared with the methane dissolution and water evaporation processes.
 \item Flow velocities of mobile gaseous and aqueous phases belong to Darcy's Law regime.
 \item Hydrate is of \textit{grain-coating} type and remains perfectly adhered to the pore walls in the soil matrix.
 \item Soil and hydrate form a composite material. 
 This material (composite-solid-matrix) is treated as a continuum phase, and the stresses are considered to be acting on this composite-solid-matrix as a whole (and not on the soil matrix alone).
 The mechanical behaviour of this composite material is described using a linear elastic stress-strain constitutive law.
 The mechanical properties of this composite material are assumed to vary with the composition and stress-state of the composite-solid-matrix.
 
\end{itemize}

\subsection{Governing equations} \label{sec:governingEquations}
\paragraph{} The transport processes characterizing the gas production from a typical sub-surface methane hydrate reservoir can be described by invoking the conservation laws for mass, momentum, and energy
described for the macroscale properties of the porous medium \cite{Helmig1997}.

\paragraph{\textbf{Mass conservation} } The mass conservation for water and gas is written componentwise, i.e., for each component $\kappa = CH_4, H_2O$, we get
\begin{align}\label{eqn:massBalance_CH4_H2O}
       &        \sum\limits_{\alpha} \left[ \partial_t \left( \phi \ \rho_{\alpha} \  \chi_{\alpha}^{\kappa} \ S_{\alpha} \right) \right]
       \ + \ 	\sum\limits_{\alpha} \left[ \nabla\cdot\left( \phi \ \rho_{\alpha} \  \chi_{\alpha}^{\kappa} \ S_{\alpha} \ {\mathbf v_{\alpha,t}} \right) \right] \notag \\
        = &	\sum\limits_{\alpha} \left[ \nabla\cdot\left(\phi \ S_{\alpha} \ {\bf J}_{\alpha}^{\kappa}\right) \right]  
       \ + \ 	\dot g^{\kappa}
       \ + \ 	\sum\limits_{\alpha} \dot q_{m_{\alpha}}^{\kappa} \ .
\end{align}

\paragraph{} The mass conservation for the hydrate phase is given by
\begin{align}\label{eqn:massBalance_Hydrate}
   \partial_t \left(\phi \ \rho_{h} \  S_{h} \right)
 \ + \ \nabla\cdot\left(\phi \ \rho_{h} \  S_{h} \ {\mathbf v_{h,t}} \right)
 \ = \ \dot g^{h}  \ .
\end{align}

\paragraph{} The mass conservation for the soil phase is given by
\begin{align}\label{eqn:massBalance_Soil}
 \partial_t \left[(1-\phi) \ \rho_{s} \right]
 \ + \ \nabla\cdot\left((1-\phi) \ \rho_{s} \ {\mathbf v_{s}} \right)
 \ = \ 0  \ .
\end{align}

\paragraph{\textbf{Momentum conservation} } For mobile phases in porous medium, under certain simplifying assumptions, the momentum conservation can be reduced to Darcy's Law \cite{Hassanizadeh1986}, 
which is stated as,
\begin{align}\label{eqn:darcysLaw}
  {\bf v}_{\alpha}=-K \ \frac{k_{r\alpha}}{\mu_{\alpha}}\left(\nabla P_\alpha - \rho_{\alpha} \ {\bf g}\right)  \ .
\end{align} 

\paragraph{} Here, $K$ is the intrinsic permeability of the composite matrix, $k_{r,\alpha}$ and $\mu_{\alpha}$ are the relative permeability and the dynamic viscosity of the phase $\alpha$ respectively. 
\newline
${\bf v}_{\alpha}$ is the velocity of the mobile phase relative to the soil matrix.
The hydrate phase is immobile relative to the soil-matrix, i.e., $\mathbf{v}_{h} = 0$.
The total velocity of any phase $\beta$ occupying the pores is given by $\ \phi S_{\beta} {\bf v}_{\beta,t} = {\bf v}_{\beta} + \phi S_{\beta} {\bf v}_{s}$.
The soil phase velocity is the rate of deformation of the soil matrix, and is given by $\ \mathbf{v}_s = \partial_t \mathbf{u}$.

\paragraph{} Momentum conservation for the composite solid matrix is given by
\begin{align}\label{eqn:momentumBalance_soil-1}
\nabla \cdot \tilde{\sigma} + \rho_{sh} \ {\bf g} = 0 \ .
\end{align}

      
\paragraph{\textbf{Energy conservation} } For describing the energy conservation in the porous medium, 
one energy balance equation is sufficient since local thermal equilibrium has been assumed \cite{Helmig1997}. 
The energy balance equation is thus given by
\begin{align}\label{eqn:energyBalance}
    & \partial_t \left[ \left( 1-\phi \right) \rho_s u_s + \sum_{\beta} \left( \phi \ S_{\beta} \ \rho_{\beta} \ u_{\beta} \right)  \right] \notag \\
    & + \ \sum_{\alpha} \left[ \nabla \cdot \left( \phi \ \rho_{\alpha} \ \chi_{\alpha}^{\kappa} \ S_{\alpha} \ {\bf v }_{\alpha,t} \ h_{\alpha} \right) \right] \notag \\
    & = \ \nabla \cdot k^c_{eff} \nabla T 
    \ + \ \dot Q_h 
    \ + \ \sum\limits_{\alpha} \left( \dot q_{m_{\alpha}}^{\kappa} \  h_{\alpha} \right) 
\end{align}
where, 
\begin{align}
  & k^c_{eff} = \left( 1 - \phi \right) k^c_s + \sum_{\alpha}\sum_{\kappa} \left( \phi \ \chi_{\alpha}^{\kappa} \ S_{\alpha} \ k_{\alpha}^c \right) + \phi \ S_{h} \ k_h^c \notag \\
  & h_{\alpha} = \int_{T_{ref}}^T Cp_{\alpha}\ dT 	\notag\\
  & u_{\gamma} = \int_{T_{ref}}^T Cv_{\gamma}\ dT  \ .	\notag
\end{align}

\paragraph{\textbf{Closure relations} } The saturations of the phases occupying the pores satisfy the summation condition \newline $\quad \sum_{\beta} S_{\beta} = 1$\ .
\newline
Additionaly, for each mobile phase $\alpha$, the constituting component mole-fractions also satisfy the summation condition $\quad \sum_{\kappa} \chi_{\alpha}^{\kappa} = 1$\ .

\paragraph{} Furthermore, the pressures of the fluid phases are related through a capillary pressure $P_c$ as $P_g - P_w = P_c $ . 
This pressure difference occurs across the gaseous and aqueous phase interface due to balancing of cohesive forces within the liquid 
and the adhesive forces between the liquid and soil-matrix.  
The parametrization used for approximating $P_c$ is further elaborated upon in Section \ref{sec:hydraulicProperties}.

\subsection{Constitutive relationships} \label{sec:constitutiveRelations}
\paragraph{} In the mathematical model described in Section \ref{sec:governingEquations} following variables can be identified,
\begin{align}
 S_{\beta} \ , \ \chi_{\alpha}^{\kappa}\ , \ P_{\alpha} \ , \ P_c \ , \ T \ , \ \phi \ , \ \tilde{\sigma} \ , \ {\bf u } \ , \ {\bf J}_{\alpha}^{\kappa} \ , \ {\bf v}_{\alpha} \ , \ \dot g^{\beta} \ , \ \dot Q_h
\end{align}
i.e., the total number of variables is $24$. (The vectors and tensors are considered as single variables and $\kappa$ denotes the mobile components, i.e. $\kappa=CH_4,H_2O$.)
However, the number of governing equations add up to only $12$.
To close this system $12$ additional constitutive relationships are defined for $\chi_{g}^{\kappa}\ , \ P_c \ , \ \tilde{\sigma} \ , \ {\bf J}_{\alpha}^{\kappa} \ , \ \dot g^{\beta} \ , \text{and} \ \dot Q_h$ in this section.
Some other properties which are important for modelling hydrate reservoirs are also discussed.

\subsubsection{Vapor-liquid equilibrium}
 \paragraph{} The two-phase $CH_4 - H_2O$ fluid system is assumed to be in a state of vapor-liquid-equilibrium. 
  To calculate the concentrations of $CH_4$ and $H_2O$ in both gaseous and aqueous phases Henry's Law \cite{AtkinsPhysicalChemistry} and Raoult's Law \cite{AtkinsPhysicalChemistry} for ideal gas-liquid solutions are invoked:
  \begin{align}
    &\text{For dissolved methane,}		\quad \chi_w^{CH_4} = H (T) \ \chi_g^{CH_4} \ P_g  \label{eqn:DissolvedMethane} \\
    &\text{For water vapor,}	\qquad \qquad \chi_g^{H_2 O} = \chi_w^{H_2 O} \ \frac{P^{sat}_{H_2 O} (T) }{P_g}	 \label{eqn:WaterVapour}
  \end{align}

\paragraph{} Using Eqn. \eqref{eqn:DissolvedMethane}, Eqn. \eqref{eqn:WaterVapour}, and the summation conditions for each $\alpha$, $\quad \sum_{\kappa} \chi_{\alpha}^{\kappa} = 1$\ , the mole fractions can be calculated explicitly.

\paragraph{} In Eqn. \eqref{eqn:DissolvedMethane}, $H$ is Henry's constant for methane gas dissolved in water 
and is calculted in our model using the empirical relation from NIST standard reference database \cite{NISTChemHenryConst}.
In Eqn. \eqref{eqn:WaterVapour}, $P^{sat}_{H_2 O}$ is the saturated vapour pressure for water in contact with methane gas. 
We use Antoine's equation \cite{AtkinsPhysicalChemistry} to calculate $P^{sat}_{H_2 O}$.

\subsubsection{Fick's law for diffusive mass-transfer flux}
      \paragraph{} The diffusive solute flux through sediment ${\bf J}_{\alpha}^{\kappa}$ is calculated using Fick's Law \cite{Hassanizadeh1986}, as stated below
        \begin{align}\label{eqn:FicksLaw}
         {\bf J}_{\alpha}^{\kappa} = - \tau D^{\alpha} \left( \rho_{\alpha} \nabla \chi_{\alpha}^{\kappa} \right) \ .
        \end{align}


  \paragraph{} The gas-phase binary diffusion coefficient for low density binary $CH_4 - H_2O$ system is estimated using the empirical relationship proposed by Stattery and Bird \cite{StatteryBird1958}.
  For the aqueous phase binary diffusion coefficient, Wilke-Chang correlation \cite{Himmelblau1964} for dilute associated liquid mixtures is used.

\subsubsection{Hydrate phase change kinetics}

\begin{figure}
  \centering
    \includegraphics[scale=0.5]{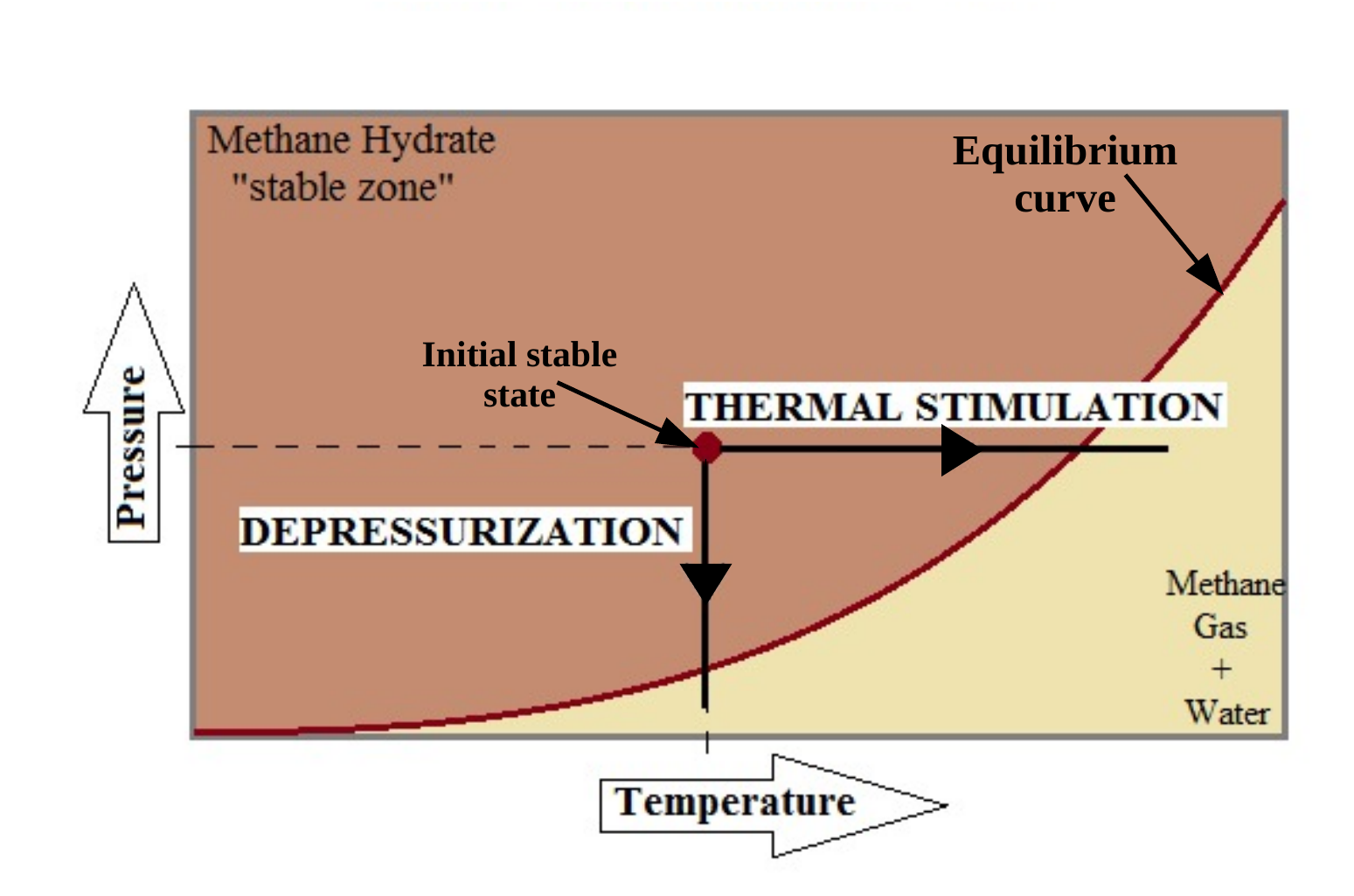}
  \caption{Methane hydrate P-T phase curve}  \label{fig:methaneHydrateEquilibriumCurve}
\end{figure}

\paragraph{} Methane hydrates, upon heating or depressurization, decompose to produce methane gas and water as shown in Fig. \ref{fig:methaneHydrateEquilibriumCurve}.
Chemically, this phase change process can be expressed as $\quad CH_4\cdot \left(H_2O\right)_{N_{Hyd}} \rightleftarrows CH_4 + N_{Hyd} \cdot H_2O$,
where, $N_{Hyd}$ is the hydration number.

\paragraph{} This non-equilibrium phase change of methane hydrate is modeled by the Kim-Bishnoi kinetic model \cite{KimBishnoi1987}. 
The rate of gas generated or consumed on hydrate decomposition or reformation is given by
\begin{align}\label{eqn:methaneGenerationRate}
 \dot g^{CH_4} = k_{reac}\ M_{CH4}\ A_{rs}\ \left( P_e - f_g \right) \ .
\end{align}

\paragraph{} Correspondingly, the rates of water and methane hydrate generated/consumed are given by
\begin{align}
 & \dot g^{H_2O} = \dot g^{CH_4} \ N_{Hyd} \ \frac{M_{H_2O}}{M_{CH_4}} \quad \text{, and } \notag \\
-& \dot g^{Hyd}  = \dot g^{CH_4}\ \frac{M_{Hyd}}{M_{CH_4}}  \ . \label{eqn:rateOfWaterAndHydrateGeneration}
\end{align}

\paragraph{} In Eqn. \eqref{eqn:methaneGenerationRate}, $k_{reac}$ is the kinetic rate constant given by 
\begin{align}\label{eqn:kineticRateConstant}
   k_{reac} = k_{reac}^0 \ \exp\left( - \frac{\Delta E_{a}}{R\ T} \right)
\end{align}
where, $E_{act}$ is the activation energy, and $k_{reac}$ is the intricsic rate constant for the kinetic phase change.

\paragraph{} $A_{rs}$ is the specific reaction area available for the kinetic reaction to occur, given by 
\begin{align}\label{eqn:reactionSurfaceArea}
 A_{rs} &= \Gamma_r A_s
\end{align}

\paragraph{} where, $\Gamma_r$ is the fraction of the pore surface area that is active in hydrate kinetics \cite{SunMohanty2006}.
The pore surface area $A_s$ is a hydraulic property of the solid matrix (Refer Sec. \ref{sec:hydraulicProperties}, Eqn. \ref{eqn:specificSurfaceArea}).


\paragraph{} The methane gas fugacity $f_g$ in the kinetic model is computed based on the Peng-Robinson's thermodynamic equation of state for methane \cite{PengRobinson1970}. 
The equilibrium pressure for the methane hydrate $P_e$ is determine using the Kamath and Holder correlation \cite{KamathHolder1987} as given below,
\begin{align}\label{eqn:hydrateEquilibriumPressure}
 P_e &= A_1 \ \exp\left( A_2 - \frac{A_3}{T \ [\text{K}]} \right)  \ .
\end{align}

\paragraph{} Methane hydrate dissociation reaction is an endothermic process. Conversely, methane hydrate re-formation is an exothermic reaction. 
The heat of reaction for hydrate phase change is modelled by
\begin{align}\label{eqn:hydrateHeatOfDissociation}
 \dot Q_h &= \frac{\dot g^{Hyd}}{M_{Hyd}}\left( B_1 - \frac{B_2}{T[\text{K}]} \right) \ .
\end{align}

\subsubsection{Properties of the fluid-matrix interaction}\label{sec:hydraulicProperties}
\paragraph{\textbf{Capillary pressure} } On the macro scale, the capillary pressure in a porous medium is an average pressure depending on the pore-size distribution and the aqueous phase saturation. 
Several parameterizations exist which relate the capillary pressure and effective aqueous phase saturation using soil specific parameters. 
Our model uses one of the most common parameterizations proposed by Brooks and Corey \cite{BrooksCorey1964}.

\paragraph{}
For an un-deformed, un-hydrated soil matrix the capillary pressure is expressed as a function of effective (or normalized) aqueous phase saturation, as given below
\begin{align}\label{eqn:brooksCoreyRelationship}
 P_{c0} = P_{entry} \ S_{we}^{\ -1 / \lambda_{BC}}
\end{align}
where, $P_{entry}$ is the gas entry pressure, 
$\lambda_{BC}$ is the soil specific parameter depending on the pore-size distribution, and 
$S_{we}$ is the normalized aqueous phase saturation given by $S_{we} = \dfrac{S_w - \left(S_{wr}+S_{gr}\right)}{1 - S_h - \left(S_{wr}+S_{gr}\right)}$. 

\paragraph{}
The effect of presence of hydrate in the soil matrix and the changing hydrate saturation on the capillary pressure $P_c$ is modelled by scaling $P_{c0}$ 
with a scaling-factor $f^{Pc}_{S_h}$ which is a function of $S_h$ \cite{ClementHookerSkeen1996,Rockhold2001}.
Also, the effect of changing porosity due to deformation of the porous matrix is accounted for by scaling $P_{c0}$ with a scaling-factor $f^{Pc}_{\phi}$ 
which is a function of $\phi$ using Civan's power-law correlation \cite{Civan2000}.
Thus, the capillary pressure is given by 
\begin{align}
\label{eqn:Pc_Scaled}
& P_c = P_{c0}\cdot f^{Pc}_{S_h}\left(S_h\right) \cdot f^{Pc}_{\phi}\left(\phi\right)
\end{align}
where, 
\begin{align*}
 f^{Pc}_{S_h} = \left(1-S_h\right)^{-\frac{m \lambda_{BC} - 1}{m \lambda_{BC}}} \ \text{, and } \ f^{Pc}_{\phi} = \frac{\phi_0}{\phi} \left( \frac{1-\phi}{1-\phi_0}\right)^{a} .
\end{align*}
where $m$ and $a$ are model parameters.

\paragraph{\textbf{Intrinsic permeability} } 
The intrinsic permeability, $K$, is related to the \textit{connectivity of the pore spaces} and the \textit{grain size} of the soil. 
It is a property of the soil matrix and is independent of the pore-fluids. 
An estimate of the intrinsic permeability can be made using mathematical expressions such as those proposed by Bear \cite{Bear1972}, or Mualem \cite{Mualem1976}.
Usually, however, the intrinsic permeability is evaluated experimentally as part of the characterization of the soil sample.

\paragraph{} The effect of changing hydrate saturation on the intrinsic permeability is modelled by scaling the initial or reference intrinsic permeability of the sediment $K_{0}$ 
with a scaling-factor $f^{K}_{S_h}$ which is a function of $S_h$ \cite{ClementHookerSkeen1996,Rockhold2001},
and, the effect of changing porosity due to deformation of the porous matrix is accounted for by scaling $K_{0}$ 
with a scaling-factor $f^{K}_{\phi}$ which is a function of $\phi$ using Civan's power-law correlation \cite{Civan2000}.
Thus, the intrinsic permeability for the hydrate sample is modelled as, 
\begin{align}
K &= K_{0}\cdot f^{K}_{S_h}\left(S_h\right) \cdot f^{K}_{\phi}\left(\phi\right) 
\end{align}
where,
\begin{align*}
 f^{K}_{S_h} = \left(1-S_h\right)^{\frac{5m + 4}{2m}} \quad \text{, and } \quad f^{K}_{\phi} = \frac{\phi}{\phi_0} \left( f^{Pc}_{\phi} \right)^{-2} \ .
\end{align*}

\paragraph{} These scaling factors (for both, $P_c$ and $\kappa$) are derived based on the assumption that hydrate phase sticks uniformly at the pore surface. For the ideal case of a spherical pore geometry, $m=3$.
 
\paragraph{\textbf{Relative permeabilities} } The relative permeability factors for both mobile phases are evaluated using the Brooks-Corey model 
in conjunction with the Burdine theorem \cite{Burdine1953}, as,
\begin{align}\label{eqn:BCRelativePermeabilities}
 &k_{rw} = S_{we}^{\frac{2+3\lambda_{BC}}{\lambda_{BC}}} \quad \text{, and } \notag \\ 
 &k_{rg} = \left( 1-S_{we} \right)^{2} \left( 1 - S_{we}^{\frac{2+\lambda_{BC}}{\lambda_{BC}}}\right) \ .
\end{align}

\paragraph{\textbf{Specific surface area} } An important property of the porous matrix is it's specific surface area, $A_{s}$, which is defined as the ratio of the total internal surface area of the pores enclosed within an REV 
to the total volume of the REV. The correlation proposed by Yousif \cite{Yousif1991} is used for estimating the specific surface area of the porous matrix,   
\begin{align}\label{eqn:specificSurfaceArea}
 & A_{s} = \sqrt{\frac{\phi_{eff}^3}{2\ K}} \quad
 \text{where, \quad } \phi_{eff} = \phi \left( 1-S_h \right) \ .
\end{align}

\paragraph{\textbf{Hydraulic tortuosity} } 
Tortuosity is empirically related to porosity, as,
\begin{align} \label{eqn:hydraulicTortuosity}
  & \tau = \phi^n \quad \text{where, } 1 \leq n \leq 3  \ .
\end{align}

\subsubsection{Poro-elasticity}
      \paragraph{\textbf{Principle of effective stress} } When a porous fluid-filled soil encounters an external load, the stress is partly supported by the soil matrix and partly by the pore-fluids. 
      The deformation of the porous medium is effected by only that part of the total stress that is supported by the soil matrix. 
      This stress, introduced by Terzaghi \cite{Terzaghi1925}, is called the effective stress. 
      Using this concept, the total stress $\tilde \sigma$ appearing in Eqn. \eqref{eqn:momentumBalance_soil-1} can be decomposed as,
      \begin{align} \label{biotStress}
      \tilde{\sigma} = \tilde{\sigma}^{\prime} + \alpha_{biot} P_{eff} \tilde{I}
      \end{align}
      where, $\tilde{\sigma}$ is the total stress acting on the bulk porous medium, $\tilde{\sigma}^{\prime}$ is the effective stress acting on the composite skeleton,
      and $P_{eff}$ is the effective pore-pressure exerted by the mobile phases, given by 
      \begin{align*}
       P_{eff} = \dfrac{S_w}{S_w + S_g}P_w + \dfrac{S_g}{S_w + S_g}P_g \ .
      \end{align*}
      
      \paragraph{} $\alpha_{biot}$ is Biot's parameter. 
      One of the generally accepted expressions for $\alpha_{biot}$ in rock-mechanics applications is $\quad \alpha_{biot} = 1 - \dfrac{B_{m}}{B_{sh}} \ $ \cite{BiotWillis1957}.
      Here, $B_{sh}$ is the bulk modulus of the composite matrix, and $B_m$ is the bulk modulus of the porous medium.

      \paragraph{\textbf{Stress-strain behaviour} }


      The deviatoric stress vs. axial strain response obtained from tri-axial tests of distributed hydrate-bearing sands reported by 
      Masui et al. \cite{Masui2006} and Miyazaki et al. \cite{Miyazaki2011} show that the hydrate-sand specimens are not elastic bodies,
      but it is possible to consider the stress-strain relationship to be elastic if the range of application is sufficiently limited to small-strain cases far away from the critical state. 
      
      \paragraph{} In our model, we use the linear-elastic constitutive law to describe the stress-strain response of the hydrate-soil composite matrix, given as,
      \begin{align}\label{eqn:linearElasticStressStrainLaw}
       &\tilde{\sigma}^{\prime} = 2 \ G_{sh} \ \tilde{\epsilon} + \lambda_{sh} (tr \ \tilde{\epsilon})\ \tilde{I}
      \end{align}
      where, $G_{sh}$ and $\lambda_{sh}$ are the Lame's parameters for the elastic composite-matrix 
      and can be defined in terms of the \textit{apparent} elastic mechanical properties (e.g. Young's modulus $E_{sh}$, and Poisson's ratio $\nu_{sh}$) as,
      \begin{align}\label{eqn:LamesParameters}
       G_{sh} = \frac{E_{sh}}{2\left(1 + \nu_{sh}\right)} \ , \quad \lambda_{sh} = \frac{E_{sh}\ \nu_{sh}}{\left(1 + \nu_{sh}\right)\left(1 - 2\ \nu_{sh}\right)} \ .
      \end{align}
      
      \paragraph{} $\tilde{\epsilon}$ is the linearized strain, given by $\quad \tilde{\epsilon} = \dfrac{1}{2}\left( \nabla {\bf u} + \nabla^T {\bf u}\right)$.

      \paragraph{Elastic properties } From the tri-axial tests it is observed that the presence of methane hydrate, in general, increases stiffness and leads to higher strength.
      Also, the effect of methane hydrate saturation on the Poisson ratio appears to be small. 
      Further details on the general trends of mechanical properties of methane hydrates can be found in Waite et al. \cite{Waite2009}, and Soga et al. \cite{Soga2006}.
      
      \paragraph{} To make the model consistent with these observations, we assume the Poisson ratio $\nu_{sh}$ to be a constant, 
      and we define the Young's modulous using the expression proposed by Santamarina and Ruppel \cite{SantamarinaRuppel2010}, given as,
      \begin{align}\label{eqn:SantamarinaRuppelEsh}
       E_{sh} = E_{s0} \left( \frac{\sigma_{c}}{\sigma_{c0}}\right)^b + c \ E_h \ \left(S_h\right)^d
      \end{align}
      where, $E_{s0}$ is the isothermal Young's modulus of hydrate-free sand at the reference confining stress of $\sigma_{c0} = 1$ kPa,
      $b$ is the sensitivity of the Young's modulus of hydrate-free sand to confining stress $\sigma_c$, 
      $c$ is the contribution of the isothermal Young's modulus of hydrate for a given pore habit, i.e., pore filling, cemented (grain-coating), or patchy,
      and $d$ is the nonlinear effect of hydrate saturation.


      \paragraph{\textbf{Compressibility} } 
      The grains of the composite material are assumed to be incompressible, but the bulk material as a whole is compressible. 
      This compressibility can be attributed to the fact that due to the pore pressure variations or isotropic loading the grains redistribute, which, on macro scale, manifests as change in density of the solid material.
      This change in density can be modelled as,
      \begin{align}\label{eqn:RateOfChangeOfSolidDensity}
      \DDt \ \rho_{sh} = \frac{\rho_{sh}}{G_{sh} \left(1-\phi_{eff}\right)} \left( \DDt \ \sigma - \phi\ \DDt \ P_{eff} \right)
      \end{align}
      where, $\sigma$ is the isotropic stress.

\section{Numerical solution} \label{sec:numericalSolution}
\begin{figure*}
  \centering
  \includegraphics[scale=0.48]{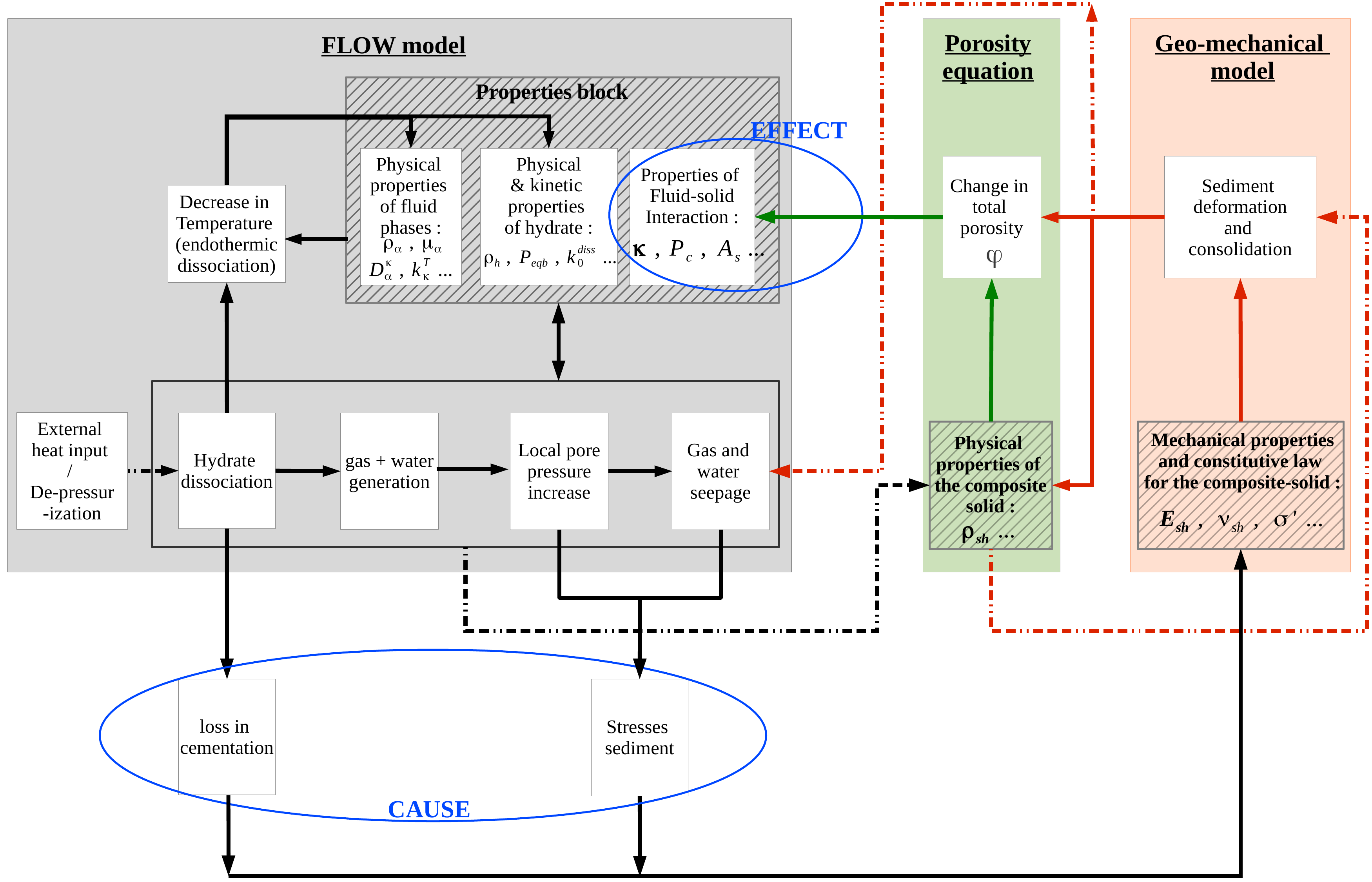}
 \caption{'Cause-effect' based interaction between the flow and geomechanical models. \\ 
	  The bold arrows represent the direct forward coupling, (e.g. $F_1-F_1$ or $F_1-F_2$, $F_2-F_3$, $F_3-F_1$). The dashed arrows represent indirect coupling (e.g. $F_1-F_3$) or backward coupling (e.g. $F_2-F_1$, $F_3-F_2$)}
  \label{fig:decouplingStrategy}
\end{figure*}

\paragraph{} 
The mathematical model describing the hydromechanical processes in a hydrate system actually contains two sub-classes of models, 
the \textit{flow and transport model} comprising of the mass, momentum, and energy balance equations for the phases occupying the pore spaces in the hydrate formation, 
i.e., equations \eqref{eqn:massBalance_CH4_H2O}, \eqref{eqn:massBalance_Hydrate}, \eqref{eqn:darcysLaw}, and \eqref{eqn:energyBalance},
and the \textit{geomechanical model}, comprising of momentum balance equation \eqref{eqn:momentumBalance_soil-1}. 
The soil phase mass balance, Eqn \eqref{eqn:massBalance_Soil}, can be seen as a \textit{glue} between these two sub-models (Fig. \ref{fig:decouplingStrategy}).

\paragraph{}
The interaction between these two models manifests physically in the form of, 
a) changes in the hydraulic properties (total porosity and permeability) due to deformation of the solid matrix structure, 
b) shift in the seepage velocity of the pore fluids due to the rate of deformation of the solid matrix, and,
c) changes in the mechanical properties of the solid matrix (strength, bulk modulus, density, etc.) resulting from the flow dynamics of the pore-fluids and the loss in cementation due to melting of the hydrate phase.
In other words, each model affetcs the other model by altering it's \textit{properties}. Thus, the nature of the coupling between these two models is dynamic, but weak.

\paragraph{} We use this observation to our advantage to devise a decoupled iterative solution strategy.
Broadly speaking, we first decouple the flow model and the geomechanical model, solve them separately for a given time-step, 
and then iteratively re-introduce the coupling through a block Gauss-Seidel solution scheme.

\paragraph{} We chose the following set of variables as the primary variables: 
the gas phase pressure $P_g$, the aqueous phase saturation $S_w$, the hydrate phase saturation $S_h$, the temperature $T$, the total-porosity $\phi$, and, the composite-matrix displacement $\mathbf{u}$.
All other variables can be derived (explicitly or implicitly) from this set of variables.

\paragraph{} The system of PDEs comprising the flow-system are solved fully implicitly for the variables $P_g$, $S_w$, $S_h$, and $T$. 
The spatial discretization is done using the fully up-winded classical cell centered finite volume method. 
Orthogonal grids aligned with the principal axes are defined and a control-volume formulation with two-point flux approximation (TPFA) is used.
For time-stepping an implicit Euler method is used. 

\paragraph{} The geomechanical system is solved for the primary variable $\mathbf{u}$. 
The soil momentum balance equation comprising the geomechanical model is spatially discretized using the Galerkin finite element (FEM) scheme defined on $Q1$ elements.
The FEM formulation described by Lewis and Schrefler \cite{LewisSchrefler2000} is used.

\paragraph{} The soil-phase mass-balance equation is solved separately for $\phi$. 
It is spatially discretized using the cell centered finite-volume method, and is marched forward in time using the implicit Euler method.

\paragraph{} 
The discretized model, can be represented as a system of algebraic equations as

\begin{align}
 &\mathbf{F}_1 : \mathbf{A}_1 \left( \mathbf{X}^{n+1}, \mathbf{X}^{n}\right) \mathbf{X}_1^{n+1} \ - \ \mathbf{B}_1 \left( \mathbf{X}^{n+1}, \mathbf{X}^{n}\right) \ = \ 0 \notag \\ & \text{which comes from the flow-model, }\notag\\
 &\mathbf{F}_2 : \mathbf{A}_2 \left( \mathbf{X}^{n+1}, \mathbf{X}^{n}\right) \mathbf{X}_2^{n+1} \ - \ \mathbf{B}_2 \left( \mathbf{X}^{n+1}, \mathbf{X}^{n}\right) \ = \ 0 \notag \\ & \text{which comes from the geomechanical-model, and }\notag\\
 &\mathbf{F}_3 : \mathbf{A}_3 \left( \mathbf{X}^{n+1}, \mathbf{X}^{n}\right) \mathbf{X}_3^{n+1} \ - \ \mathbf{B}_3 \left( \mathbf{X}^{n+1}, \mathbf{X}^{n}\right) \ = \ 0 \notag \\ & \text{which comes from the total-porosity equation. }\notag
\end{align}

\begin{figure}
  \centering
    \includegraphics[scale=0.27]{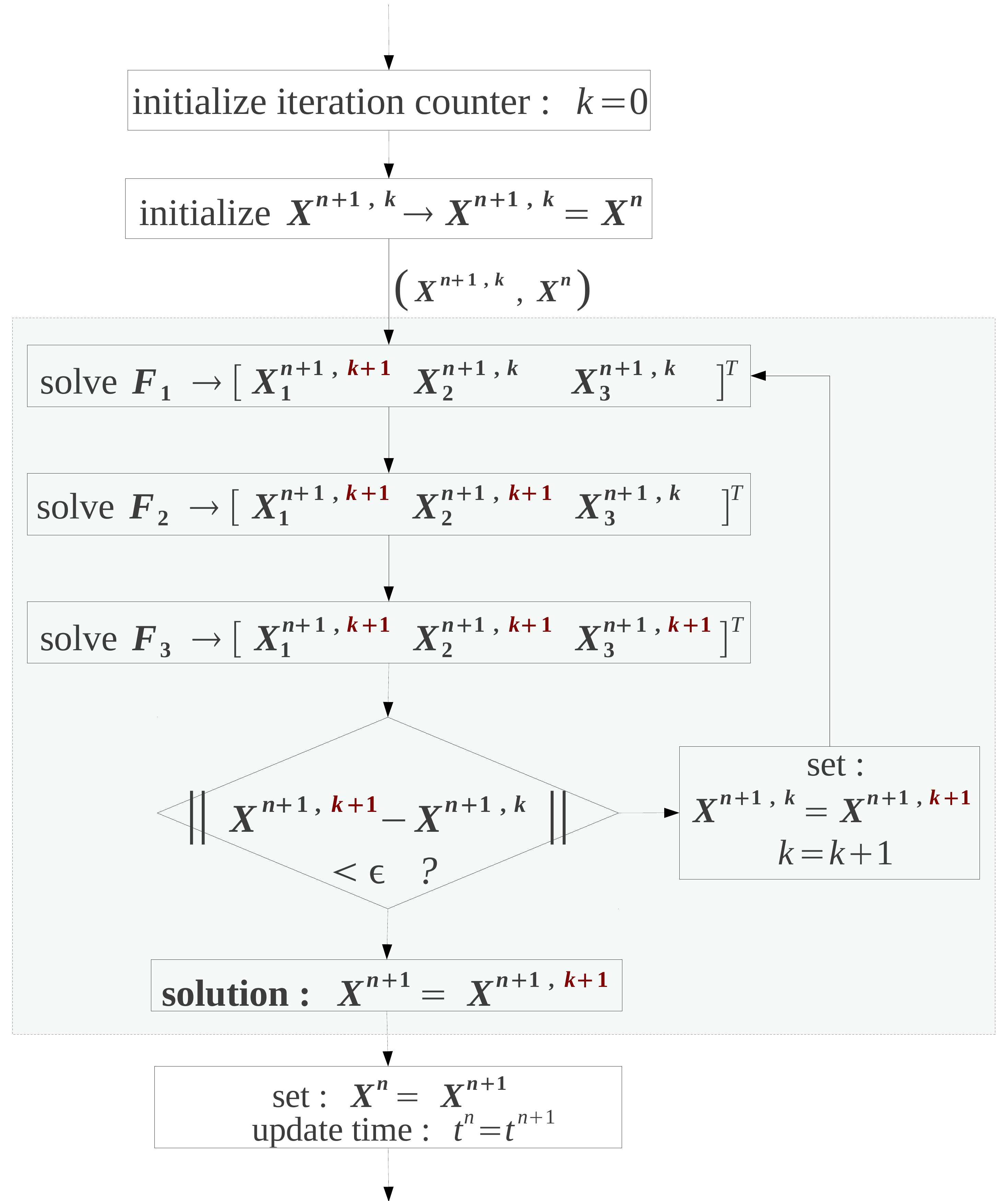}
  \caption{Block Gauss-Seidel \textit{outer} iterative solution loop}
  \label{fig:gaussSeidelIterationScheme}
\end{figure}

$\mathbf{X}$ is the solution vector given as $\mathbf{X} = \left[ \mathbf{X}_1 \ \ \mathbf{X}_2 \ \ \mathbf{X}_3 \right]^T$
where, $\mathbf{X}_1 = \left[ P_g \ \ S_w \ \ S_h \ \ T \right]^T$, $\mathbf{X}_2 = \left[ \mathbf{u} \right]^T$, and $\mathbf{X}_3 = \left[ \phi \right]^T $.
The indices $n$ and $n+1$ denote the solution at time $t^{n}$ and $t^{n+1}$ respectively. 
The strong non-linearities in each of the sub-systems $\mathbf{F}_1$, $\mathbf{F}_2$, and $\mathbf{F}_3$ are dealt with using a damped Newton-Raphson linearization. 
Each of the resulting linear sub-systems are solved using the SuperLU linear solver \cite{SuperLU}. 
This forms the \textit{inner iterative loop} which takes care of the decoupled solution.
The coupling between $\mathbf{F}_1$, $\mathbf{F}_2$, and $\mathbf{F}_3$ is re-introduced through an \textit{outer iterative loop} using a Gauss-Seidel scheme, 
as shown in Fig. \ref{fig:gaussSeidelIterationScheme}.

\paragraph{}
The numerical scheme is implemented in the C++ based DUNE-PDELab framework \cite{DUNEBastian2007}.
The numerical code is flexible in it's dimensionality and is capable of solving problems in 1D, 2D and 3D.
Furthermore, the \textit{block} structure of the decoupling scheme makes the code modular,
and gives the flexibility to solve only the flow and transport model or the full geo-mechanical model, depending on the problem at hand.

\section{Numerical examples} \label{sec:numericalExample}
\paragraph{} 
The mathematical model described in Section \ref{sec:mathematicalModel} consists of three important parts that are specific to the gas-production application. 
These are the methane hydrate \textit{dissociation kinetics}, the \textit{poroelastic coupling}, and the \textit{poroelastic-kinetics coupling}.
In this section, test problems which focus on each of these parts separately are presented.

\paragraph{}
In Section \ref{subsec:verificationReactionKinetics}, we simulate hydrate dissociation in a depressurized lab-scale hydrate sample. 
In this problem the geo-mechanical effects are negligible and reaction kinetics dominates over fluid flow. Thus this problem effectively isolates dissociation kinetics from the other processes.
Next, in Section \ref{subsec:verificationThreePhaseHydrate}, we consider a problem similar to the five-spot problem in a diagonal flow setting, containing a melting hydrate block in the domain.
The geo-mechanical effetcs are negligible. This problem focuses on the coupling between hydrate phase change and fluid flow in the three-phase hydrate model.
In Section \ref{subsec:verificationTerzaghiConsolidation}, we simulate the classical $1D$ consolidation problem to ensure a correct implementation of the poroelastic coupling in our numerical scheme.
Following this, as an extension to the $1D$ consolidation problem, in Section \ref{subsec:verificationHPE}
we present a $1D$ test setting where a hydrate sample is depressurized while being subjected to an external vertical loading. 
This problem focuses on the kinetics-poroelastic coupling.
Under simplifying assumptions, an analytical solution for the phase-pressure evolution is derived for this setting, 
which is then used to verify the numerical scheme.

\paragraph{}
In the final example in Section \ref{subsec:3Dexample} we simulate a $3D$ hydrate reservoir which is destabilized through depressurization using a low pressure gas well.
In this example we put together all the important model components including hydrate phase change, non-isothermal effects, multi-phase multi-component flow, and poroelastic soil deformation.

\subsection{\textbf{Test 1}: Dissociation kinetics model} \label{subsec:verificationReactionKinetics}
\paragraph{}
We consider 1D and 2D experiments on hydrate dissociation by depressurization by Tang et al. (2007) \cite{Tang2007} and Yuhu et al. (2009) \cite{Yuhu2009}, respectively.

\subsubsection{1D Case} \label{subsubsec:TANG-1D}
\paragraph{\textbf{Experimental set-up} }
A cylindrical, stainless steel cell with internal diameter of $38$ mm
and length of $500$ mm was used as the main pressure vessel. The cell
was jacketed with an insulating, impermeable layer and was immersed in
an air-bath. During each experimental run, the dry sands were sieved
into size range of $300-450$ $\mu m$ and were pushed tightly into the
vessel, resulting in a sediment with porosity of $33 \%$ and a
permeability of $300$ mdarcy.

\paragraph{}
The sediment was saturated with distilled water, and the methane hydrate
was formed in-situ by slowly injecting methane gas at a pressure higher
than the equilibrium pressure. The hydrate was formed in two stages to
obtain a homogeneous distribution.

\paragraph{}
To perform the dissociation experiment, the back pressure regulator was
set to a pressure lower than the hydrate equilibrium pressure at the
working temperature, and the outlet valve was opened quickly. The gas
released through the outlet valve was continously recorded and a \textit{cumulative gas-production curve} was plotted.

\begin{figure} 
  \begin{center}
  \includegraphics[scale=0.775]{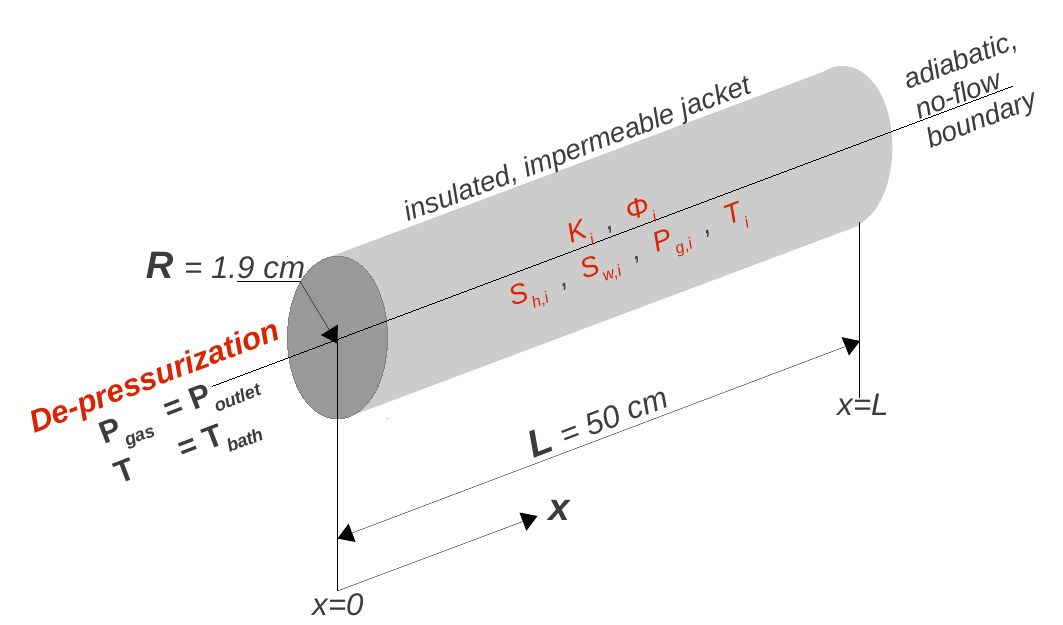}
  \caption{Test setting for 1D experiment (Tang et al. (2007))} \label{fig:TangSchematic}
  \end{center}
\end{figure}

\paragraph{\textbf{Numerical simulation} }
A schematic of the test domain for this experiment is shown in Fig. \ref{fig:TangSchematic}. 
The domain is discretized into $100$ cells along the X-axis and the simulation is performed in 1D. 
The geo-mechanical block is switched off and only the flow-transport block is solved.
Gravity is neglected.
The \textit{depressurization} (i.e. back pressure regulation) is considered at the left boundary. 
The left boundary also serves as the \textit{gas outlet}. 
Two depressurization modes have been considered for this test. 
In the first (\textit{Test-ID:Run2}), the pressure is decreased from $3.535$ MPa to $0.93$ MPa at $T_{bath} = 1.54^0C$. 
In the second (\textit{Test-ID:Run3}), the pressure is decreased from $3.584$ MPa to $1.94$ MPa at $T_{bath} = 2.08^0C$.
The total dissociation process for \textit{Run2} and \textit{Run3} was reported by Tang et al. to last $40$ and $110$ minutes respectively. 
So, the $t_{end}$ for the numerical simulations was chosen accordingly.
The initial and boundary conditions are listed in Table \ref{table:TangICs} and Table \ref{table:TangBCs}, respectively.

\paragraph{}
Tang et al. (2007) used TOUGH-Fx/Hydrate to simulate the experimental data. 
The value of the intrinsic rate constant $k_d^0$, as reported by Tang et al. was back calculated to fit the experimental data.
We have used the same method to calibrate our model and obtain the best value of $k_d^0$ for each depressurization mode (\textit{Run2} and \textit{Run3}).
The resulting values for kinetic-parameters $\Delta E_{a}/R$ and $k_d^0$ are listed and compared in Table \ref{table:TangKineticParameters}.
"\textit{reported}" refers to the values reported by Tang, and
"\textit{fitted}" refers to the values obtained from our calculations.
Also, since the parameterization for $P_{eqb}$ was not reported by Tang et al., 
we have used the standard relationship proposed by Kamath and Holder \cite{KamathHolder1987} for pure methane dissolved in distilled water.

\begin{figure}
  \centering
    \subfloat[Run 2]{%
    \includegraphics[scale=0.62]{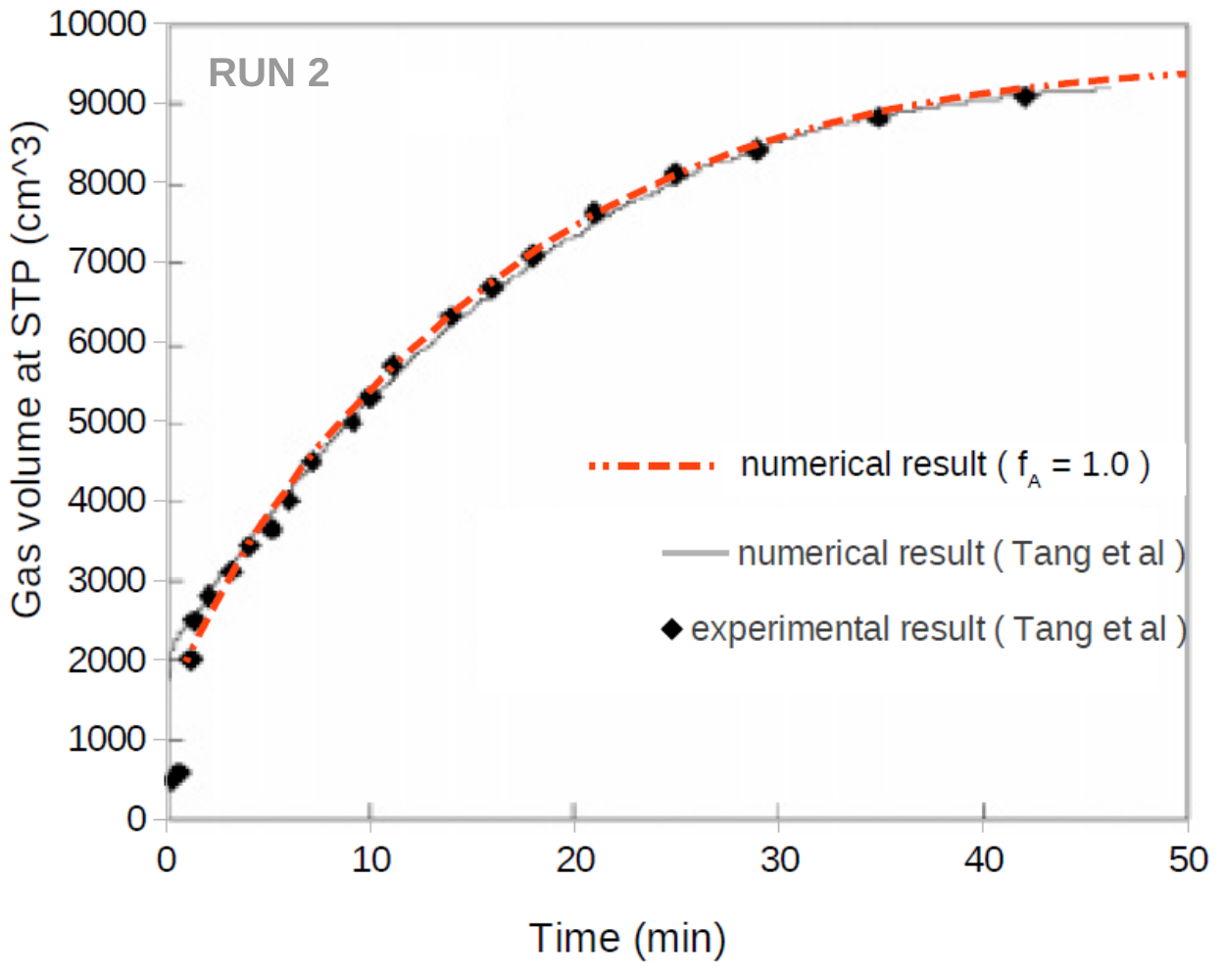}
    \label{fig:TangRun2}}
    \vfill
    \subfloat[Run 3]{%
    \includegraphics[scale=0.64]{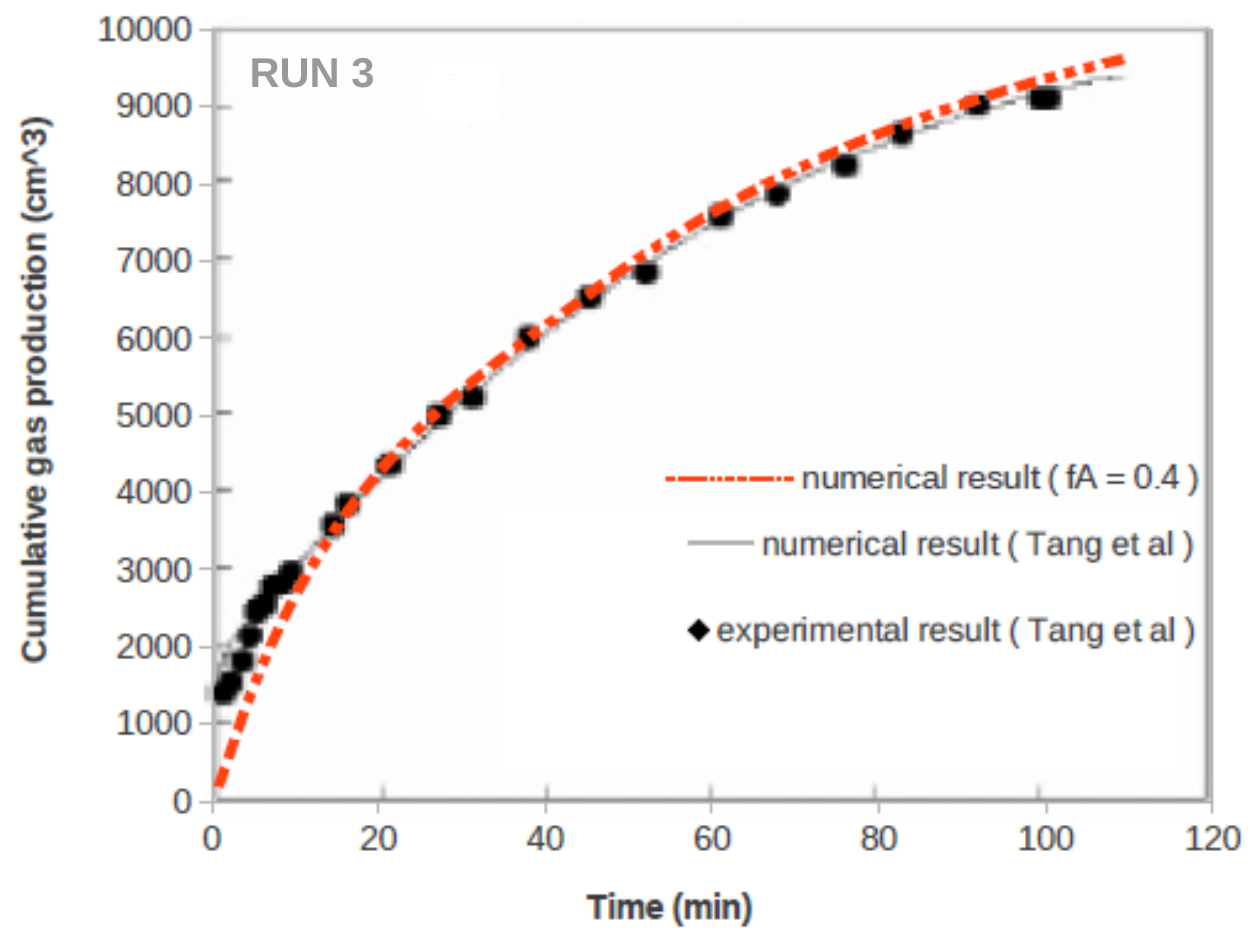}
    \label{fig:TangRun3}}
    \caption{Cumulative gas production curves (comparison with Tang et al. (2007)) Here, \ $f_A = \dfrac{\text{fitted } k^{0}_{d}}{\text{reported } k^{0}_{d}}$ .}
    \label{fig:TangResults}
\end{figure}

\paragraph{\textbf{Results} }
Fig. \ref{fig:TangRun2} and Fig. \ref{fig:TangRun3} show the comparison between cumulative gas volume curves obtained experimentally and numerically for \textit{Run2} and \textit{Run3}, respectively.
Our numerical results show a very close overall match with the experimental results. 


\begin{table}
     \begin{center}
      \caption{Initial conditions (Tang et al. (2007))}
      \label{table:TangICs}
      \begin{tabular}{L{0.7cm} L{1.5cm}| C{1.5cm} | C{1.5cm}}
	\hline
	ICs 		&		& Run 2		& Run 3 	\\ \hline
	$P_{g,i}$ 	& [ MPa ]	& $3.535$ 	& $3.584$  	\\
	$T_{i}$		& [ $^0C$ ]	& $1.54$	& $2.08$	\\
	$S_{w,i}$	& [ $\%$ ]	& $29.61 $	& $19.25 $ 	\\
	$S_{h,i}$	& [ $\%$ ]	& $21.83 $	& $25.44 $	\\
	$K_{i}$		& [ $mD$ ]	& $300$		& $300$		\\
	$\phi_{i}$	& [ $\%$ ]	& $30.8 $	& $30.8 $	\\ [.2ex]\hline
      \end{tabular}
    \end{center}
\end{table}

\begin{table}
     \begin{center}
      \caption{Boundary conditions (Tang et al. (2007))}
      \label{table:TangBCs}
      \begin{tabular}{L{2.5cm} L{0.1cm}| C{1.25cm} | C{1.25cm}}
	\hline
	at $x = 0$, $ t>0$		&& Run 2	& Run 3		\\ \hline
	$P_{outlet}$ \quad [MPa]	&& $0.93$	& $1.94$ 	\\
	$T_{bath}$ \qquad [$^0$C]	&& $1.54$	& $2.08$	\\ \hline
	at $x = L$, $t>0$ 	\\ \hline
	$\dot m_g = 0$ 		\\
	$\dot m_w = 0$ 		\\
	$\nabla T = 0$		\\ \hline
      \end{tabular}
    \end{center}
\end{table}

\begin{table}
  \begin{center}
  \caption{Kinetic parameters (Tang et al. (2007))}
  \label{table:TangKineticParameters}
  \begin{tabular}{L{2.5cm} | C{1.5cm} | C{1.5cm}}
  \hline
			& $\Delta E_{a}/R$ 	& $k_{d}^{0}$			\\ [.5em]
			& $\left[K\right]$	& $\left[ \frac{\text{mol}}{m^2\cdot Pa \cdot s}\right]$ \\ [.5em] \hline 
   Run 2 (reported)	& 9400			& $1.7 \times 10^4$		\\ [0.2em] 
   Run 2 (fitted)	& 9400			& $1.7 \times 10^4$		\\ [0.2em] \hline
   Run 3 (reported)	& 9400			& $1.4 \times 10^4$		\\ [0.2em]
   Run 3 (fitted)	& 9400			& $0.8 \times 10^4$		\\ [0.2em] \hline
  \end{tabular}
 \end{center}
\end{table}

\subsubsection{2D Case}

\paragraph{} 
The experimental set-up for hydrate formation and dissociation processes used in this experiment are very similar to the 1D experiment by Tang et al. described above. 
The main difference is the \textit{sample geometry}, which is cylindrical in the 1D case and square wafer-like in this case. 
Reaction kinetics is essentially only a time dependent process, and the number of spatial dimensions do not directly affect the kinetics. 
However, testing the kinetics model in both 1D and 2D geometries ensures that the spatial coupling between the different model components are correctly resolved, 
and that no spurious spatial effects manifest in the simulation of dissociation process.

\paragraph{\textbf{Experimental set-up} }
The hydrate formation and dissociation unit was a stainless steel vessel with length, width, and thickness of $380$ mm, $380$ mm, and $18$ mm respectively, and was immersed in an air-bath. 
The procedure for sand sample preparation and in-situ hydrate formation were similar to that described in Sec. \ref{subsubsec:TANG-1D}. 
The resulting sediment porosity and permeability were $40 \%$ and $1.97$ Darcy, respectively, and the hydrate saturation was $17.6\%$.

\paragraph{}
To perform the dissociation experiment, the back pressure was reduced from an initial pressure of $3.24$ MPa to $2.25$ MPa, and the outlet valve was opened quickly.
The bath temperature was maintained at $1.7^0C$.
The gas released through the outlet valve was continously recorded and a \textit{gas production rate} was plotted.

\begin{figure} 
  \begin{center}
  \includegraphics[scale=0.6]{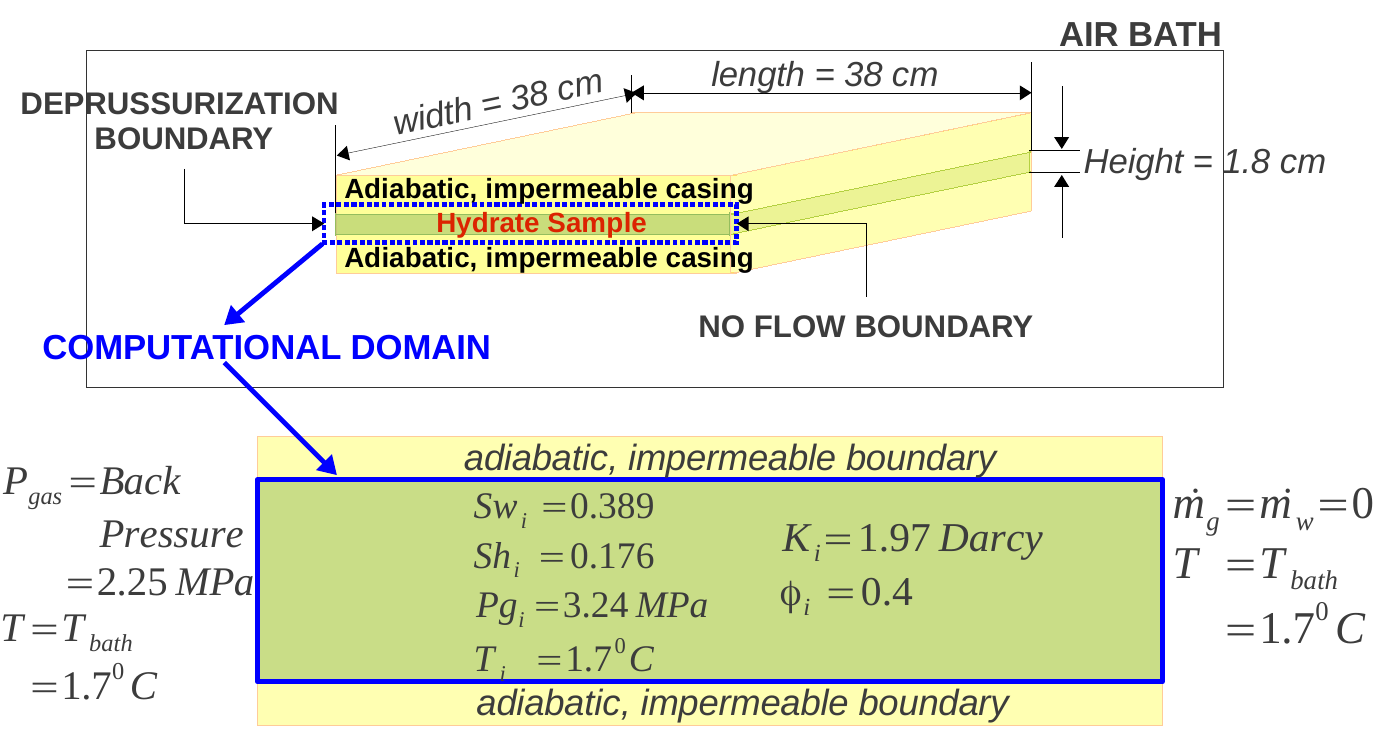}
  \caption{Test setting for 2D experiment (Yuhu et al.(2009))} \label{fig:YuhuSchematic}
  \end{center}
\end{figure}

\paragraph{\textbf{Numerical simulation} }
A schematic of the test domain for this experimental set up along with the initial and boundary conditions are shown in Fig. \ref{fig:YuhuSchematic}. 
The domain is discretized in $10 \times 100$ cells, and the simulations are performed in 2D.
The geo-mechanical block is switched off and only the flow-transport block is solved.
Gravity is neglected.
Depressurization and gas outlet are prescribed at the left boundary. 
The reaction-kinetics parameters are listed in Table \ref{table:YuhuKineticParameters}.

\begin{figure}
  \centering
    \subfloat[Gas pressure profile]{%
    \includegraphics[scale=0.445]{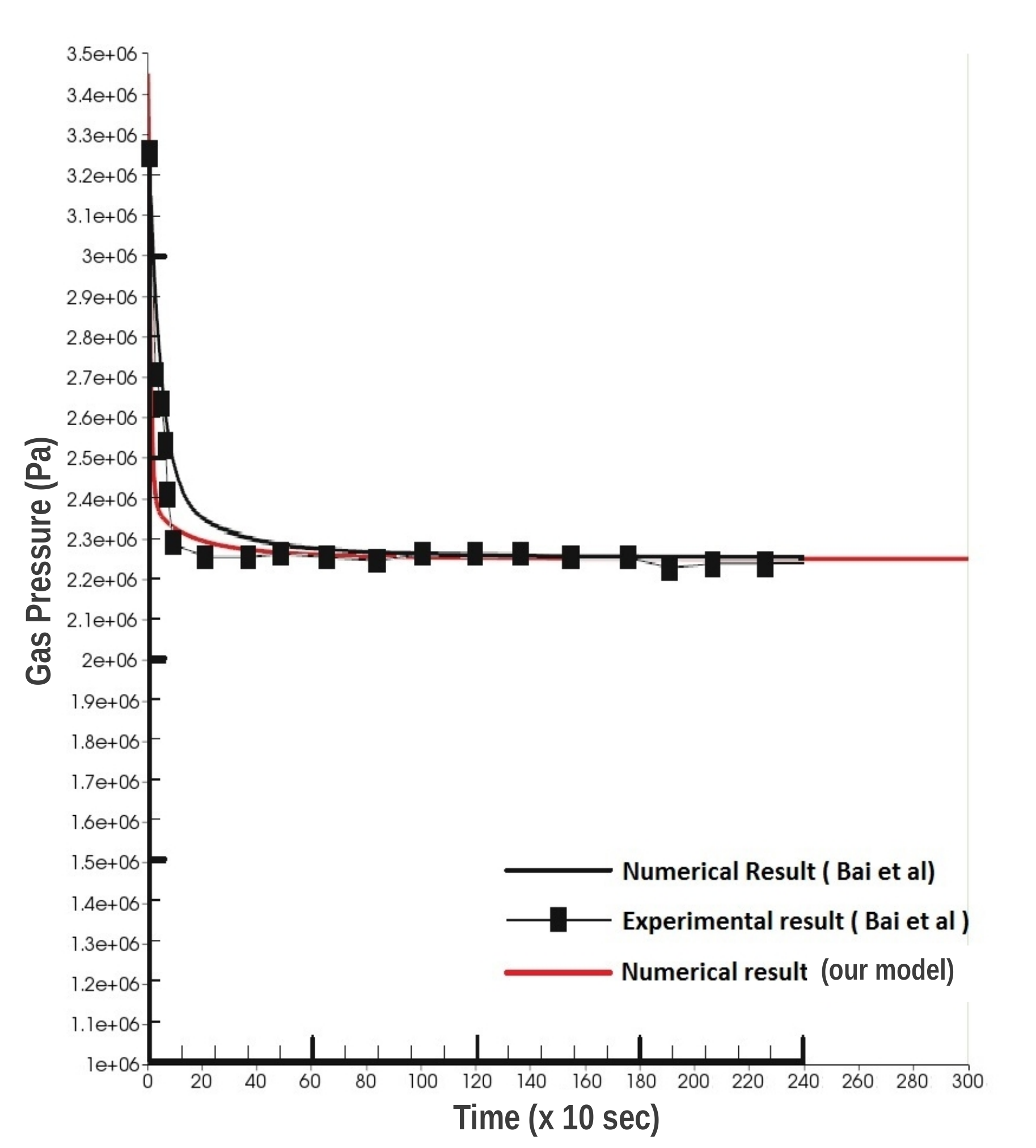}
    \label{fig:YuhuGasPressure}}
    \vfill
    \subfloat[Gas generation rate]{%
    \includegraphics[scale=0.43]{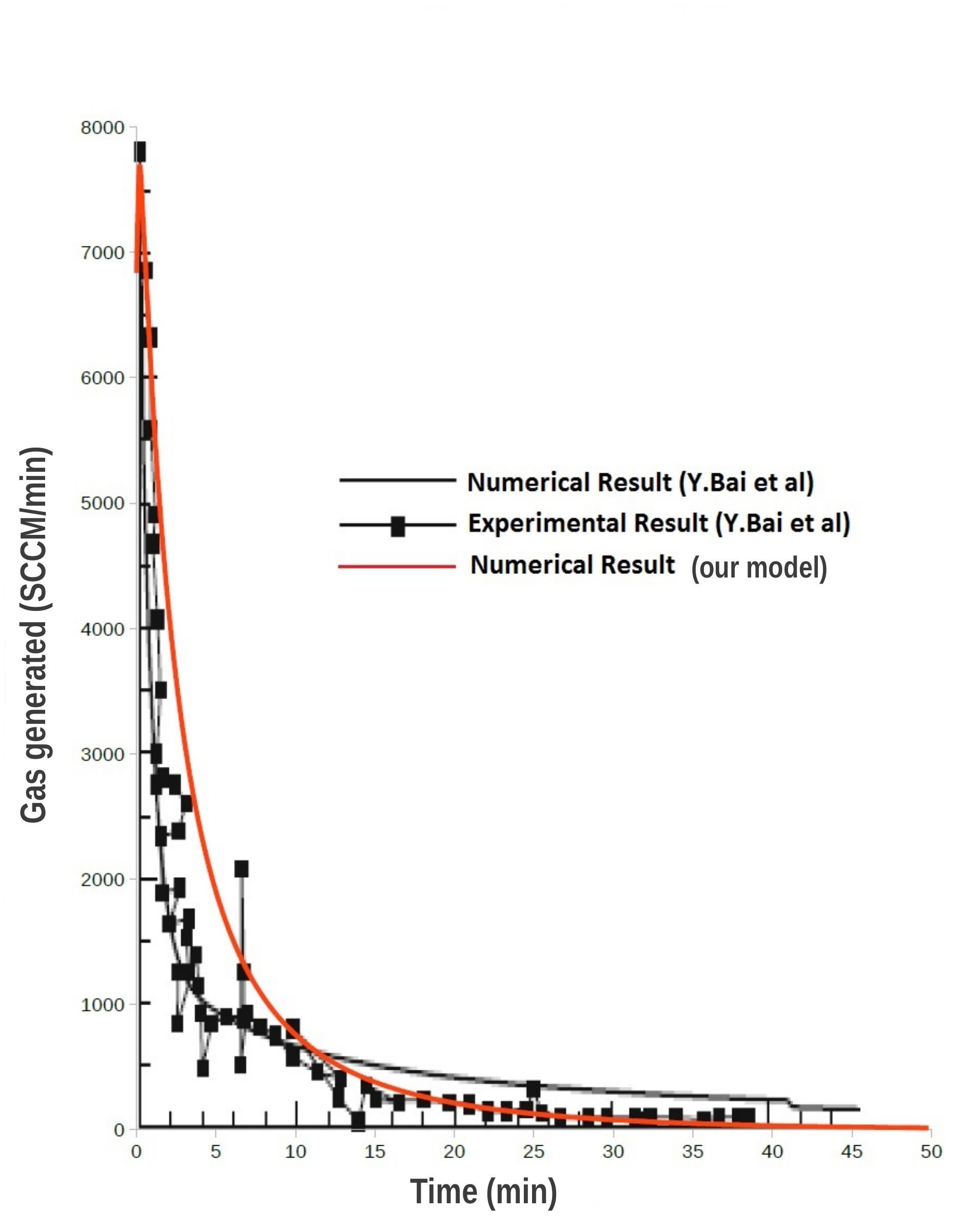}
    \label{fig:YuhuGasRate}}
    \caption{Results (comparison with Yuhu et al. (2009))}
    \label{fig:YuhuResults}
\end{figure}

\paragraph{\textbf{Results} }
Fig. \ref{fig:YuhuGasRate} and Fig. \ref{fig:YuhuGasPressure} show the experimental and numerical comparisons of the gas generation rate and the gas pressure for this setting.
Our numerical results show good agreement with the experimental values, especially towards steady state. 

\begin{table}
 \begin{center}
  \caption{Kinetic parameters (Yuhu et al. (2009))}
  \label{table:YuhuKineticParameters}
  \begin{tabular}{L{1cm} L{1.5cm}| C{4cm} }
  \hline
  $\Delta E_{a}/R$ 	& $\left[K\right]$						& $9752.73$			\\ 
  $k_{d}^{0}$		& $\left[ \frac{\text{mol}}{m^2\cdot Pa \cdot s}\right]$	& $3.6 \times 10^4$		\\ [1ex]
  $\dot Q_h$		& $\left[ \frac{J}{kg}\right]$					& $-1050 \ T + 3527000$		\\ [1ex]\hline
  \end{tabular}
 \end{center}
\end{table}

\subsection{\textbf{Test 2}: Three-phase hydrate model} \label{subsec:verificationThreePhaseHydrate}
\paragraph{}
In this section, we test the coupling between the kinetics model and the two phase flow model in our numerical scheme.


\paragraph{}
For this, we use an artificial setting similar to the five-spot problem, with the addition of a gas source in the domain in the form of a dissociating block of hydrate. 
This test ensures a correct implementation of the convection, diffusion, and the reaction terms in the 2D numerical scheme. 
It also ensures that the numerical scheme does not produce any spurious grid-based effects.

\begin{figure}
  \begin{center}
  \includegraphics[scale=0.6]{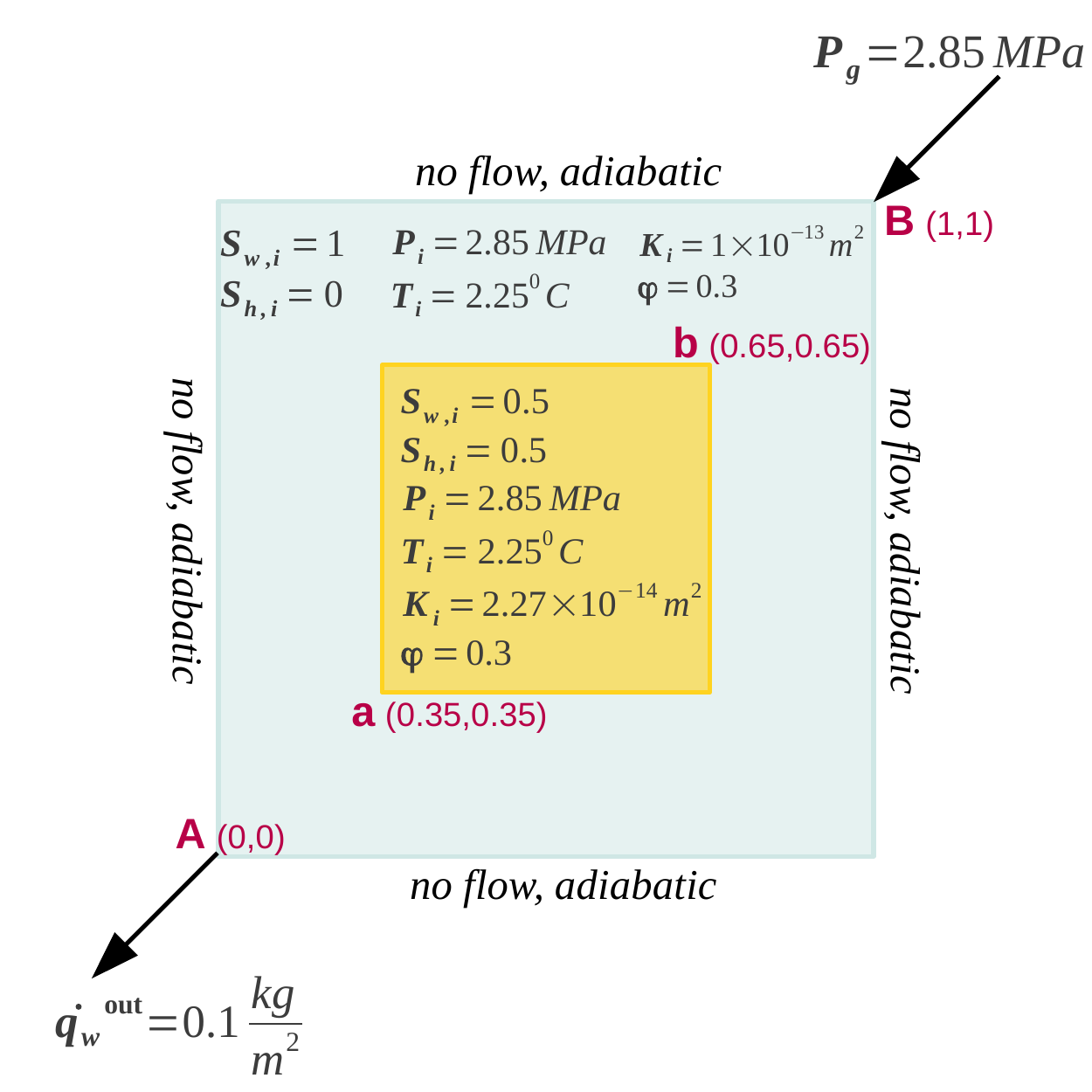}
  \caption{Test setting for 5-Spot test for hydrate model} \label{fig:fiveSpotHydrateSchematic}
  \end{center}
\end{figure}

\paragraph{\textbf{Problem set-up} } 
A schematic of the test domain is shown in Fig. \ref{fig:fiveSpotHydrateSchematic}.
The domain is a unit square with a $0.3$m $\times$ $0.3$m hydrate block located in the center. 
The domain is initially saturated with water. The hydrate saturation in the block is $50\%$.
Point \textit{A} at $(0,0)$ is the gas well. 
\textit{Neumann water-outflux} B.C. is prescribed at \textit{A}.
The diagonally opposite point \textit{B} at $(1,1)$ is held constant at initial pressure. 
The rest of the domain boundaries are closed and adiabatic.
The depressurization caused by water outflow at \textit{A} is expected to destabilize the hydrate block causing it to dissociate. 
The released gas must then get drawn towards the low pressure in the gas well at point \textit{A}.

\paragraph{\textbf{Numerical simulation} } 
For this problem, the geomechanical block is switched off and only the flow-transport block is solved.
The initial and the boundary conditions for the problem are specified in Fig. \ref{fig:fiveSpotHydrateSchematic}.
The hydraulic properties, hydrate stability curve parameters, and the dissociation kinetics parameters used in the numerical simulation are listed in Table \ref{table:parametersHydrateModel}.
The end-time for this problem is chosen as $t_{end}=500$ minutes.
For the base test (\textit{run0}), the domain is discretized uniformly into $20 \times 20$ cells. The time-step size is kept constant at $120$ seconds.
To check the mesh dependency of the numerical scheme, the mesh is successively halved, 
i.e. $\left(\Delta x\right)_{run1} = \frac{1}{2} \left(\Delta x\right)_{run0}$ , $\left(\Delta x\right)_{run2} = \frac{1}{2} \left(\Delta x\right)_{run1}$ , and $\left(\Delta x\right)_{run3} = \frac{1}{2} \left(\Delta x\right)_{run2}$ .
The time-step size is also successively halved so that $\Delta x / \Delta t $ ratio remains constant for each of the test-runs.

\begin{figure*}
 \centering
%
	  
	  \subfloat[\textbf{Run1} , $\left(\Delta x = \frac{1}{40}\text{m} \ , \ \Delta t = 60\text{s}\right)$][\textbf{Run1} \\ $\left(\Delta x = \frac{1}{40}\text{m} \ , \ \Delta t = 60\text{s}\right)$]{%
	  \includegraphics[trim=260mm 170mm 10mm 9mm,clip=true,scale=0.275]{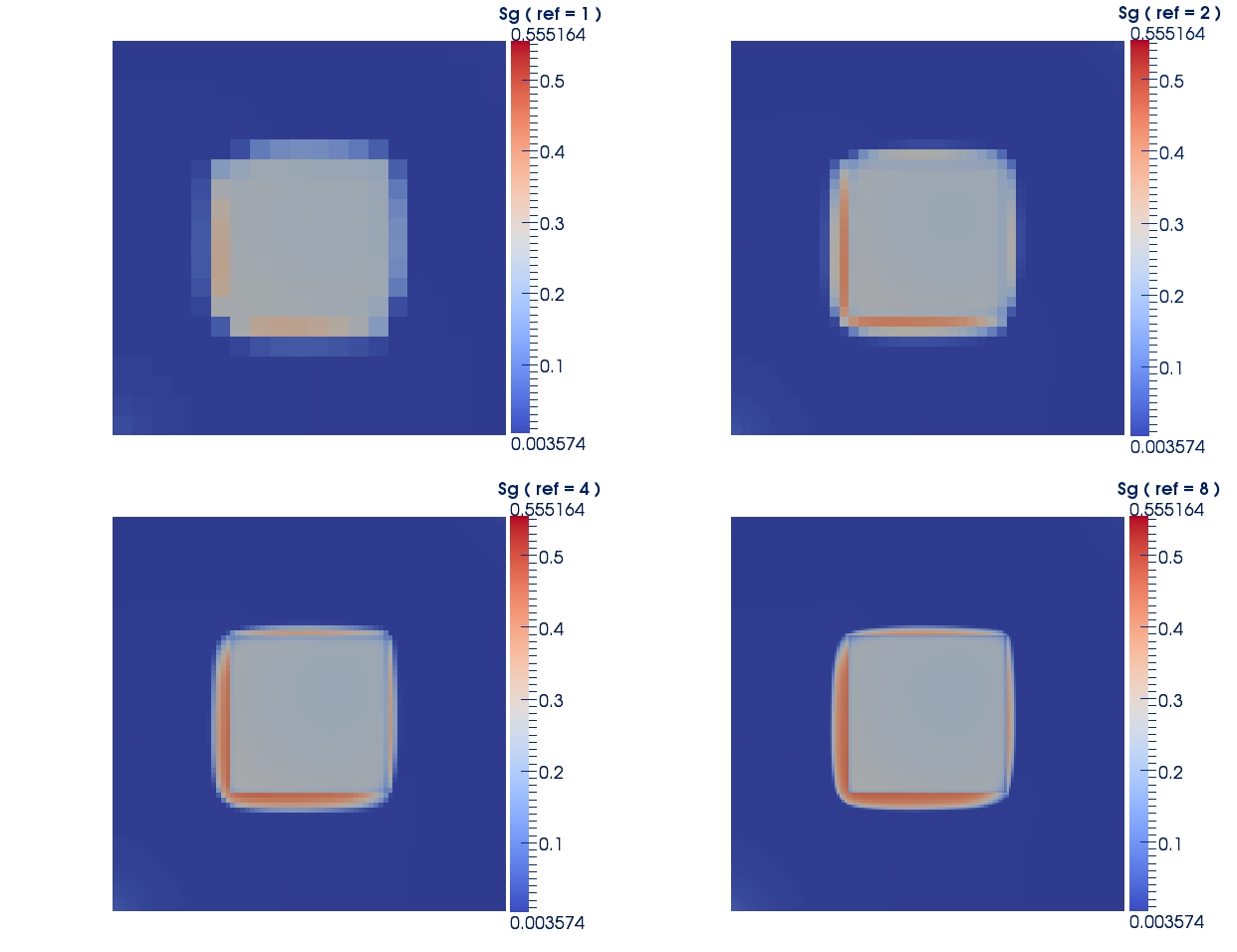}
	  \label{fig:Hydrate5Spot-Run1,t=100min}}
	  \hspace{0.7cm}
	  \subfloat[\textbf{Run1} , $\left(\Delta x = \frac{1}{40}\text{m} \ , \ \Delta t = 60\text{s}\right)$][\textbf{Run1} \\ $\left(\Delta x = \frac{1}{40}\text{m} \ , \ \Delta t = 60\text{s}\right)$]{%
	  \includegraphics[trim=260mm 170mm 10mm 9mm,clip=true,scale=0.275]{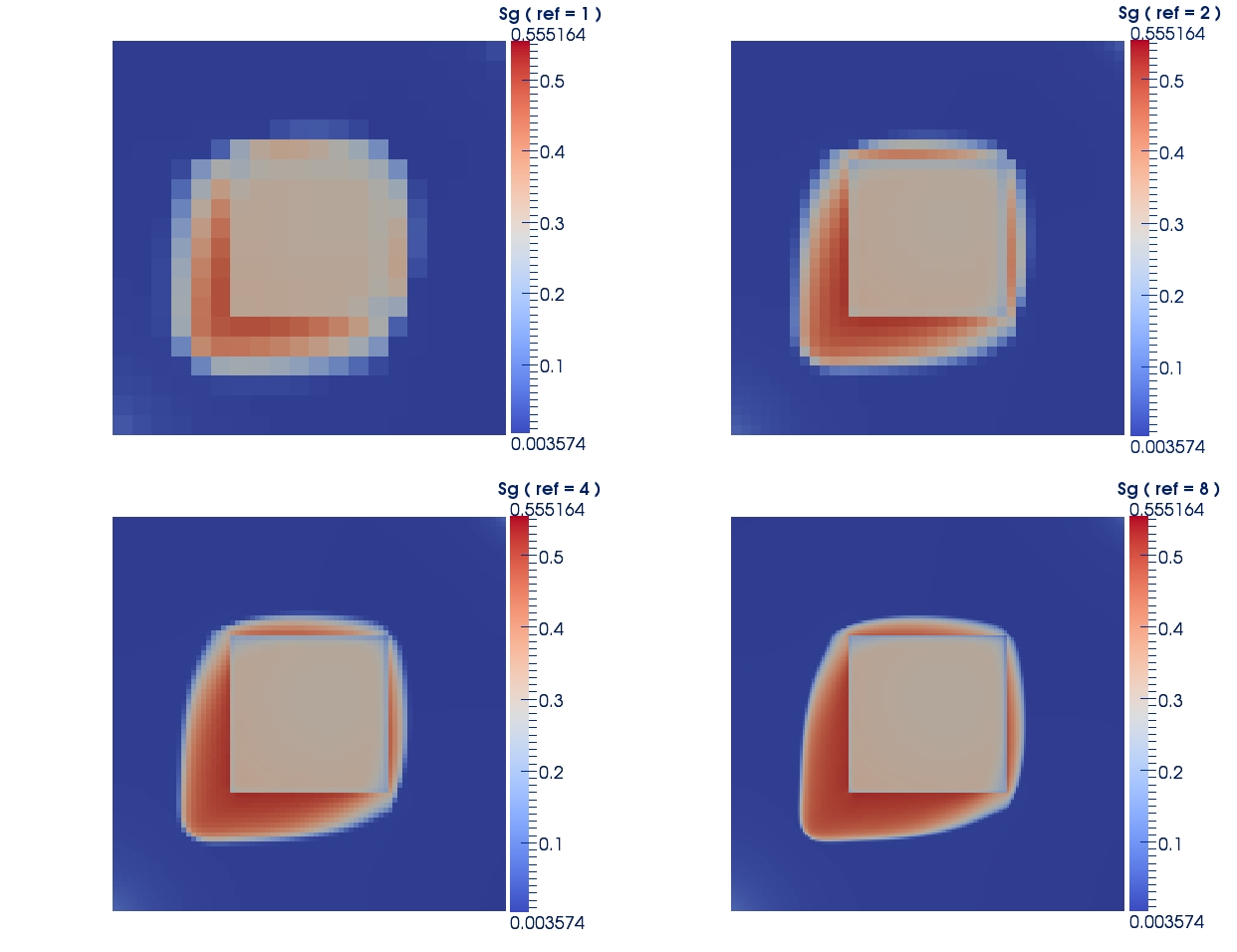}
	  \label{fig:Hydrate5Spot-Run1,t=200min}}
	  \hspace{0.7cm}
	  \subfloat[\textbf{Run1} , $\left(\Delta x = \frac{1}{40}\text{m} \ , \ \Delta t = 60\text{s}\right)$][\textbf{Run1} \\ $\left(\Delta x = \frac{1}{40}\text{m} \ , \ \Delta t = 60\text{s}\right)$]{%
	  \includegraphics[trim=260mm 170mm 10mm 9mm,clip=true,scale=0.275]{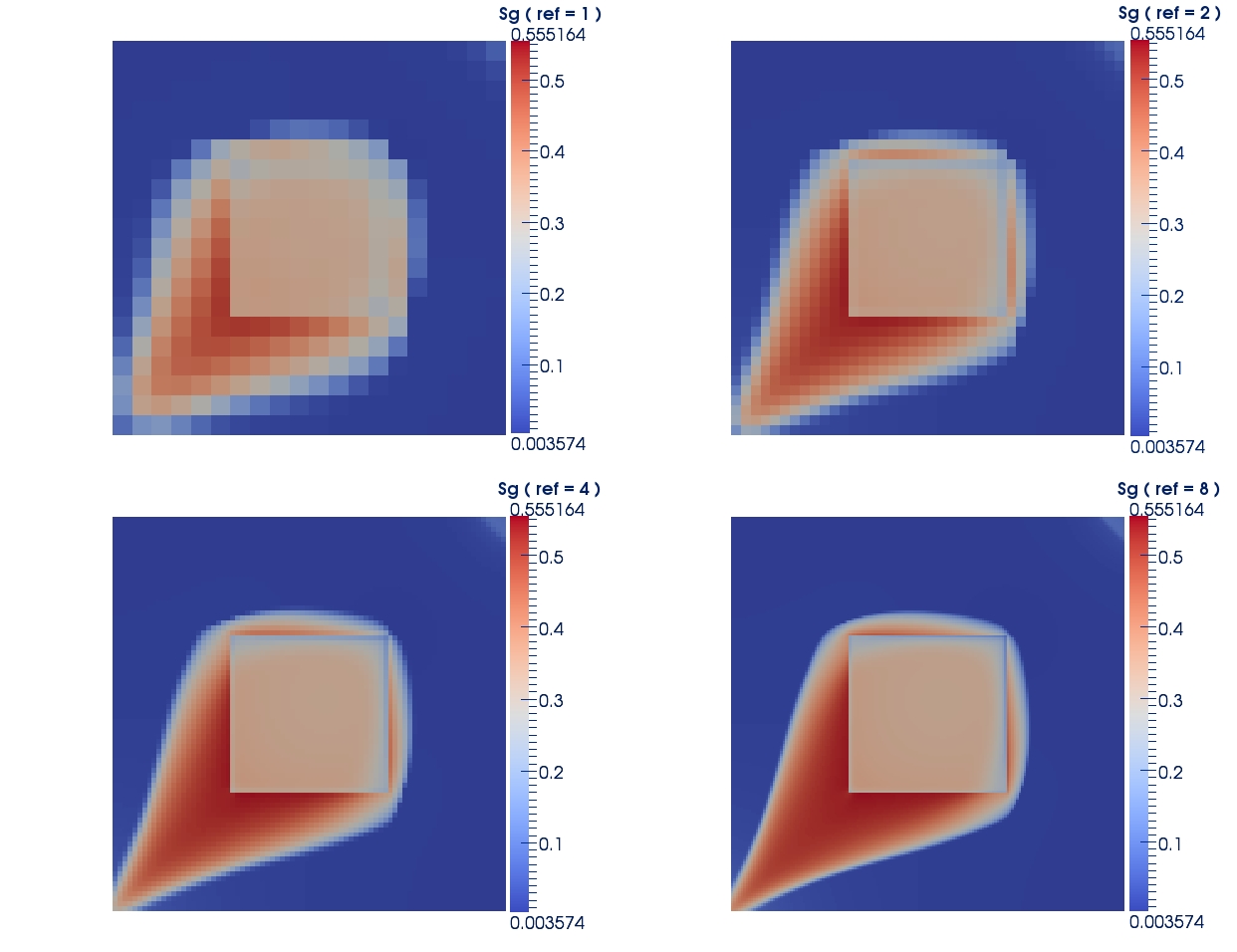}
	  \label{fig:Hydrate5Spot-Run1,t=300min}}
	  
	  \vfill
	  
%
	  
	  \subfloat[\textbf{Run3} , $\left(\Delta x = \frac{1}{160}\text{m} \ , \ \Delta t = 15\text{s}\right)$][\textbf{Run3} \\ $\left(\Delta x = \frac{1}{160}\text{m} \ , \ \Delta t = 15\text{s}\right)$]{%
	  \includegraphics[trim=260mm 8mm 10mm 177mm,clip=true,scale=0.275]{fiveSpotHydrate-screenshot1_time=100min.jpg}
	  \label{fig:Hydrate5Spot-Run3,t=100min}}
	  \hspace{0.7cm}
	  \subfloat[\textbf{Run3} , $\left(\Delta x = \frac{1}{160}\text{m} \ , \ \Delta t = 15\text{s}\right)$][\textbf{Run3} \\ $\left(\Delta x = \frac{1}{160}\text{m} \ , \ \Delta t = 15\text{s}\right)$]{%
	  \includegraphics[trim=260mm 8mm 10mm 177mm,clip=true,scale=0.275]{fiveSpotHydrate-screenshot1_time=200min.jpg}
	  \label{fig:Hydrate5Spot-Run3,t=200min}}
	  \hspace{0.7cm}
	  \subfloat[\textbf{Run3} , $\left(\Delta x = \frac{1}{160}\text{m} \ , \ \Delta t = 15\text{s}\right)$][\textbf{Run3} \\ $\left(\Delta x = \frac{1}{160}\text{m} \ , \ \Delta t = 15\text{s}\right)$]{%
	  \includegraphics[trim=260mm 8mm 10mm 177mm,clip=true,scale=0.275]{fiveSpotHydrate-screenshot1_time=300min.jpg}
	  \label{fig:Hydrate5Spot-Run3,t=300min}}
	  
	  \vspace{0.3cm}
	  
	  \makebox[0.2\linewidth]{\centering \textbf{Time} $\mathbf{= 100}$ \textbf{min}}
	  \hspace{1.8cm}
	  \makebox[0.2\linewidth]{\centering \textbf{Time} $\mathbf{= 200}$ \textbf{min}}
	  \hspace{1.8cm}
	  \makebox[0.2\linewidth]{\centering \textbf{Time} $\mathbf{= 300}$ \textbf{min}}
	  
	  \caption{Five-spot test - $S_g$ profile}
	  \label{fig:hydrate5Spot}
	  
\end{figure*}

\begin{figure}[h]
  \centering
    \subfloat[Time = $200$ minutes]{%
    \includegraphics[scale=0.315]{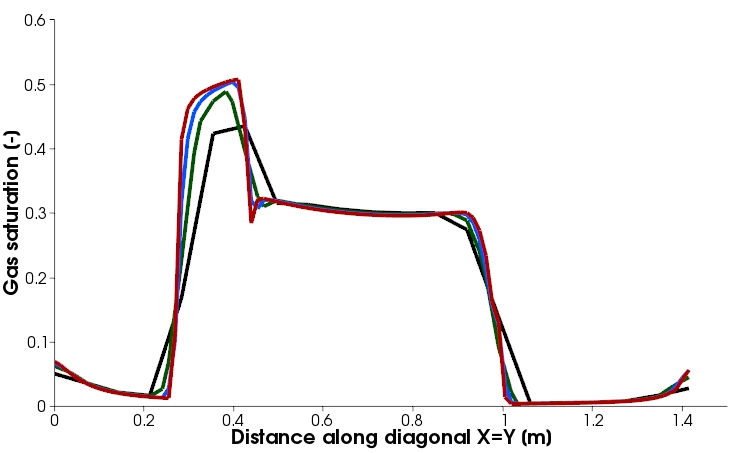}
    \label{fig:hydrate5Spot-linePlot-t=200min}}
    \quad
    \subfloat[Time = $300$ minutes]{%
    \includegraphics[scale=0.315]{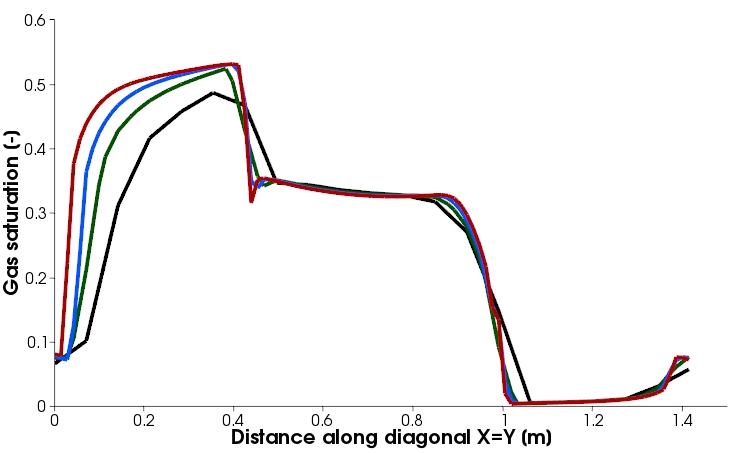}
    \label{fig:hydrate5Spot-linePlot-t=290min}}
    
        \vspace{0.4cm}

    \makebox[\linewidth]{\centering\textbf{Legend:} }\vfill
    \makebox[\linewidth]{ \raisebox{.4ex}{ \makebox{ \textcolor{black}{ \linethickness{0.3mm}\line(1,0){15} } } } ref=1, 
			  \raisebox{.4ex}{ \makebox{ \textcolor{darkolivegreen}{ \linethickness{0.3mm}\line(1,0){15} } } } ref=2,
			  \raisebox{.4ex}{ \makebox{ \textcolor{blue}{ \linethickness{0.3mm}\line(1,0){15} } } } ref=4,
			  \raisebox{.4ex}{ \makebox{ \textcolor{red}{ \linethickness{0.3mm}\line(1,0){15} } } } ref=8}
						     
    \caption{Five-spot test for hydrate model\\ $S_g$ profile along diagonal line $X=Y$}
    \label{fig:hydrate5Spot-linePlots}
\end{figure}

\paragraph{\textbf{Results} }
The gas plume takes about $300$ minutes to reach the gas well at \textit{A}.
Fig. \ref{fig:hydrate5Spot} shows the screenshots of gas saturation in the domain at $100$, $200$, and $300$ minutes.
Fig. \ref{fig:hydrate5Spot-linePlots} shows the line-plot of $S_g$ at $200$ and $300$ minutes along the diagonal aligned with flow direction, i.e. line $X-Y=0$.

\paragraph{}
In Fig. \ref{fig:hydrate5Spot-linePlot-t=200min} the solution shows convergence with mesh-refinement. 
The flow in the right half of the domain (i.e. $ Y+X-1 > 0 $) is diffusion dominated, whereas, 
that in left half (i.e. $ Y+X-1 < 0 $) is convection dominated being strongly influenced by the low pressure in the gas well. 
The gas front is more diffusive on a coarse mesh, but gets sharper as refinement is increased. 
Fig. \ref{fig:hydrate5Spot-linePlot-t=290min} shows the saturation of the gas plume that reaches the gas well at $300$ minutes. 

\begin{table}
 \begin{center}
  \caption{Five-spot test for hydrate model\\ \textbf{Model parameters} }
  \label{table:parametersHydrateModel}
      \begin{tabular}{L{1.cm} L{0.75cm} C{5.4cm} }
	\hline
	\multicolumn{3}{c}{Brooks-Corey parameters}			\\ \hline
	\\[-1.5ex]
	$P_{entry}$ 		& [Pa]		& $5000$ 		\\
	$\lambda_{BC}$ 		& [$^0$C]	& $1.5$		\\[2ex] \hline
	\multicolumn{3}{c}{Hydrate stability curve}			\\ \hline
	\\[-1.5ex]
	$P_{eqb}$		& [kPa]		& $\exp\left( 38.980 - \frac{8533.80}{T\left[\text{K}\right]}\right) \ \ \text{if } T > 273.15 $ 		\\
				& 		& $\exp\left( 14.717 - \frac{1886.79}{T\left[\text{K}\right]}\right) \ \ \text{if } T < 273.15 $ 		\\ [2ex]\hline
	\multicolumn{3}{c}{Dissociation kinetics parameters}		\\ \hline
	\\[-1.5ex]
	$\Delta E_{a}/R$ 	& $\left[K\right]$						& 9400			\\ 
	$k_{d}^{0}$		& $\left[ \frac{\text{mol}}{m^2\cdot Pa \cdot s}\right]$	& $3.6 \times 10^4$	\\ [2ex] \hline
      \end{tabular}
 \end{center}
\end{table}

\begin{figure}
 \centering
  \includegraphics[scale=0.2]{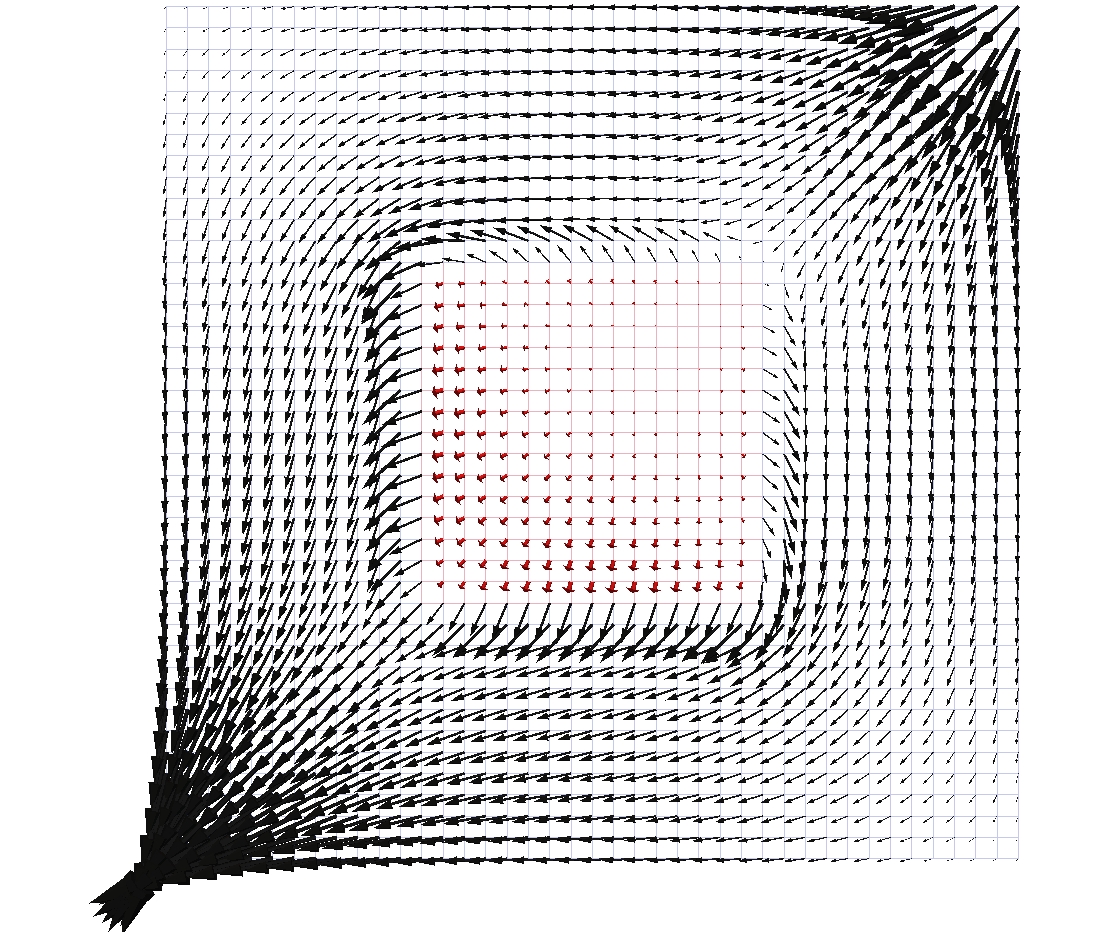}
  \caption{Five-spot test - $\mathbf{v}_g$ profile (at $t=30$ minutes)\\ Red vectors represent $\mathbf{v}_g$ in the hydrate zone \\ Black vectors represent $\mathbf{v}_g$ in hydrate free zone. }
  \label{fig:gasVelProfile_t=20min}
\end{figure}

\begin{figure}[h]
  \centering
    \subfloat[Time = $300$ minutes]{%
    \includegraphics[scale=0.315]{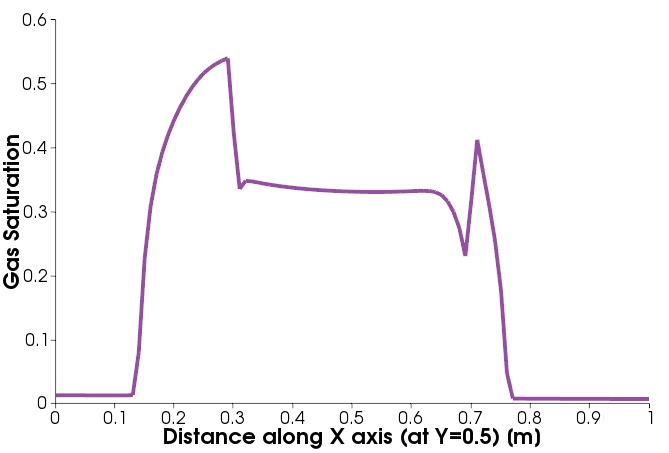}
    \label{fig:hydrate5Spot-XAxislinePlot-t=300min}}
    \quad
    \subfloat[Time = $500$ minutes]{%
    \includegraphics[scale=0.315]{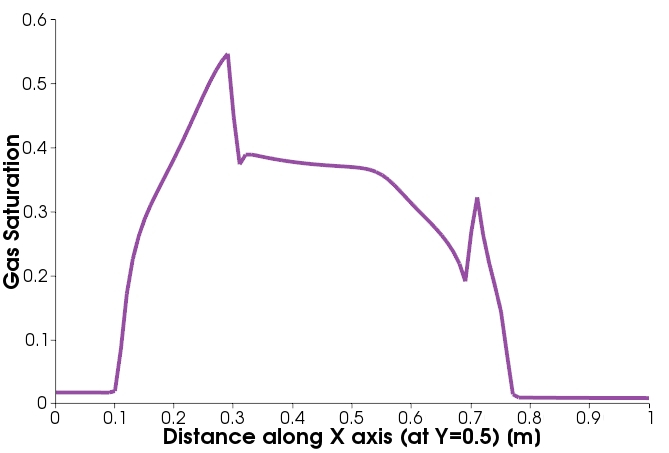}
    \label{fig:hydrate5Spot-XAxislinePlot-t=500min}}
    \caption{Five-spot test\\ $S_g$ profile along X-axis (at $Y=0.5$ m)}
    \label{fig:hydrate5Spot-linePlotsAlongXaxis}
\end{figure}

\paragraph{} We would like to point out that in Fig. \ref{fig:hydrate5Spot-linePlot-t=200min} and Fig. \ref{fig:hydrate5Spot-linePlot-t=290min}, 
what appears to be a \textit{kink} in gas saturation at the corner of the hydrate zone is \textit{not} a numerical artifact. 
This \textit{kink} is caused because of the following physical effects:
the gas velocity in the hydrate free zone  is higher than that in the hydrate zone (due to difference of almost an order of magnitude in the permeabilities) 
(See Fig. \ref{fig:gasVelProfile_t=20min}). 
This causes the gas to be sucked out of the hydrate zone faster than the time required by gas to equilibriate inside the hydrate zone. 
So, the gas begins to \textit{deplete} along the edges of the hydrate zone. 
Further, the extent of the depletion is higher where $\dot g^{CH_4}$  is lower (i.e., where $P_g$ is higher).
This effect is more clearly visible in Fig. \ref{fig:hydrate5Spot-XAxislinePlot-t=300min} and Fig. \ref{fig:hydrate5Spot-XAxislinePlot-t=500min} 
which show the $S_g$ profiles along X-axis (at $Y=0.5$ m) at times $t=300$ min and $t=500$ min respectively.
The pressure in the right half along the X-axis is higher than that in the left half, causing $\dot g^{CH_4}$ to be lower in the right half. 
Therefore, the extent of gas depletion is higher in the right half of the hydrate zone.

\subsection{\textbf{Test 3}: Poroelastic coupling} \label{subsec:verificationTerzaghiConsolidation}
\paragraph{} 
In this example we ignore the methane hydrates in the medium, thus reducing the model to a simple two-phase hydro-mechanical system. 
We consider the classical 1D consolidation problem by Terzaghi \cite{Terzaghi1925} to test the fluid pressure response generated by the mechanical compression of the soil.
This test was originally formulated by Terzaghi for analyzing the time delay observed when compressing clay layers 
and is now considered as a standard benchmark test for the coupling relationships between fluid and mechanical systems. 

\begin{figure}
 \centering
     \includegraphics[scale=0.5]{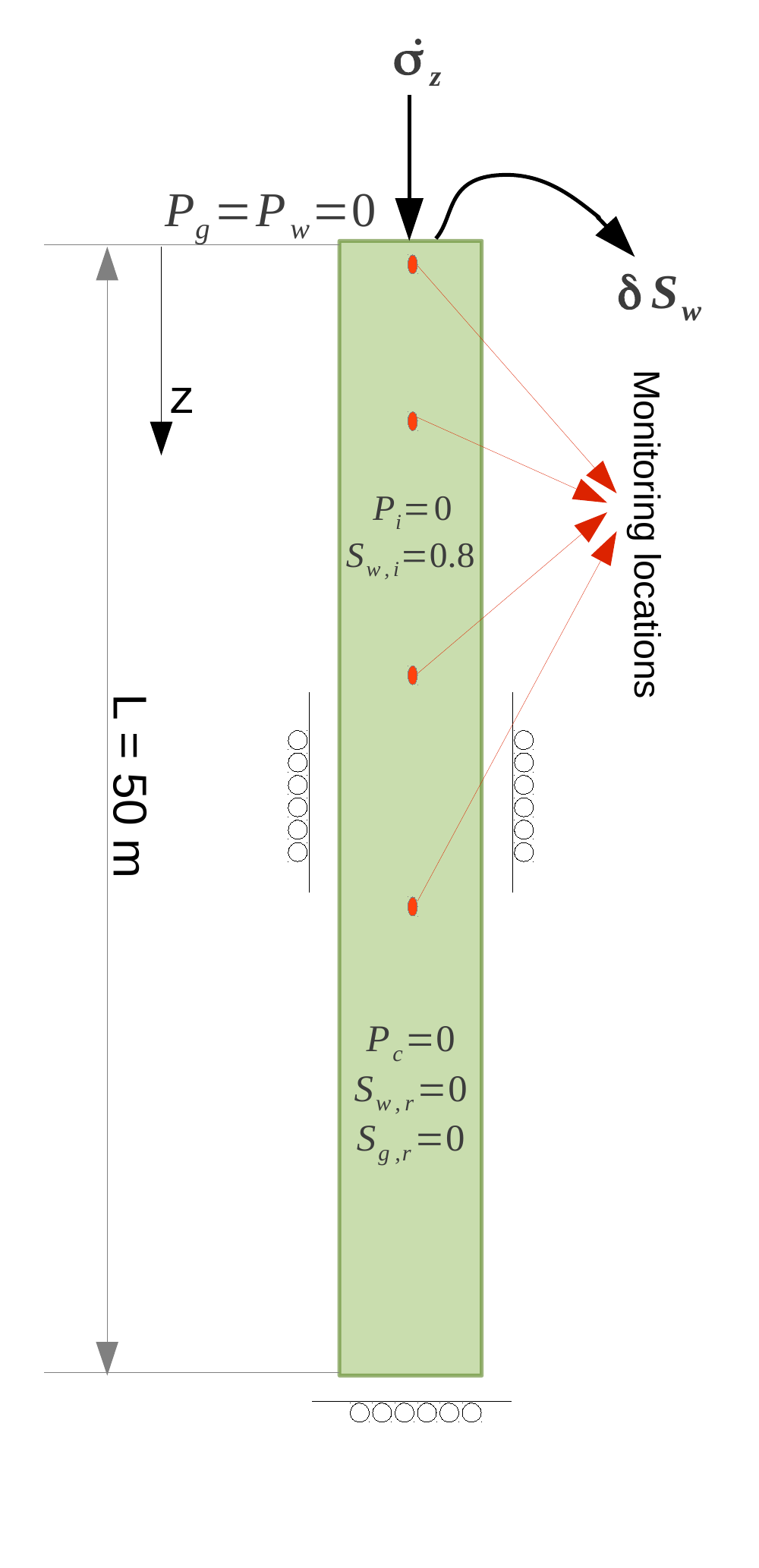}
     \captionof{figure}{Schemmatic of Terzaghi problem}
     \label{fig:terzaghiSchemmatic}
\end{figure}

\paragraph{\textbf{Problem statement} }
The problem set-up consists of a confined soil sample surrounded by a circular ring and placed in a container filled with water. 
The sample is loaded by a constant or ramped vertical stress at its upper surface, and the deformation is measured. 
The lower boundary is impermeable, and the upper boundary is fully drained. 
This is called a \textit{confined compression test} or an \textit{oedometer test}. 
Fig. \ref{fig:terzaghiSchemmatic} shows a schematic for this problem.
It is expected that the compression of a soil sample will be accompanied by an expulsion of pore fluids from the sample. 
Also, if the soil permeability is low, this may take considerable time.
In Terzaghi's original work, the pore fluid and the soil particles were both assumed to be incompressible,
so that the only mechanism of deformation was a rearrangement of the particles. 
However, Biot's more generalized consolidation framework, which is also the basis of our poroelasticity model, accounts for both fluid and soil compressibilities.
So, in the further discussion, the fluid and the soil are treated as compressible mediums.

  \begin{table}
    \centering
    \captionof{table}{Solid-matrix properties for Terzaghi problem}\label{table:solidProperties}
    \begin{tabular}{L{2cm} C{1cm} C{1.25cm} L{2cm}}\hline
      Property 			& Symbol	& Unit		& Value			\\ \hline
      Drained bulk modulus 	& $B_s$ 	& GPa		& $8.0$			\\
      Poisson ratio 		& $\nu$ 	& - 		& $0.20$		\\
      Porosity 			& $\phi$ 	& - 		& $0.19$		\\
      Permeability 		& $K$ 		& $m^2$		& $1.9\times 10^{-13}$	\\
      Biot constant	 	& $\alpha$ 	& -		& $0.8$			\\ \hline
    \end{tabular}
  \end{table}

  \begin{table}
    \centering
    \captionof{table}{Fluid properties for Terzaghi problem}\label{table:fluidProperties}
    \begin{tabular}{L{2cm} C{1cm} C{1.25cm} L{2cm}}\hline
      Property 		& Symbol	& Unit			& Value			\\ \hline
      \multicolumn{4}{l}{\textit{Wetting fluid}}					\\
      Bulk modulus 	& $B_w$ 	& GPa			& $2.933$		\\
      Density 		& $\rho_w$ 	& $\frac{kg}{m^3}$ 	& $997.05$		\\
      Viscosity		& $\mu_w$ 	& $Pa\cdot s$ 		& $8.9008\times 10^{-4}$\\
      \multicolumn{4}{l}{\textit{Non-wetting fluid}}					\\
      Bulk modulus 	& $B_{nw}$ 	& GPa			& $1.187$		\\
      Density 		& $\rho_{nw}$ 	& $\frac{kg}{m^3}$ 	& $997.05$		\\
      Viscosity		& $\mu_{nw}$ 	& $Pa\cdot s$ 		& $8.9008\times 10^{-4}$\\ \hline
    \end{tabular}
  \end{table}

\paragraph{} The mathematical description of such a problem in 1D reduces to a fluid diffusion equation of hydrogeology,
\begin{align}\label{eqn:terzaghiDiffusionEqn}
 \partial_t \overline P - c\ \frac{\partial^2}{\partial z^2} \overline P = 0
\end{align}
where, $\overline P$ is the mean fluid pressure given by $\overline P = S_w P_w + S_{nw}P_{nw}$, and $c$ is 1-D fluid diffusivity. 

\paragraph{} For a vertical load $\sigma_{zz}$ ramped linearly at the top boundary at a rate $d\sigma_{zz}/dt = \dot \sigma_z$, the analytical solution for the pore pressure response is given as
\begin{align} \label{eqn:terzaghiSolution}
 &\frac{\overline P\left(z,t \right)}{\overline P_0} 
  =  1 - \left( \frac{L-z}{L} \right)^2 \notag\\
 &- \frac{32}{\pi^3} \ \sum\limits_{m=0}^{\infty}\left( \frac{-1^m}{\left( 2m + 1\right)^3} \ e^{ \left(-\psi^2 c t\right)} \cos\left[\psi\left( L-z\right)\right]  \right)
\end{align}

\paragraph{} where, $\overline P_0$ is the total pressure generation given as,
\begin{align}
 \overline{P}_0 = \frac{L^2}{2c}\left( H_v \dot \sigma_z \right) \ .
\end{align}

\paragraph{} $\psi = \left( 2m+1 \right) \pi / \left( 2L \right)$, $L$ is the total column length, and $z$ is the location in the column downward from the applied stress.
$H_v$ is the 1-D Skempton coefficient given by,
\begin{align}
 H_v = - \left. \frac{\delta \overline P}{\delta \sigma_{zz}}\right|_{\epsilon_{xx} = \epsilon_{yy} = \zeta = 0} = \frac{\alpha}{B_{sv} S_v}
\end{align}
where, $B_{sv}$ is the uniaxial drained bulk modulus, and $S_v$ is the 1-D specific storage given by $S_v = K / (\mu c)$.\newline
$\mu$ is the fluid mobility given as $\dfrac{1}{\mu}=\dfrac{1}{2}\left( \dfrac{1}{\mu_{w}} + \dfrac{1}{\mu_{nw}} \right)$

\paragraph{}
The complete derivation of the analytical solution (Eqn. \ref{eqn:terzaghiSolution}) can be found in many of the textbooks on soil mechanics, 
for example, Verrujit (2013) \cite{Verruijt2013}.

\paragraph{\textbf{Numerical simulation} } 
For benchmarking, we use the test setting described by Kolditz et al. \cite{Kolditz2012}. 
A soil column of $50$ m is chosen. The properties of the rock material are listed in Table \ref{table:solidProperties}, and that of the two fluid phases are given in Table \ref{table:fluidProperties}.
The column is discretized uniformly into $200$ grid cells. The initial fluid pressure in the column is null, and the initial fluid saturations are $S_w = 0.8$ and $S_{nw} = 0.2$. 
The hydraulic properties are chosen as $P_c = 0$ and $k_{r,w} = k_{r,nw} = 0.5$.
As shown in Fig. \ref{fig:terzaghiSchemmatic}, a load $\sigma_{zz} = 10$ MPa is applied at the top boundary at a loading rate of $\dot \sigma_z = 0.01 \text{and } 0.001 $ MPa/s.
The bottom boundary of the column is subjected to roller displacement BC and the top boundary is allowed to compress freely under the applied load $\dot \sigma_z$.
The top boundary is a free-drainage boundary. All other boundaries are no-flow.

\paragraph{\textbf{Results} }
The results of the numerical simulation for the two loading rates ($\dot\sigma_z = 0.01,\ 0.001$ MPa/s) are presented in Fig. \ref{fig:alpha=0,8_dsigmaBydt=0,01} and Fig. \ref{fig:alpha=0,8_dsigmaBydt=0,001}.
Compression of the column leads to a rapid pressure increase followed by dissipiation of pressure over time from top of the column. 
It can also be observed that for a lower loading rate, the pore-pressure equilibrates faster, whereas, for a higher loading rate the pore-pressure takes longer to dissipate.
The numerical results show a very good agreement with the analytical solution.

\begin{figure}
  \centering
    \subfloat[$\dot \sigma_z = 0.01$ MPa]{%
    \includegraphics[scale=0.281]{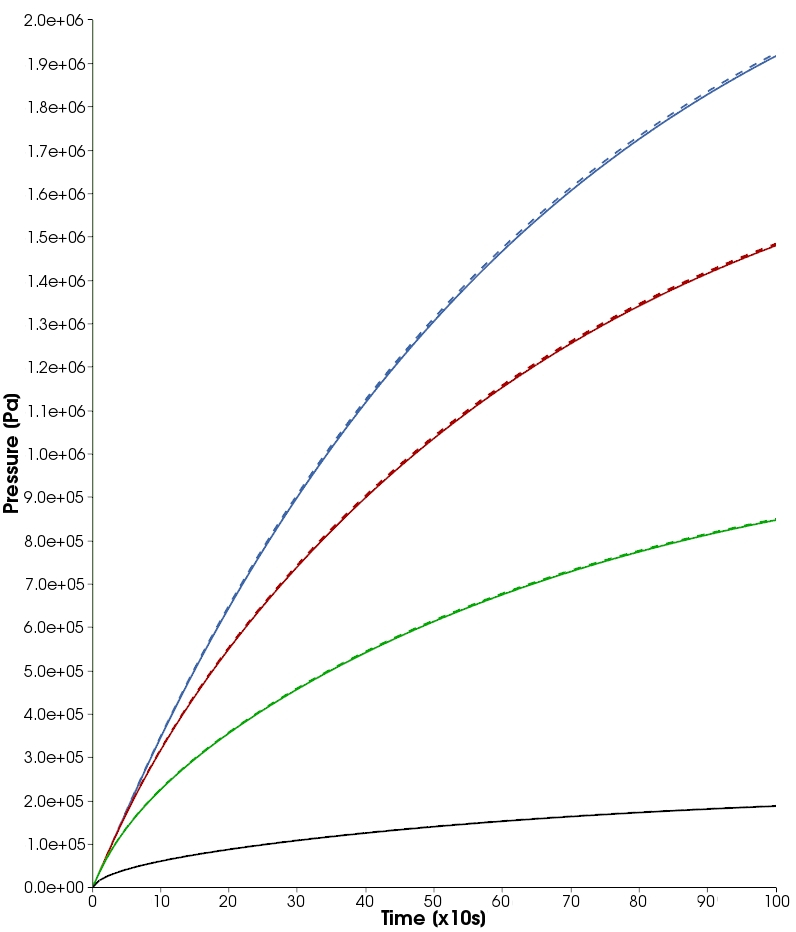}
    \label{fig:alpha=0,8_dsigmaBydt=0,01}}
    \vfill
    \subfloat[$\dot \sigma_z = 0.001$ MPa]{%
    \includegraphics[scale=0.281]{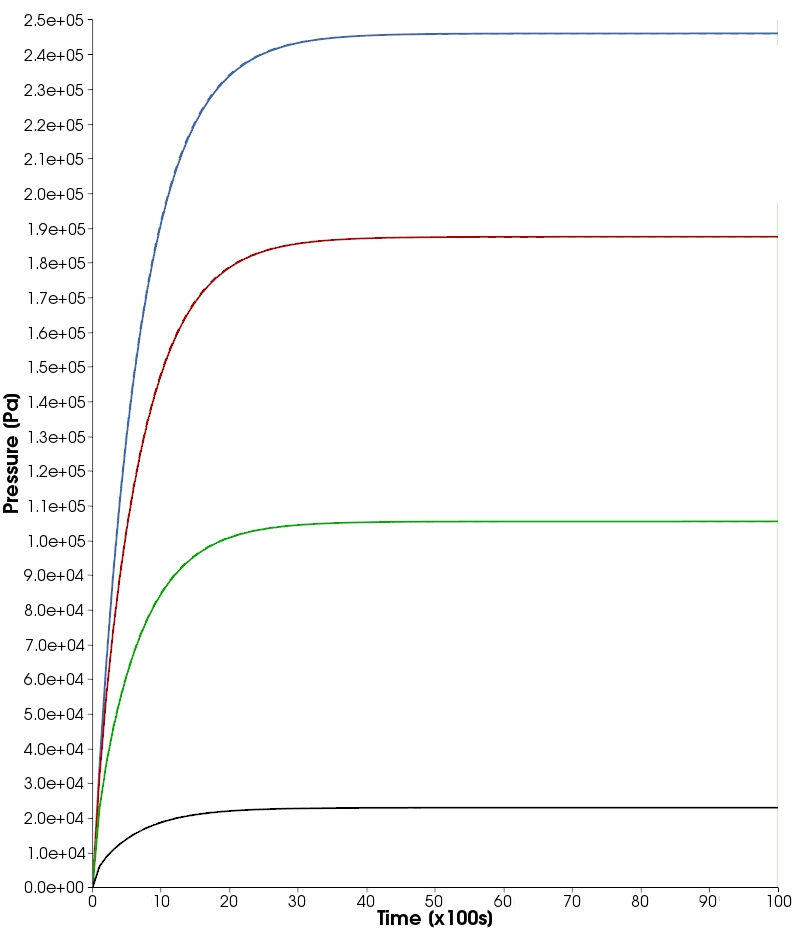}
    \label{fig:alpha=0,8_dsigmaBydt=0,001}}
    
    \vspace{0.4cm}
    \makebox[\linewidth]{\centering\textbf{Legend:} }\vfill
    \makebox[\linewidth]{ \raisebox{.4ex}{ \makebox{ \textcolor{blue}{ \linethickness{0.3mm}\line(1,0){25} } } } numerical $\left(\frac{z}{L} = 0.60 \right)$ \quad
						     \textcolor{blue}{\dashedrule }\ analytical $\left(\frac{z}{L} = 0.60 \right)$ }
    \makebox[\linewidth]{ \raisebox{.4ex}{ \makebox{ \textcolor{red}{ \linethickness{0.3mm}\line(1,0){25} } } } numerical $\left(\frac{z}{L} = 0.40 \right)$ \quad
						     \textcolor{red}{\dashedrule }\ analytical $\left(\frac{z}{L} = 0.40 \right)$ }
    \makebox[\linewidth]{ \raisebox{.4ex}{ \makebox{ \textcolor{green}{ \linethickness{0.3mm}\line(1,0){25} } } } numerical $\left(\frac{z}{L} = 0.20 \right)$ \quad
						     \textcolor{green}{\dashedrule }\ analytical $\left(\frac{z}{L} = 0.20 \right)$ }
    \makebox[\linewidth]{ \raisebox{.4ex}{ \makebox{ \textcolor{black}{ \linethickness{0.3mm}\line(1,0){25} } } } numerical $\left(\frac{z}{L} = 0.04 \right)$ \quad
						     \textcolor{black}{\dashedrule }\ analytical $\left(\frac{z}{L} = 0.04 \right)$ }
	     
    \caption{Terzaghi benchmark test results }
    \label{fig:terzaghiResultsCompressible}
\end{figure}

\subsection{\textbf{Test 4}: Kinetics-poroelastic coupling (KPE test)} \label{subsec:verificationHPE}

\begin{figure}
  \centering
  \includegraphics[scale=0.5]{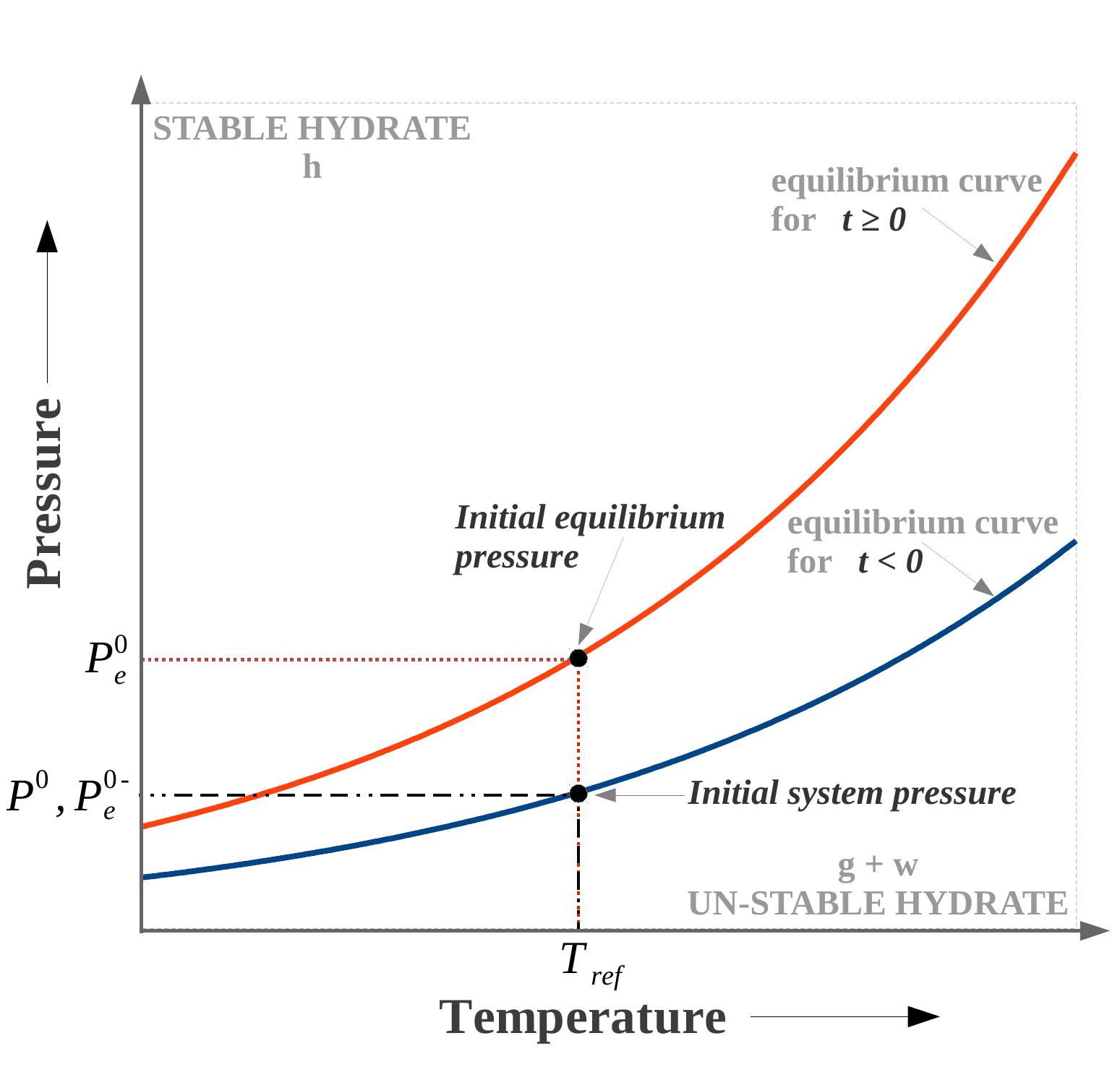}
  \captionof{figure}{KPE test - Hydrate stability curve shift at $t=0$}%
  \label{fig:stabilityCurveShift}
\end{figure}

\paragraph{}
We now extend Terzaghi's $1D$ consolidation problem to include hydrate kinetics in the poroelastic coupling.
We consider a confined soil sample which is uniformly hydrated and fully saturated with water.
A constant vertical stress is applied at the top boundary while the lower boundary is held fixed.
The upper boundary is fully drained, while at the lower boundary the initial pressure is maintained at all times. 
For time $t\leq 0^-$, the thermodynamic state of the sample lies on the hydrate stability curve, so that $P^{0-} = P_e^{0-}$, and the hydrate in the sample is stable (see Fig. \ref{fig:stabilityCurveShift}). 
Here, $P$ indicates the phase pressure, and $P_e$ indicates the equilibrium pressure for hydrate stability.
At time $t=0$, the hydrate stability curve experiences an instantaneous shift such that $P_e^{0}>P_e^{0-}$ (while $P^{0} = P^{0-}$). The hydrate becomes unstable and begins to dissociate.  
This generates excess pore-pressure which prevents the \textit{full consolidation} of the sample.

\begin{figure}
  \centering
  \includegraphics[scale=0.5]{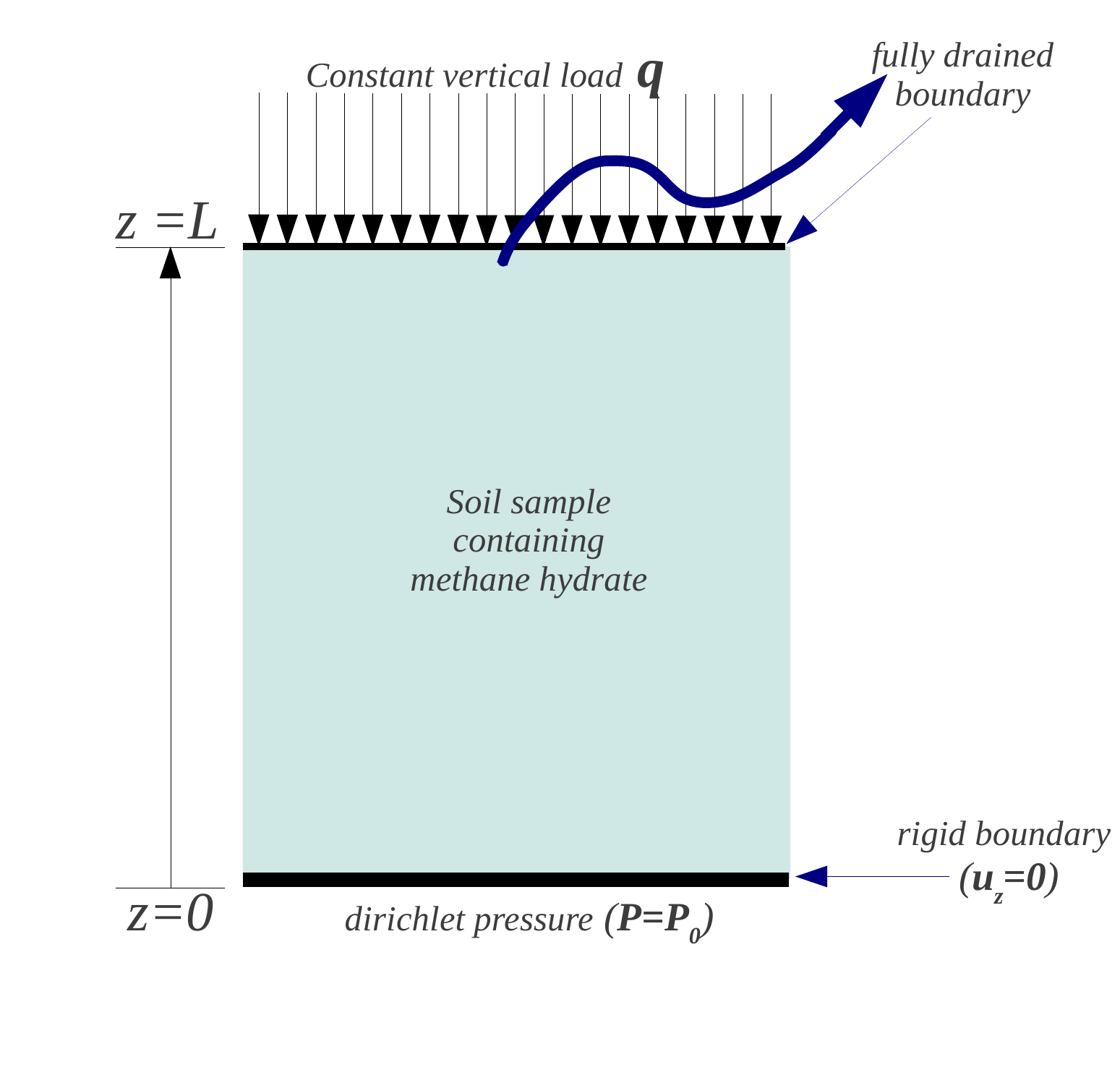}
  \captionof{figure}{KPE test - Problem schematic}%
  \label{fig:HPEproblemschematic}
\end{figure}

\paragraph{}
The schematic for this problem is shown in Fig. \ref{fig:HPEproblemschematic}.
Although highly simplified, this problem helps us to isolate the \textit{poroelastic-kinetic coupling}, 
thus providing a framework for validating the numerical implementation of our hydro-mechanical code for the hydrate reservoir model.  

\subsubsection{Problem statement}\label{subsec:HPEproblemStatement}
\paragraph{} For this problem we make the following additional assumptions:

\begin{itemize}
 \item Gas does not dissolve in water phase, and water vapor is not formed, i.e., $\chi_g^{CH_4} = 1$, and $\chi_w^{H_2O} = 1$.
       Based on this assumption, we rewrite the mass balance equations for water and methane in Section \ref{sec:mathematicalModel} phase-wise instead of component-wise.
 \item All the processes, including hydrate phase change, are isothermal.
 \item There is no suction pressure between the two mobile phases, so that $P_g - P_w = 0$.
       Since the phase pressures are now equal, we drop the subscript and assign the symbol $P$ to the phase pressures throughout this section.
 \item Relative permeabilities are $\ k_{r,g} = k_{r,w} = 0.5$.
 \item Effect of gravity is neglected.
 \item Effect of porosity and hydrate saturation on intrinsic permeability $K$ is neglected, i.e. $K$ is constant.
\end{itemize}

\paragraph{}
Further, we simplify the hydrate reaction kinetics model as
\begin{align}
 & \dot q_g = k_0 \ A_{s} \ M_{g} \ \left(P_{e}-P\right) \notag \\
 & A_{s} = A_{s,0} \ S_h
\end{align}
where $k_0$ is the rate of hydrate dissociation, and $A_{s,0}$ is the specific surface area of the hydrate free sample.
Both $k_0$ and $A_{s,0}$ are assumed to be constant.

\paragraph{}
Using these assumptions, we can reduce the mathematical model described in Section \ref{sec:mathematicalModel} to an ODE for the pressure $P$, given by
\begin{align}\label{eqn:storageEquation}
 \alpha \frac{d}{dt}\epsilon_v + S \frac{d}{dt}P = \nabla \cdot \frac{K}{\mu_{f}} \nabla P + C \ S_h \ \left(P_e - P \right) \ .
\end{align}
A detailed derivation is given in \ref{appendix:A}.

\paragraph{}
The term $S$, called the \textit{Storativity}, is given as \newline
$S = \phi_e \left( \dfrac{S_{w,e}}{ B_{w}} + \dfrac{S_{g,e}}{ B_{g}}\right) + \dfrac{\left( \alpha - \phi_e \right)}{ B_{sh}}$ , \newline
term $C = \left(N_{h} \dfrac{M_{w}}{\rho_{w}} + \dfrac{M_{g}}{\rho_{g}} - \dfrac{M_{h}}{\rho_{sh}}\right) k_0 \ A_{s,0}$ , \newline
$\epsilon_v$ is the volumetric strain given as $\epsilon_v=\nabla\cdot \mathbf{u}$ , and \newline
$\mu_f$ is the fluid mobility given as $\dfrac{1}{\mu_f}=\dfrac{1}{2}\left( \dfrac{1}{\mu_g}+\dfrac{1}{\mu_w} \right)$.
$\alpha$ is Biot's constant.

\paragraph{}
Eqn. \eqref{eqn:storageEquation} is the \textit{storage equation}. 
In this form it can be interpreted as: 
on the REV scale, the compression of the soil consists of compression of the pore-fluids and the compression of the solid particles, 
plus the total volume of fluid expelled from the REV and the fluid generated in the REV.

\paragraph{}
Further, in Eqn. \eqref{eqn:storageEquation}, the term $\left[\alpha \dfrac{d}{dt}\epsilon_v\right]$ is the \textit{mechanical} part, 
and the term $\left[\nabla \cdot \dfrac{\kappa}{\mu_{f}} \nabla P \right]$ is the \textit{flow} part.
The terms $\left[S \dfrac{d}{dt}P\right]$ and $\left[C \ S_h \ \left(P_e - P \right)\right]$ are the \textit{coupling} terms, 
the former for the coupling between the flow and the mechanical models, and the latter for the coupling between the flow and the reaction kinetics models.

\paragraph{}
For the 1D problem under consideration, we can rewrite Eqn. \eqref{eqn:storageEquation} as
\begin{align}\label{eqn:storageEquation1D}
 \alpha \frac{d}{dt}\epsilon_v + S \frac{d}{dt}P = \frac{K}{\mu_{f}} \ \frac{d^2}{dz^2} P + C \ S_h \left(P_e - P \right) \ .
\end{align}

\paragraph{}
In case of 1D deformation, the volumetric strain equals the vertical strain and is induced by the vertical stress $\sigma_{zz}'$,
\begin{align}\label{eqn:volStrainVerticalStress}
 \frac{d}{dt}\epsilon_v = - C_{m} \frac{d}{dt}\sigma_{zz}' = - C_{m} \left( \frac{d}{dt}\sigma_{zz} - \alpha \frac{d}{dt}P \right) \ .
\end{align}
where, $C_m$ is the compressibility of bulk porous material, such that $C_m = \dfrac{1}{B_m}$ \ .

\paragraph{}
Thus, we eliminate $\dfrac{d}{dt}\epsilon_v$ in Eqn. \eqref{eqn:storageEquation1D} using Eqn. \eqref{eqn:volStrainVerticalStress}, which gives
\begin{align}\label{eqn:generalStorageEquation1D}
  \frac{d}{dt}P & = \frac{\alpha C_m}{\alpha^2 C_m + S} \frac{d}{dt}\sigma_{zz} + \frac{\kappa}{\mu_{f}\left( \alpha^2 C_m + S \right)} \ \frac{d^2}{dz^2} P \notag \\
		& + \frac{C \ S_h}{ \alpha^2 C_m + S } \left(P_e - P \right) \ .
\end{align}

\paragraph{}
At time $t=0$, an external load $q$ is instantaneously applied, and the equlibrium pressure of hydrates is instantaneously changed from $P_e^{0-}$ to $P_e^{0}$.
Since both these processes are instantaneous, no fluid is mobilized at $t=0$, i.e., in Eqn. \eqref{eqn:generalStorageEquation1D}, $\dfrac{d^2}{dz^2}P = 0$.
So, from Eqn. \eqref{eqn:generalStorageEquation1D}, we get the initial condition of the sample as
\begin{align}
 t=0 \ : \ P = P^0 & = \frac{\alpha C_{m}}{\alpha^2 C_{m} + S + CS_{h}} \ q \notag \\
		   & + \frac{CS_{h}}{\alpha^2 C_{m} + S + CS_{h}} \ P_e^{0} \ .
\end{align}

\paragraph{}
For $t>0$, the external load remains constant, so $\dfrac{d}{dt}\sigma_{zz} = 0$. The equilibrium pressure also remains constant, i.e., $P_e^{0+} = P_e^{0} = P_e$.
Thus, from Eqn. \eqref{eqn:generalStorageEquation1D},
\begin{align}\label{eqn:PressureODE}
 t > 0 \ : \ \frac{d}{dt}P &= \underbrace{\frac{\kappa}{\mu_{f}\left( \alpha^2 C_m + S \right)}}_{\textstyle{C_v}} \ \frac{d^2}{dz^2} P \notag\\
	    &+ \underbrace{\frac{C \ S_h}{ \alpha^2 C_m + S }}_{\textstyle{C_r}} \left(P_e - P \right) \notag \\
\implies t > 0 \ : \ \frac{d}{dt}P &= C_v \ \frac{d^2}{dz^2} P + C_r \left(P_e - P \right) \ .
\end{align}
$C_v$ is the \textit{consolidation parameter} which comes from the Terzaghi's classical theory of consolidation. 
$C_r$ is the \textit{reaction parameter}. It is indicative of the damping of the normal consolidation due to dissociation kinetics.

\paragraph{}
The boundary conditions at the top and bottom of the sample are
\begin{align}
 & t > 0, \ z = L \ : \ \frac{d}{dz}P = 0  \notag \\
 & t > 0, \ z = 0 \ : \ P = P^0   \ .
\end{align}

\paragraph{}
The ODE in Eqn. \eqref{eqn:PressureODE} is a non-homogeneous ODE. We can homogenize it by choosing a new primary variable $\overline{P}$ such that $\overline{P} = P_e - P $. 
Then the initial-boundary-value problem (IBVP) can be formally summarized as
\begin{empheq}[box=\widefbox]{align}\label{eqn:PressureODE_IBVP}
 0\leq z \leq L, \  & t>0 \ : \ \frac{d}{dt}\overline{P} = C_v \ \frac{d^2}{dz^2} \overline{P} - C_r\ \overline{P} \nonumber \\
 & \nonumber \\
          z = 0, \  & t>0 \ : \ \overline{P} = P_e - P^0 \nonumber \\
          z = L, \  & t>0 \ : \ \frac{d}{dz}\overline{P} = 0 \nonumber \\
 0\leq z \leq L, \  & t=0 \ : \ \overline{P} = P_e - P^0 \nonumber \\
  \text{where, } \qquad  & \nonumber \\
  P^0 = &\frac{\alpha C^0_{m} \ q \ + \ CS_{h}^0 \ P_e }{\alpha^2 C^0_{m} + S^0 + CS_{h}^0} 
\end{empheq}
which is a homogeneous ODE with non-homogeneous boundary conditions. 
An analytical solution can be obtained for this problem using any of the standard techniques for solving ODEs.
The final solution for $P$ can be written as
  \begin{align}\label{eqn:AnalyticalSolutionP}
 \frac{P_e-P(z,t)}{P_e-P^0\qquad } & = \frac{\cosh \left( \theta \left( L-z \right) \right)}{\cosh\left( \theta L \right)} \notag\\
  & + \frac{2}{L}\ \sum \limits_{n=1}^{\infty}\left( \frac{1}{\lambda_n} \left[ 1 - \frac{\lambda_n^2}{\lambda_n^2 + \theta^2 }\right] \sin\left( \lambda_n z \right) \right.\notag\\
  & \left. \phantom{ \frac{1}{\lambda_n} } \ \exp \left[ - C_v \left( \lambda_n^2 + \theta^2 \right) t \right] \right) \ .
 \end{align}

\paragraph{}

\subsubsection{Numerical simulation}\label{subsec:HPEnumericalSolution}
\paragraph{\textbf{Test setting} }
We chose a sample of length $L=1$ m containing $30\%$ hydrate by volume. 
The sample is initially fully water saturated and is contained in a pressure vessel at $P^0 = 6$ MPa.
A constant external load $q = 10$ MPa is applied at the top boundary, i.e., at $z = L = 1$. 
At the bottom boundary, i.e., $z=0$, the pressure is held constant at the initial value.

\begin{figure}
  \centering
  \includegraphics[scale=0.5]{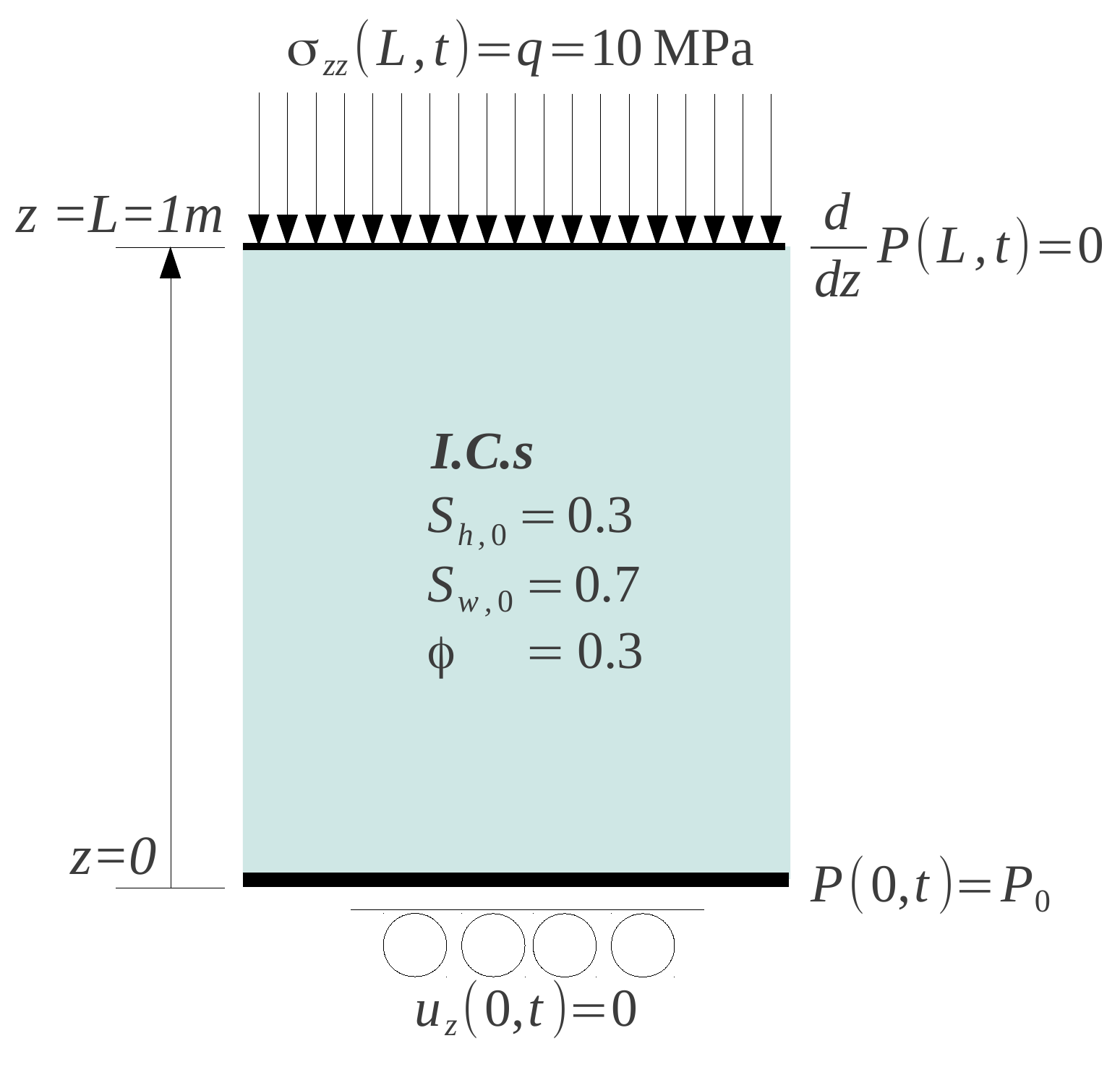}
  \captionof{figure}{KPE test - Computational domain settings}%
  \label{fig:HPE_testSchematic}
\end{figure}

\paragraph{}
Fig. \ref{fig:HPE_testSchematic} shows the domain specifications, the initial conditions, and the boundary conditions. 
The material properties are listed in Table \ref{table:HPE_materialProperties}.

\begin{table}
  \centering
    \captionof{table}[t]{KPE test - Material properties}
    \label{table:HPE_materialProperties}
  \vskip\abovecaptionskip
   \begin{tabular}[h]{L{2cm}  L{0.8cm} L{1.5cm}  C{2cm} }
   \hline 
   Property		& Symbol	& Unit			& Value 		\\ \hline
   \multicolumn{4}{c}{Water phase} 	\\ \hline
   Density		& $\rho_w$	& $kg \cdot m^{-3}$	& $997.05$ 		\\
   Molar mass		& $M_w$		& $kg \cdot mol^{-1} $	& $0.018$ 		\\
   Dynamic viscosity	& $\mu_w$	& $Pa \cdot s$		& $8.9008\times10^{-3}$ \\
   Bulk modulus		& $B_w$		& $GPa$			& $2.933$ 		\\ \hline
   \multicolumn{4}{c}{Gas phase} 	\\ \hline
   Density		& $\rho_g$	& $kg \cdot m^{-3}$	& $0.717$ 		\\
   Molar mass		& $M_g$		& $kg \cdot mol^{-1} $	& $0.016$ 		\\
   Dynamic viscosity	& $\mu_g$	& $Pa \cdot s$		& $1.0245\times10^{-5}$ \\
   Bulk modulus		& $B_g$		& $GPa$			& $0.1013$ 		\\ \hline
   \multicolumn{4}{c}{Hydrate phase} 	\\ \hline
   Density		& $\rho_h$	& $kg \cdot m^{-3}$	& $900$ 		\\
   Molar mass		& $M_h$		& $kg \cdot mol^{-1} $	& $0.119$ 		\\
   Hydration number	& $N_h$		& $-$			& $5.75$ 		\\
   Young's modulus	& $E_h$		& $GPa$			& $1.35$ 		\\ \hline
   \multicolumn{4}{c}{Soil phase} 	\\ \hline
   Density		& $\rho_s$	& $kg \cdot m^{-3}$	& $700$ 		\\
   Surface area		& $A_{s,0}$	& $m^{2}$		& $10^5$ 		\\
   Young's modulus	& $E_s$		& $GPa$			& $0.3$ 		\\ \hline
   \multicolumn{4}{c}{Solid composite} 	\\ \hline
   Poisson ratio	& $\nu_{sh}$	& $-$			& $0.2$ 		\\
   Biot constant	& $\alpha$	& $-$			& $0.8$ 		\\ \hline
  \end{tabular}
\end{table}

\paragraph{}
The storage equation governing this problem, Eqn. \eqref{eqn:PressureODE_IBVP}, contains two parameters: $C_v$ and $C_r$. 
To test the numerical implementation, we chose three different values of $C_v$ and $C_r$, each with a different order of magnitude.
Therefore, we run nine tests with all combinations of the chosen $C_v$ and $C_r$.
We control the parameter $C_v$ by varying the sample permeability $\kappa$ and the parameter $C_r$ by varying the dissociation rate constant $k_0$.
For each test, the value of $P_e$ is chosen such that the initial condition of no-drainage is satisfied.  
The control parameters for each of the nine tests are listed in Table \ref{table:HPEtestParameters}.
It can be observed in Table \ref{table:HPEtestParameters} that as the reaction rate constant $k_0$ increases the value of equilibrium pressure $P_e$ decreases.
This is due to the no-drainage condition at $t=0$. The faster the dissociation, the more the generated fluids will mobilize.
Conversely, the slower the dissociation reaction, the higher margin we get for raising $P_e$ without instantaneously mobilizing the fluids.

\begin{table}
\centering
\caption{KPE test - Control parameters}\label{table:HPEtestParameters}
 \begin{tabular}[h]{C{0.4cm} | R{0.8cm} | R{1.25cm} | R{1.35cm} | R{1cm} | R{0.9cm}}
 \hline
  test ID & $\kappa$ $ \left[ mD \right]$	& $k_0$ \scalebox{0.75}{$\times a^*$} \scalebox{1}{$\left[\frac{mol}{m^2\ Pa\ s}\right]$}	& $C_v		$	& $C_r		$ 	& $P_e $ $\left[ \text{MPa} \right] $	 \\  [0.3em]\hline
 1	 & $0.1		$			& $360		$			& $1.53755 	$	& $0.289504 	$ 	& $19.151 				$ \\ 
 2	 & $0.01	$			& $360	 	$			& $0.153755	$	& $0.289504 	$ 	& $19.151 				$ \\ 
 3	 & $0.001	$			& $360	 	$			& $0.0153755	$	& $0.289504 	$ 	& $19.151 				$ \\ 
 4	 & $0.1		$			& $3600		$			& $1.53755 	$	& $2.89504 	$ 	& $7.315 				$ \\ 
 5	 & $0.01	$			& $3600 	$			& $0.153755 	$	& $2.89504 	$ 	& $7.315 				$ \\ 
 6	 & $0.001	$			& $3600 	$			& $0.0153755 	$	& $2.89504 	$ 	& $7.315 				$ \\ 
 7	 & $0.1		$			& $36000	$			& $1.53755 	$	& $28.9504 	$ 	& $6.132 				$ \\ 
 8	 & $0.01	$			& $36000 	$			& $0.153755 	$	& $28.9504 	$ 	& $6.132 				$ \\ 
 9	 & $0.001	$			& $36000	$			& $0.0153755 	$	& $28.9504	$ 	& $6.132 				$ \\ \hline 
\end{tabular}
\vspace{0.1cm}\\
\scalebox{0.75}{${}^* a = \exp{\left(-\Delta E_{a}/(R T)\right)}$, where, $\Delta E_{a}/R=9400 \ K$ and $T = 283.15 \ K$}
\end{table}

\paragraph{\textbf{Simulation and results}}
The domain is discretized into $400$ cells in z-direction, and the problem is solved in $1D$.
The time-step is kept constant at $t=0.1$ s and the simulation is run until $t_{end} = 60$ s.

\begin{figure*}
  \centering
    \subfloat[TestID 1]{%
    \includegraphics[scale=0.265]{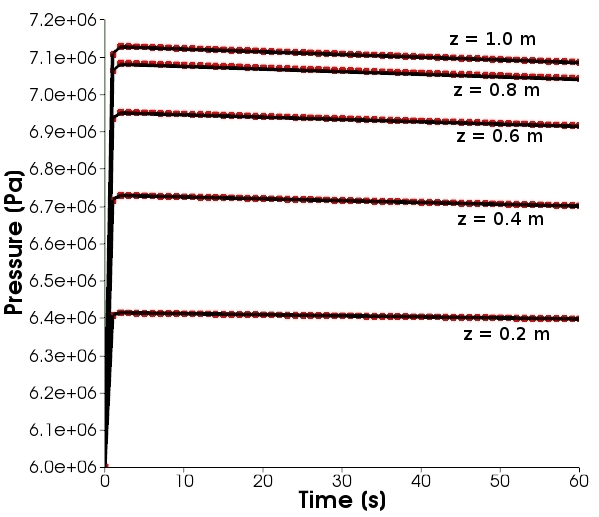}
    \label{fig:testID1_PoverTime}}
    \hfill
    \subfloat[TestID 2]{%
    \includegraphics[scale=0.265]{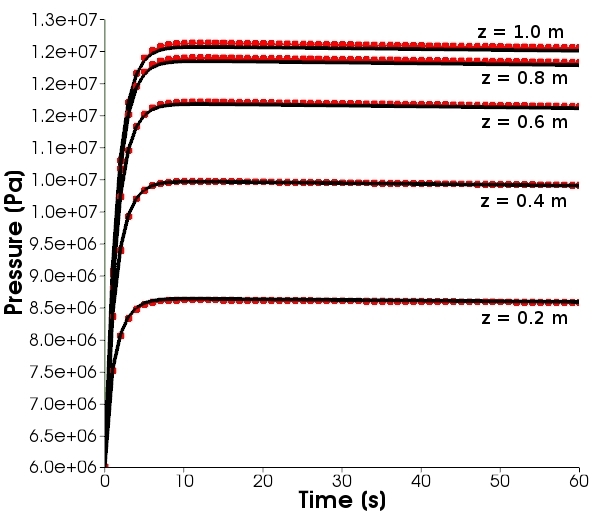}
    \label{fig:testID2_PoverTime}}
    \hfill
    \subfloat[TestID 3]{%
    \includegraphics[scale=0.265]{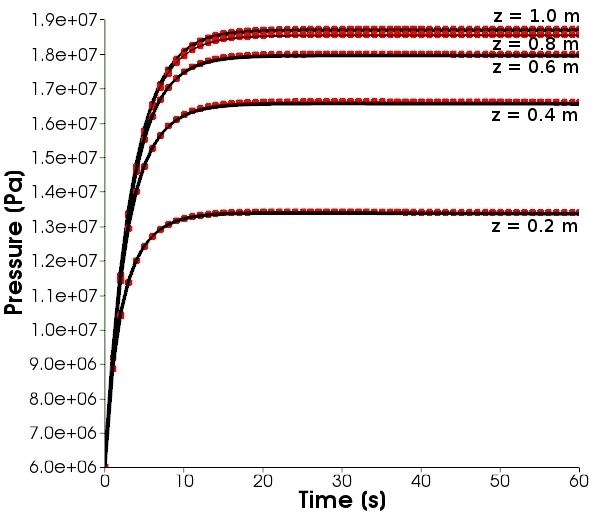}
    \label{fig:testID3_PoverTime}}
    
    \vfill
    
    \subfloat[TestID 4]{%
    \includegraphics[scale=0.265]{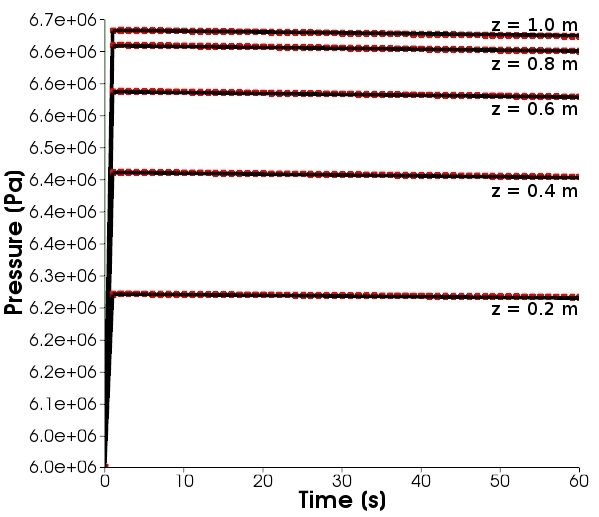}
    \label{fig:testID4_PoverTime}}
    \hfill
    \subfloat[TestID 5]{%
    \includegraphics[scale=0.265]{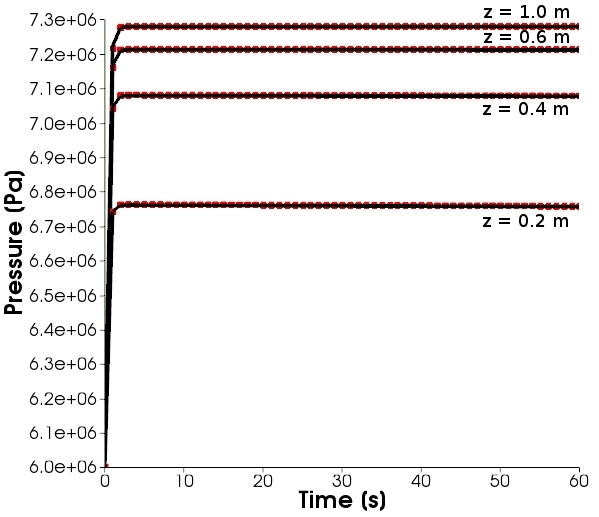}
    \label{fig:testID5_PoverTime}}
    \hfill
    \subfloat[TestID 6]{%
    \includegraphics[scale=0.265]{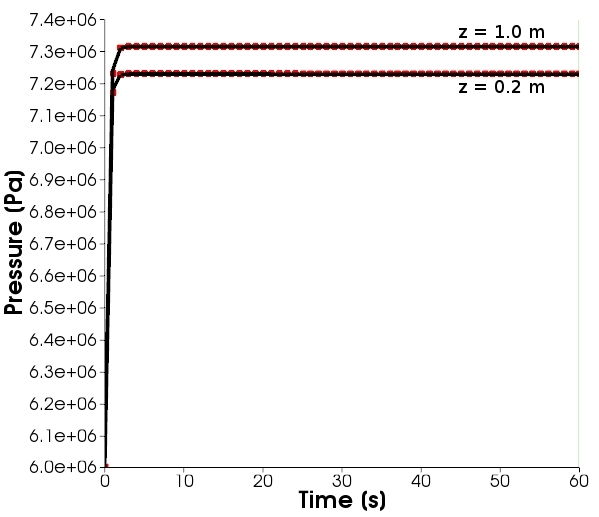}
    \label{fig:testID6_PoverTime}}
    
    \vfill
    
    \subfloat[TestID 7]{%
    \includegraphics[scale=0.265]{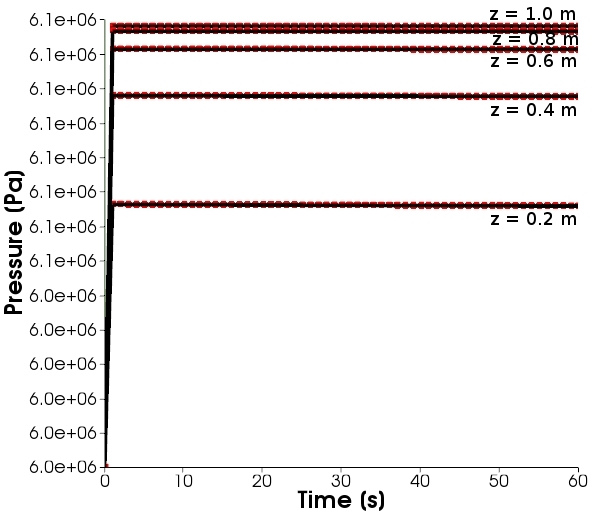}
    \label{fig:testID7_PoverTime}}
    \hfill
    \subfloat[TestID 8]{%
    \includegraphics[scale=0.265]{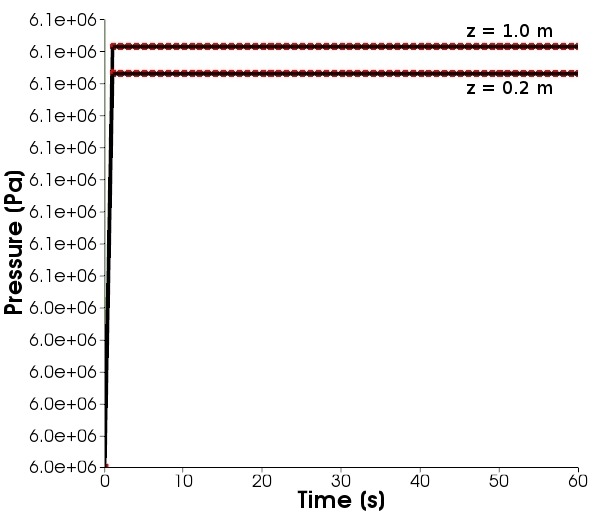}
    \label{fig:testID8_PoverTime}}
    \hfill
    \subfloat[TestID 9]{%
    \includegraphics[scale=0.265]{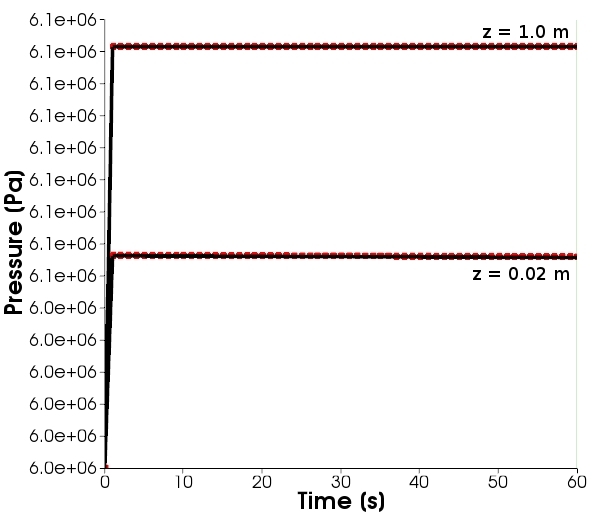}
    \label{fig:testID9_PoverTime}}
    
    \vspace{0.25cm}
    
    \makebox[\linewidth]{\textbf{Legend:} 
	     \quad \raisebox{.6ex}{\makebox{\textcolor{red}{ \circle*{3}\circle*{3}\circle*{3}\circle*{3}\circle*{3}\circle*{3}\circle*{3}\circle*{3}\circle*{3} }}} numerical solution \quad ,
	     \quad \raisebox{.4ex}{\makebox{\linethickness{0.4mm}\line(1,0){25}}} \ analytical solution}
    
    \captionsetup{justification=centering}
    \caption{KPE test - Comparison of numerical and analytical solutions  \newline
	     Note: For \textit{TestID 5}, pressure profiles at $z=1$ m and $0.8$ m are equal, for \textit{TestID 6} and \textit{TestID 8}, pressure profile at $z=1$ m, $0.8$ m, $0.6$ m, and $0.4$ m are equal,
	     and for \textit{TestID 9} pressure profile at $z=1$ m, $0.8$ m, $0.6$ m, $0.4$ m, and $0.2$ m are equal.
	     }
    \label{fig:KPE_PoverTime}
\end{figure*}

\paragraph{}
In Fig. \ref{fig:KPE_PoverTime}, the numerically computed pressure values for each test case $1-9$ are compared with the 
analytical pressure $P(z,t)$ obtained from Eqn. \eqref{eqn:AnalyticalSolutionP}. 
For each test case, the pressure solutions are plotted over time at the observation points $z = 1$ m, $0.8$ m, $0.6$ m, $0.4$ m and $0.2$ m.
The plots show a good agreement between the numerical and the analytical solutions signifying that the poroelastic - reaction kinetics coupling terms are correctly handled in the numerical code. 

\begin{figure}
 \centering
 \includegraphics[scale=0.32]{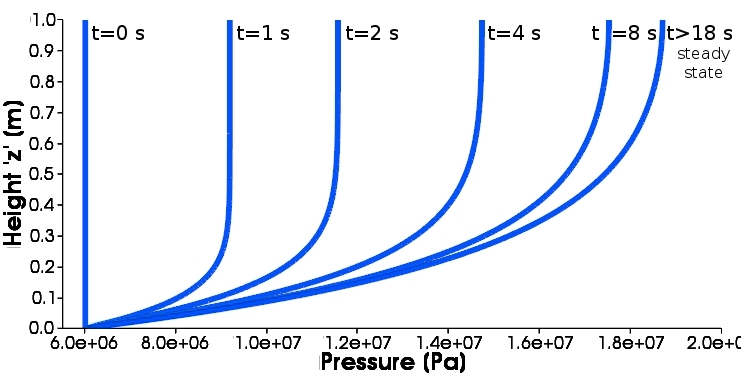}
 \caption{KPE test - testID 3: Time-evolution of the numerical solution of P along the height of the sample}
 \label{fig:testID3_POverLine}
\end{figure}

\begin{figure}
  \centering
  \includegraphics[scale=0.5]{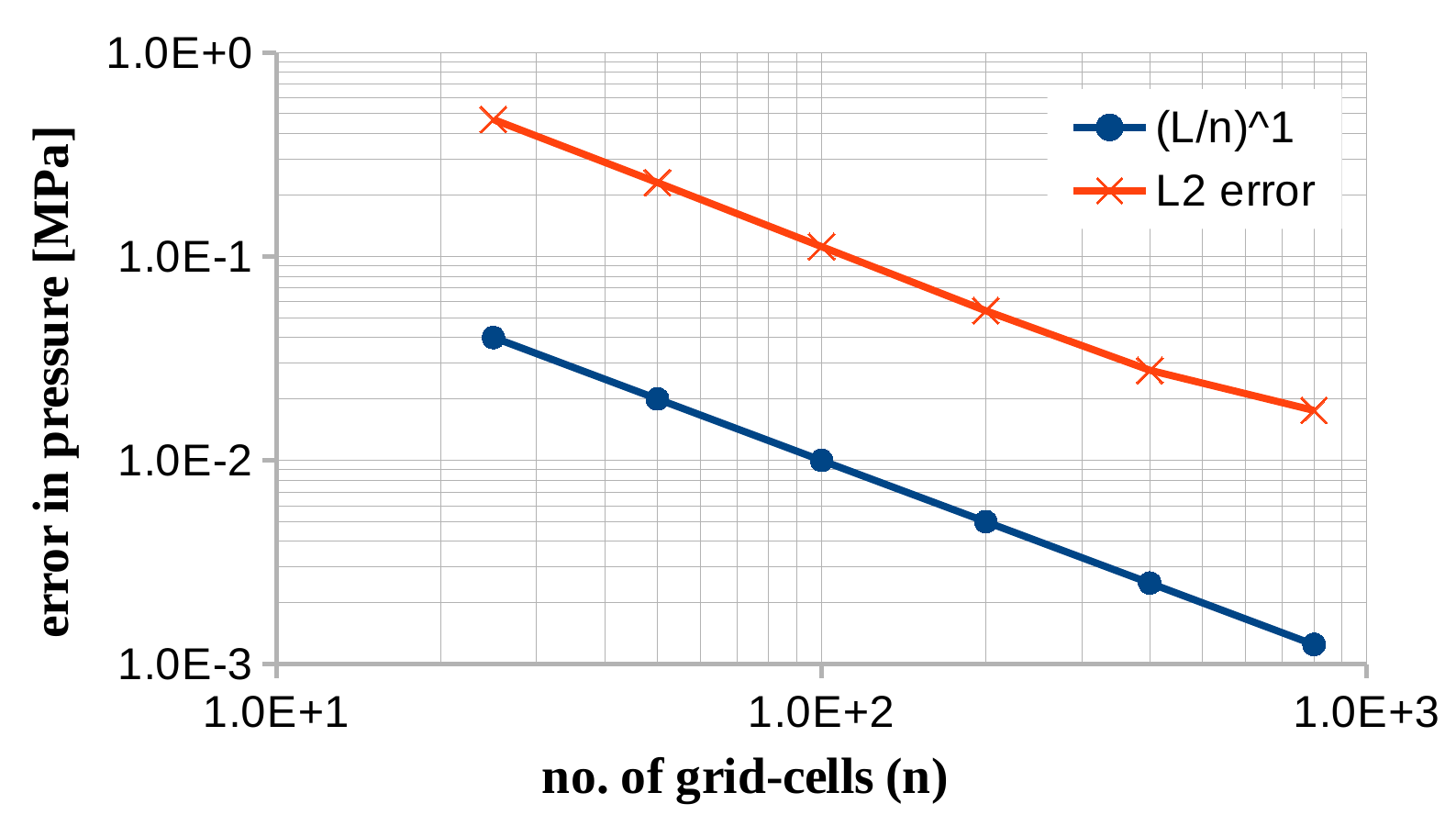}
  \captionof{figure}{KPE test - L2-error vs n}%
  \label{fig:HPE_L2Error}
\end{figure}

\paragraph{\underline{testID-3 } }
The pressure build-up along the length of the sample is plotted for \textit{testID 3} in Fig. \ref{fig:testID3_POverLine}.
Also, for \textit{testID 3} a grid-convergence study is performed. The mesh is refined from $n=25$ cells up to $n=800$ cells, and correspondingly, 
the time-step size is reduced from $\delta t=2$s down to $\delta t=0.0625$s such that the ratio $\dfrac{\delta z}{\delta t}$ remains constant.
The L2 error is calculated at $t=10$s for each refinement and is plotted against the number of grid-cells $n$ on a log-log graph in Fig. \ref{fig:HPE_L2Error}. 
It can be seen that the error decays linearly with refinement. This is in line with the finite volume discretization technique.


\subsection{\textbf{Test 5}: $3D$ hydrate reservoir problem} \label{subsec:3Dexample}
\paragraph{}
So far we considered examples which focussed on systematically isolated couplings and model components.
In this section, we present a more complex example where we simulate the hydro-geomechanical processes in a subsurface hydrate reservoir 
which is destabilized by depressurization using a low pressure gas well.
This example puts together all the important components of our model including dissociation kinets, non-isothermal effects, multi-phase multi-component fluid flow, and geo-mechanics,
and qualitatively shows the effects and counter effects of various physical processes occuring in the hydrate reservoir.
The objective of this example is to give a first idea about the capabilities of our hydrate reservoir model.
A detailed quantitative analysis of the problem and parameter sensitivity study is however beyond the scope of this paper.

  \begin{figure}
    \centering
      \includegraphics[scale=0.45]{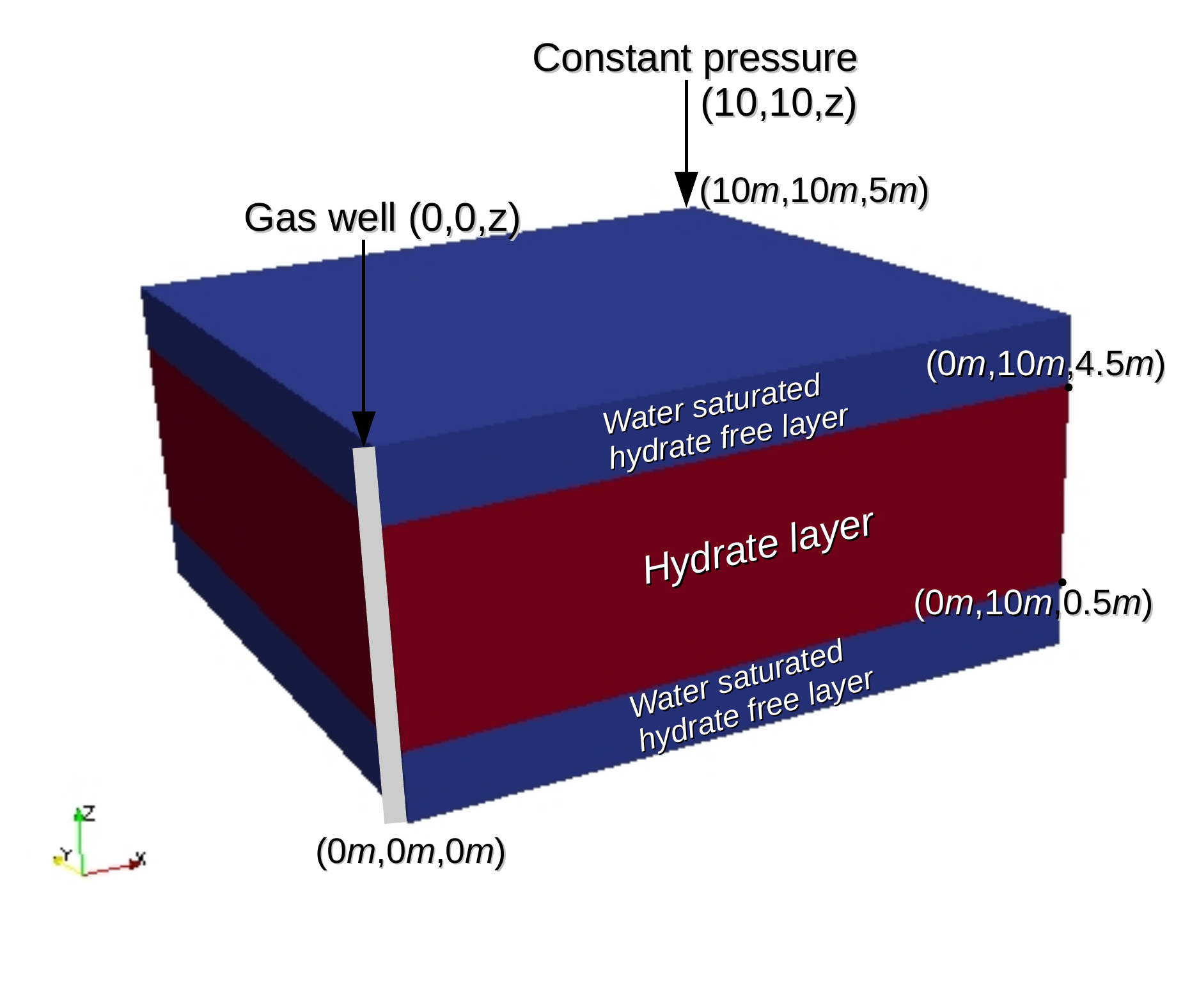}
      \captionof{figure}{Schemmatic of depressurized $3D$-reservoir example}
      \label{fig:3DtestSchemmatic}
  \end{figure}

\paragraph{\textbf{Test setting}}
We consider a scaled down $3D$ reservoir with dimensions $10m\times10m\times5m$, as shown in Fig. \ref{fig:3DtestSchemmatic}.
The hydrate is homogeneously distributed in a $4m$ thick layer lying between $0.5m\leq z \leq 4.5m$, and has a saturation of $40\%$ by volume. 
The reservoir is fully saturated with water and has an initial pressure of $10$ MPa. 
The reservoir is depressurized through a low pressure gas well located at $(0,0,z)$.
The pressure in the gas well is maintained at $P_{well}=4$ MPa.
A constant vertical load of $10$ MPa is acting on the top boundary of the reservoir (i.e. at $z=10$ m).
The initial and the boundary conditions are listed in Table \ref{table:3Dtest_ICValues} and Table \ref{table:3Dtest_BCValues}, respectively.
The material properties and other model parameters are listed in Table \ref{table:3Dtest_properties}.

  \begin{table}
    \centering
    \captionof{table}{Initial conditions for depressurized $3D$-reservoir example }\label{table:3Dtest_ICValues}
   \begin{tabular}{L{4.75cm} L{2.5cm}}	\hline \\ [-0.7em]
    \textbf{Hydrate layer} 				\\
						& $P_{eff,i} = 10$ MPa 		\\
						& $S_{h,i} = 0.4$ 		\\
    at $t=0$, for				& $S_{g,i} = 0$			\\ 
    $0\leq x,y \leq 10$ , $0.5\leq z \leq 4.5$	& $T_{i} = 10 \ ^0$C		\\
						& $K_{i} = 0.0198$ mD		\\
						& $\phi_{eff,i} = 0.18$		\\ [0.4em] \hline \\ [-0.7em]
    \textbf{Hydrate-free layers} \\
						& $P_{eff,i} = 10$ MPa 		\\
    at $t=0$, for				& $S_{h,i} = 0$ 		\\
    $0\leq x,y \leq 10$ , $z < 0.5$		& $S_{g,i} = 0$			\\ 
    and						& $T_{i} = 10 \ ^0$C		\\
    $0\leq x,y \leq 10$ , $z > 4.5$					& $K_{i} = 0.1$ mD		\\
						& $\phi_{eff,i} = 0.3$		\\ [0.4em] \hline    
   \end{tabular}
  \end{table}
  
  \begin{table}
   \centering
   \captionof{table}{Boundary conditions for depressurized $3D$-reservoir example }\label{table:3Dtest_BCValues}
   \begin{tabular}{L{5cm} L{2.25cm}}
    \multicolumn{2}{c}{FLOW model} 						\\ \hline
    Gas well at 			& $P_g = 4$ MPa 			\\
    $x=0$, $y=0$, $0\leq z \leq 5$	& $S_w = 0$ 				\\
				 	&$\nabla\cdot T = 0$			\\ \hline
    Pressure constraint at		& $P_{eff} = P_{eff,i}$ 		\\
    $x=10$, $y=10$, $0\leq z \leq 5$	& $S_w = S_{w,i}$ 			\\
					& $T = T_i$				\\ \hline
    No-flow and adiabatic conditions 	& $\mathbf{v}_g \cdot  \hat n = 0$ 	\\
    on remaining boundaries, i.e.,	& $\mathbf{v}_w \cdot  \hat n = 0$ 	\\
					& $\nabla\cdot T = 0$			\\ \hline      
    \\
    \multicolumn{2}{c}{GEOMECHANICAL model} 								\\ \hline
    \multicolumn{2}{l}{Top boundary} 									\\
    $0\leq x\leq 10$ , $0\leq y\leq 10$ , $z=5$	& $\sigma_{zz}=10$ MPa , $\sigma_{xy}=\sigma_{yx}=0$	\\ 
    \multicolumn{2}{l}{Bottom boundary} 								\\
    $0\leq x\leq 10$ , $0\leq y\leq 10$ , $z=0$ & $u_z=0$ , $\sigma_{xy}=\sigma_{yx}=0$			\\
    \multicolumn{2}{l}{Remaining boundaries}			 					\\ 
						& $u_x=u_y=0$ , $\sigma_{zz}=0$ 			\\ \hline
   \end{tabular}

  \end{table}

  \begin{table}
   \centering
    \captionof{table}{Material properties and model parameters for depressurized $3D$-reservoir example }\label{table:3Dtest_properties}
   \begin{tabular}{L{0.3cm} R{2.5cm} R{1.5cm} L{2cm}} \hline
      \multicolumn{4}{c}{\textbf{Thermal conductivities}} 		\\[0.3em]
      $k^c_g$	& \multicolumn{2}{r}{$- 0.886\times 10^{-2}\ \ \ +  \ 0.242\times 10^{-3}\ T\ $}	& $W\cdot m^{-1}\cdot K^{-1}$	\\
		& \multicolumn{2}{r}{$- 0.699\times 10^{-6}\ T^2 + 0.122 \times 10^{-8}\ T^3$} & \\
      $k^c_w$				&\multicolumn{2}{r}{$0.3834\ ln(T) - 1.581$}	& $W\cdot m^{-1}\cdot K^{-1}$	\\
      \multicolumn{2}{l}{$k^c_h$}	&$2.1$				& $W\cdot m^{-1}\cdot K^{-1}$	\\
      \multicolumn{2}{l}{$k^c_s$}	&$1.9$				& $W\cdot m^{-1}\cdot K^{-1}$	\\
      \multicolumn{4}{c}{\textbf{Specific heat capacities}}		\\
      \multicolumn{2}{l}{$Cp_w$}	&$4186$			& $J\cdot kg^{-1}\cdot K^{-1}$ \\
      $Cv_w$				&\multicolumn{2}{r}{$Cp_w+R_{H_2O}$}	& $J\cdot kg^{-1}\cdot K^{-1}$ \\[0.3em]
      \multicolumn{2}{l}{$Cv_h$}	&$2700$			& $J\cdot kg^{-1}\cdot K^{-1}$ \\
      \multicolumn{2}{l}{$Cv_s$}	&$800$			& $J\cdot kg^{-1}\cdot K^{-1}$ \\
      \multicolumn{4}{c}{\textbf{Dynamic viscosities}} \\[0.3em]
      $\mu_g$	&\multicolumn{2}{r}{$10.4 \ e{-6} \left(\frac{273.15 + 162 }{T+162}\right) \left(\frac{T}{273.15}\right)^{1.5}$}	& $Pa\cdot s$	\\[0.5em]
      $\mu_w$	& \multicolumn{2}{r}{$0.001792 \ \exp\left[ - 1.94 - 4.80 \frac{273.15}{T} \right.$}	& $Pa\cdot s$	\\
		&  \multicolumn{2}{r}{$\left. + 6.74 \left(\frac{273.15}{T}\right)^2\right]$} & \\
      \multicolumn{4}{c}{\textbf{Densities}} \\
      \multicolumn{2}{l}{$\rho_g$}	&$\frac{P_g}{z R_g T}$	& $kg\cdot m^{-3}$ \\[0.6em]
      $\rho_w$	&vapour:		&$0.0022\frac{P_g}{T}$	& $kg\cdot m^{-3}$ \\[0.6em]
		&liquid:		&$1000$			& $kg\cdot m^{-3}$ \\[0.6em]
      $\rho_h$	& \multicolumn{2}{r}{$900 \left[ 1 + \left(\frac{\alpha_{biot}-\phi}{1-\phi}\right)\frac{\Delta P_{eff}}{B_{sh}} \right. $} & $kg\cdot m^{-3}$ \\
		& \multicolumn{2}{r}{$\left. - \left(\frac{\lambda_{sh} + (2/3)G_{sh} }{B_{sh}}\right) \frac{\nabla\cdot\mathbf{u}}{1-\phi_{eff}} \right]$} & \\[0.6em]
      $\rho_s$	& \multicolumn{2}{r}{$2100 \left[ 1 + \left(\frac{\alpha_{biot}-\phi}{1-\phi}\right)\frac{\Delta P_{eff}}{B_{sh}} \right. $} & $kg\cdot m^{-3}$ \\
		& \multicolumn{2}{r}{$\left. - \left(\frac{\lambda_{sh} + (2/3)G_{sh} }{B_{sh}}\right) \frac{\nabla\cdot\mathbf{u}}{1-\phi_{eff}} \right]$} & \\[1em]
      \multicolumn{4}{c}{\textbf{Hydraulic properties}} \\
      \multicolumn{2}{l}{$\lambda_{BC}$}				& $1.2$		&	\\
      \multicolumn{2}{l}{$P_{entry}$}				& $50$		& $kPa$	\\
      \multicolumn{2}{l}{$m$,$a$ (Eqn. \eqref{eqn:Pc_Scaled})}	& $3$, $2$	& \\
      \multicolumn{4}{c}{\textbf{Hydrate kinetics}} \\
      \multicolumn{2}{l}{$k^0_{reac}$}								& $3.6\times 10^{4}$		& $\text{mol}\cdot m^{-2}\cdot Pa^{-1}\cdot s^{-1}$\\
      \multicolumn{2}{l}{$\Delta E_a / R$} 							& $9752.73$			& $K$ \\
      \multicolumn{2}{l}{$N_{Hyd}$}								& $5.75$			& \\
      \multicolumn{2}{l}{$A_1$, $A_2$, $A_3$ (Eqn. \eqref{eqn:hydrateEquilibriumPressure})} 	& $1000$, $38.98$, $8533.8$	& \\
      \multicolumn{2}{l}{$B_1$, $B_2$ (Eqn. \eqref{eqn:hydrateHeatOfDissociation})} 		& $56599$, $16.744$		& \\
      \multicolumn{2}{l}{$\Gamma_r$ (Eqn. \eqref{eqn:reactionSurfaceArea})} 			& $\phi \ S_h$			& \\
      \multicolumn{4}{c}{\textbf{Poroelasticity parameters}} \\
      \multicolumn{2}{l}{$\alpha_{biot}$}	& $0.6$		& \\
      \multicolumn{2}{l}{$\nu_{sh}$}		& $0.2$		& \\
      \multicolumn{2}{l}{$E_{s0}$}		& $0.3$		& $GPa$\\
      \multicolumn{2}{l}{$E_{h}$}		& $1.35$	& $GPa$\\ 
      \multicolumn{2}{l}{$b$, $c$, $d$ (Eqn. \eqref{eqn:SantamarinaRuppelEsh})} & $0$, $1$, $1$	& \\ \hline
   \end{tabular}

  \end{table}

\begin{figure*}
  \centering
    \subfloat[Hydrate saturation]{%
      \includegraphics[scale=0.2]{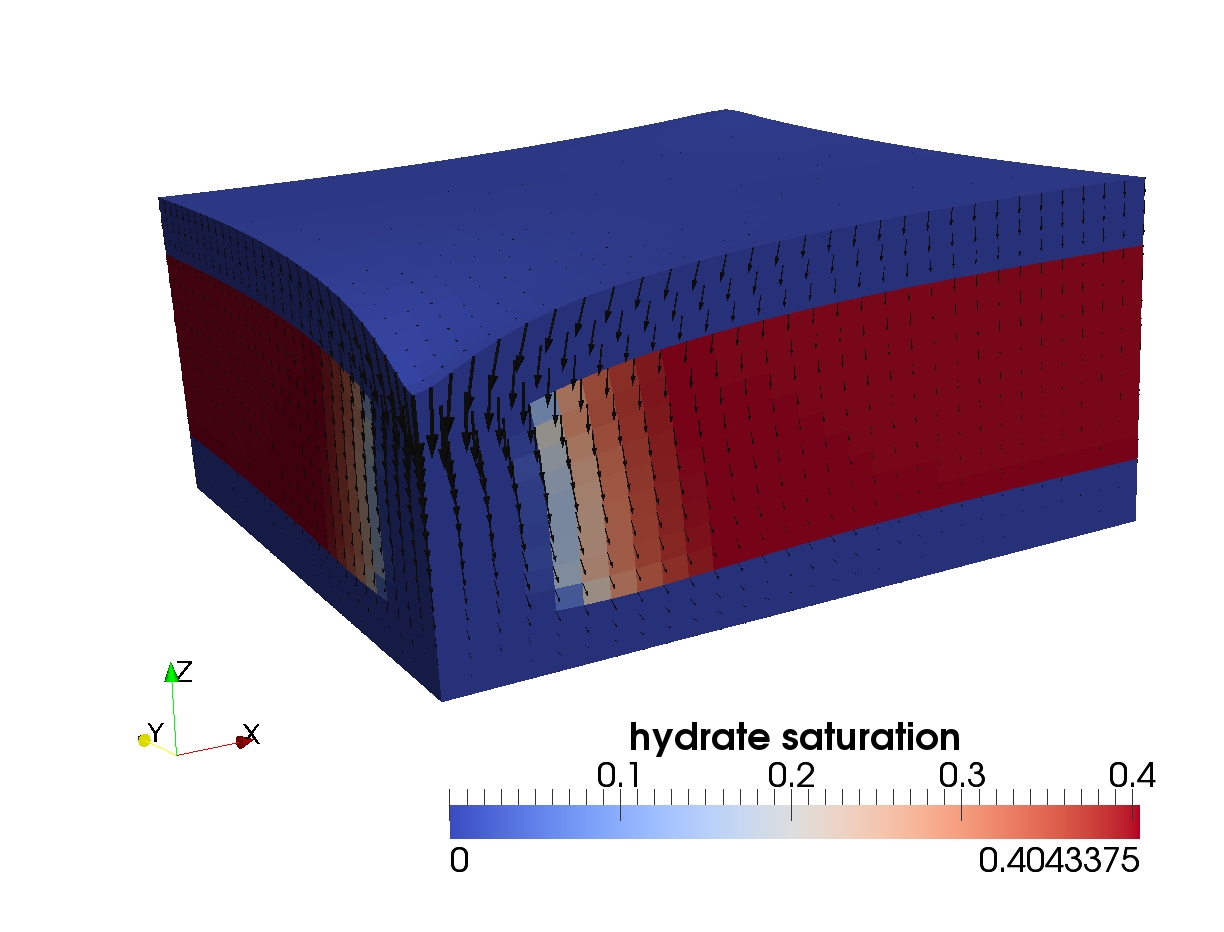}
    \label{fig:3Dtest-ShProfile}}
    \hspace{-0.2cm}
    \subfloat[Gas saturation]{%
      \includegraphics[scale=0.2]{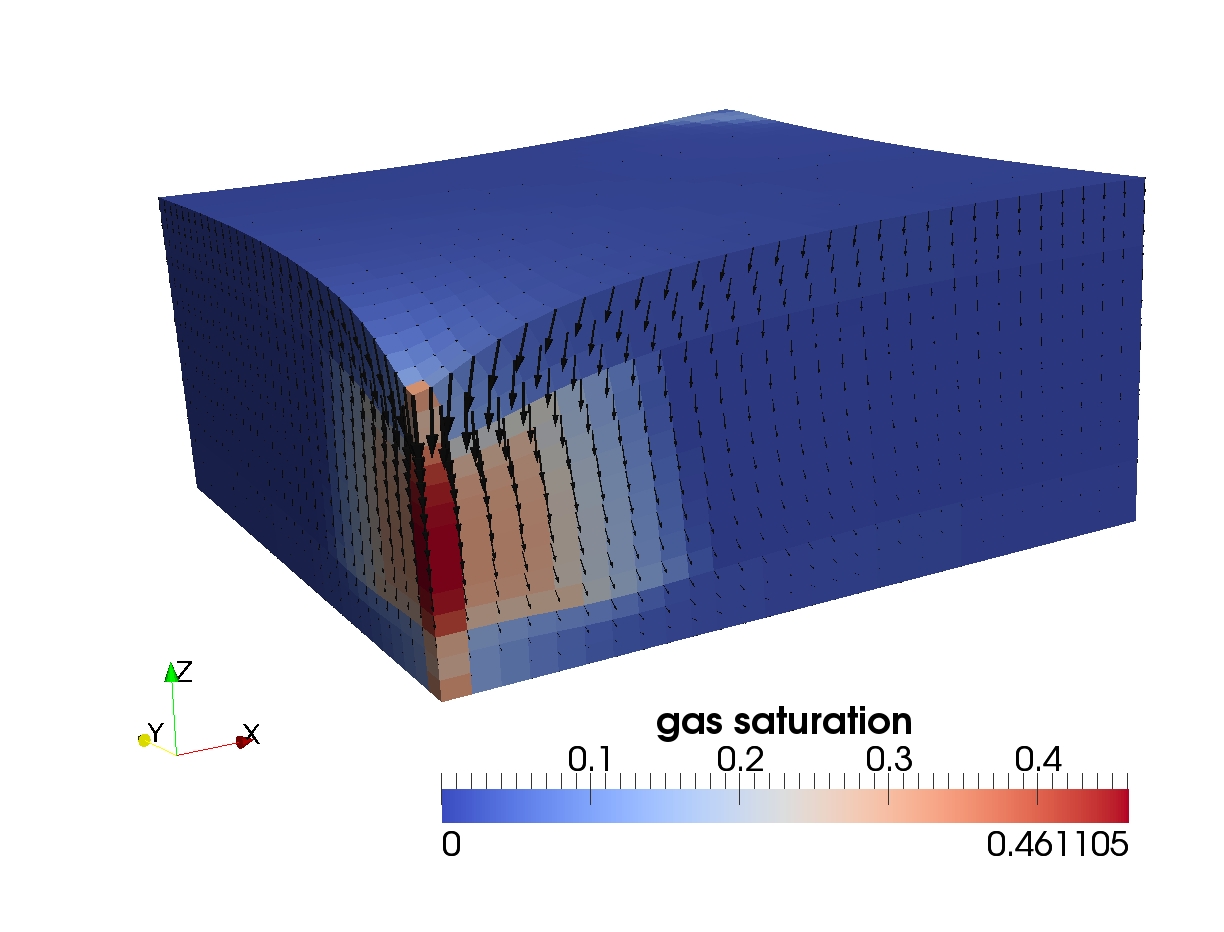}
    \label{fig:3Dtest-SgProfile}}
    
    \vfill
    
    \subfloat[Temperature]{%
    \includegraphics[scale=0.2]{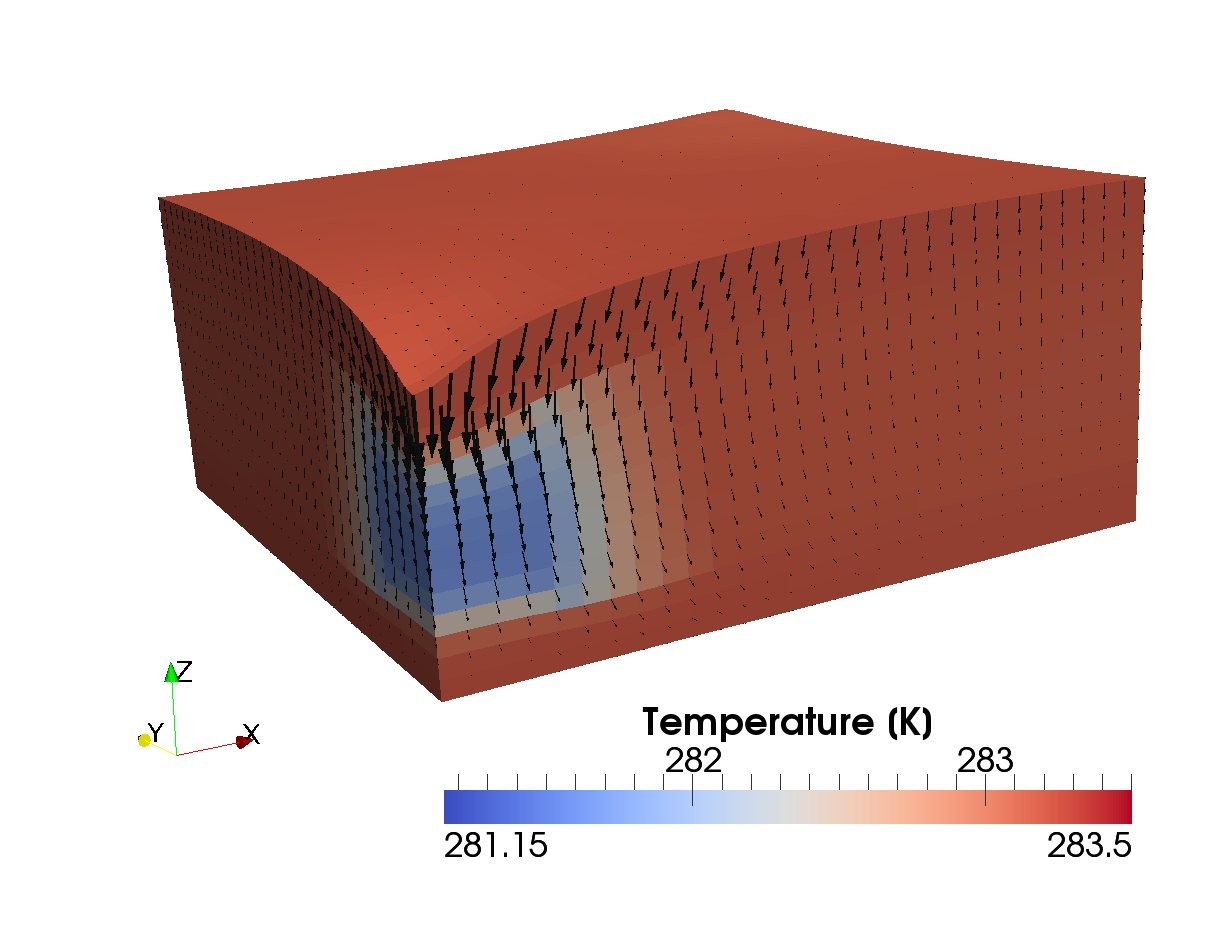}
    \label{fig:3Dtest-TempProfile}}
    \hspace{-0.2cm}
    \subfloat[Deviatoric stress]{%
      \includegraphics[scale=0.2]{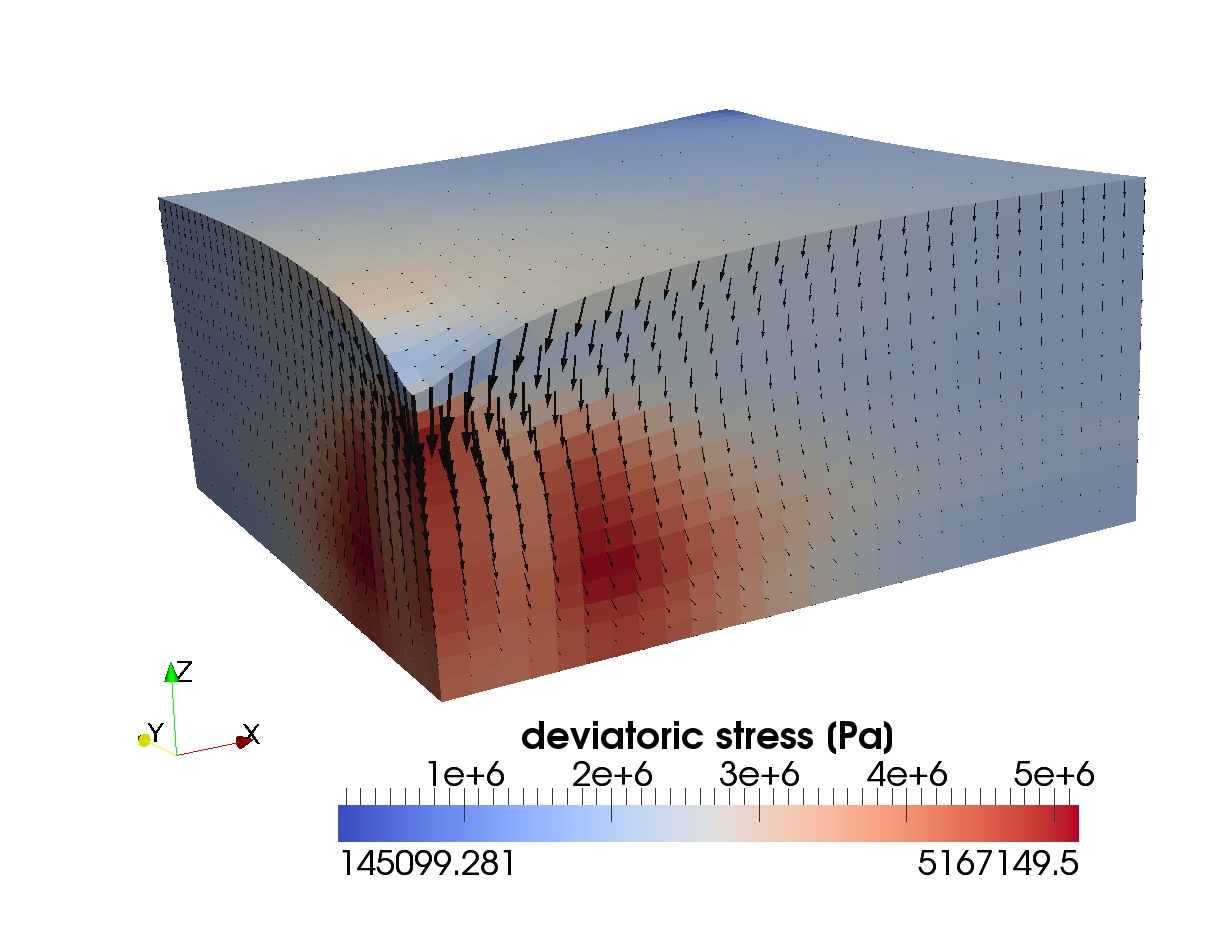}
    \label{fig:3Dtest-DevStressProfile}}
    
    \vfill
    
    \subfloat[Effective porosity]{%
    \includegraphics[scale=0.2]{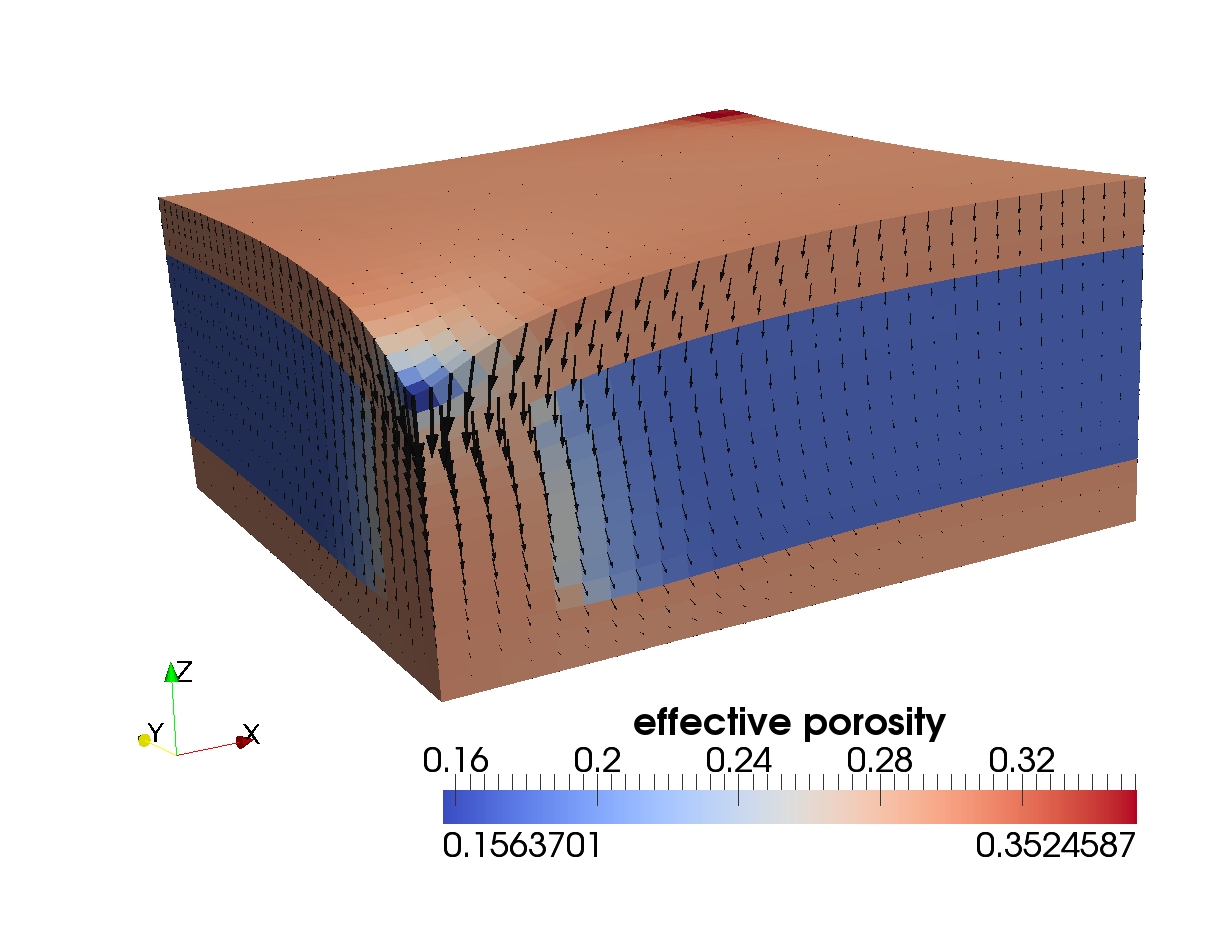}
    \label{fig:3Dtest-PhiEffProfile}}
    \hspace{-0.2cm}
    \subfloat[Intrinsic permeability]{%
      \includegraphics[scale=0.2]{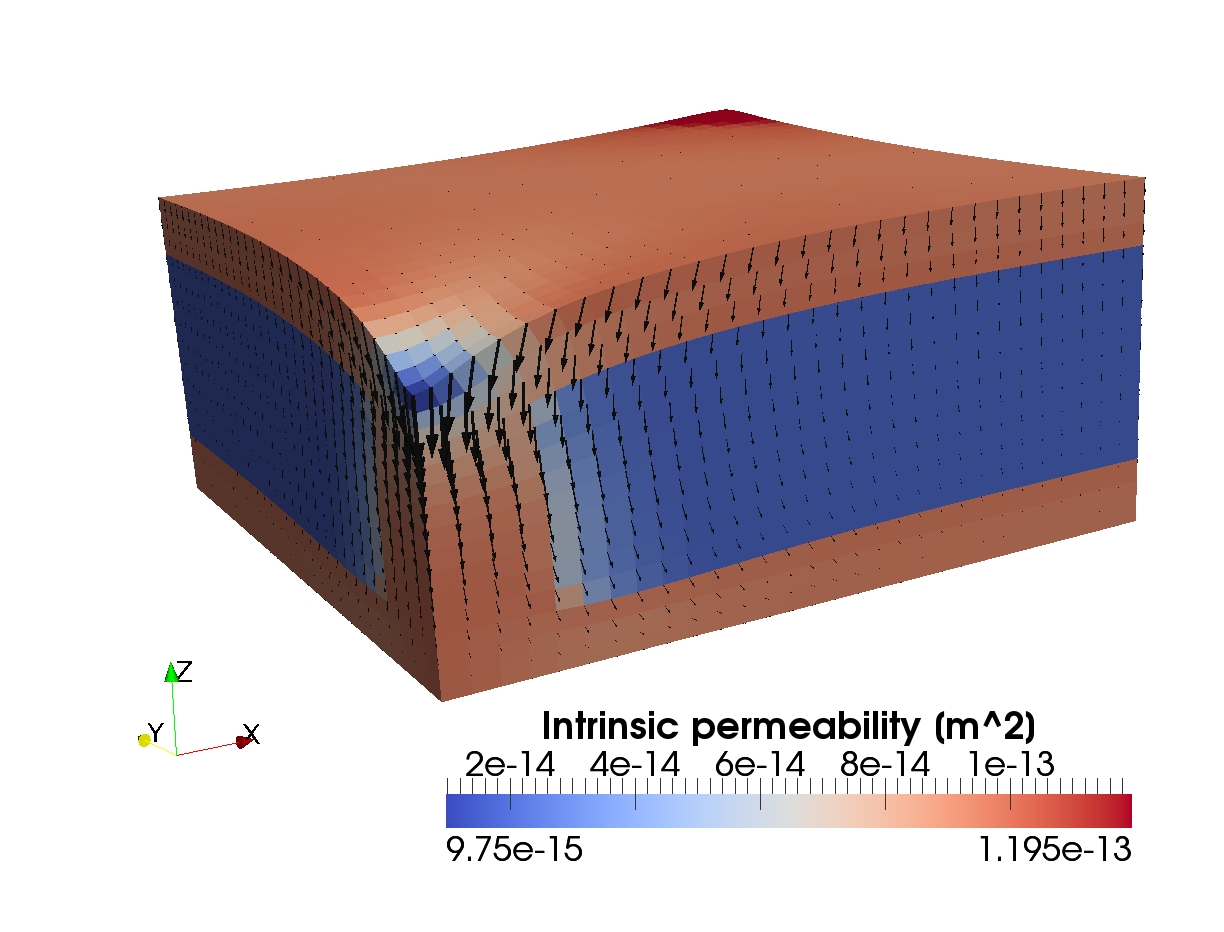}
    \label{fig:3Dtest-PermeabilityProfile}}
    
    \caption{Numerical example - Depressurization of a $3D$ hydrate reservoir \\ Selected profiles at $time=2.5$ hours}
    \label{fig:3Dtest-results}
\end{figure*}

\paragraph{\textbf{Numerical simulation and results}}
The domain is discretized into $30\times30\times15$ cells.
Full hydro-geomechanical model is solved.
To reduce the computational costs the decoupling strategy and iterative solution scheme described in Section \ref{sec:numericalSolution} is used. 
The primary variables being solved for are: gas phase pressure $P_g$, aqueous phase saturation $S_w$, hydrate phase saturation $S_h$, temperature $T$, total porosity $\phi$, and displacements $\mathbf{u}$. 
Some of the other important secondary variables which are calculated as post process include gas saturation $S_g$, effective porosity $\phi_{eff}$, intrinsic permeability $K$, stresses $\tilde \sigma'$, 
strains $\tilde \epsilon$, etc.
The simulation is run until $t_{end}=24\ hrs$ with a time step size of $dt = 200$ s.
Selected profiles showing the state of the reservoir at $t_{end}$ are shown in Fig. \ref{fig:3Dtest-results}.
Fig. \ref{fig:3Dtest-ShProfile} and Fig. \ref{fig:3Dtest-SgProfile} show the melted hydrate and the accumulated gas in the vicinity of the gas well.
Fig. \ref{fig:3Dtest-TempProfile} shows the decrease in temperature due to the endothermic nature of hydrate dissociation.
Fig. \ref{fig:3Dtest-DevStressProfile} shows the stress built up in the region aroung the well where the hydrate is dissociating.
Fig. \ref{fig:3Dtest-PhiEffProfile} and Fig. \ref{fig:3Dtest-PhiEffProfile} show the change in effective porosity and intrinsic permeability as a result of hydrate dissociation and soil deformation.
The vectors in Fig. \ref{fig:3Dtest-results} show the displacements $\mathbf{u}$. 
The domain is warped with respect to displacement to show the ground subsidence around the well clearly. The warping of the domain is achieved through post-processing using PARAView \cite{Paraview}.

\section{Concluding remarks}
\paragraph{} In this article, we have presented a model concept for multi-phase, multi-component flow through deformable methane hydrate reservoirs.
This forms the \textit{core} of our hydrate numerical model and contains only those model components which are necessary to simulate the most important hydro-geomechanical processes observed in a subsurface hydrate reservoir.
The structure of this model is such that the \textit{core} can be modularly extended to enhance each of the model components to the desired level of complexity depending on the application at hand.

\paragraph{}
Our focus so far has been to ensure that the dynamics of the hydro-geomechanical interactions are consistently accounted for in our mathematical model.
We have identified the important physical processes and the \textit{cause-effect} based couplings. 
Based on this we have presented our decoupling strategy, where we have discussed how we breakdown our complex multiphysics model into relatively simpler sub-models.
We have also discussed our solution strategy, which involves first solving the sub-models separately to obtain a decoupled solution, and then reintroducing the couplings iteratively.
Through the numerical examples, each of which isolates an important physical process in a hydrate reservoir, 
we have shown that our model is versatile and is capable of capturing the important couplings effectively.

  \renewcommand{\thesection}{Appendix \Alph{section}}
  \setcounter{section}{0}
  
  \renewcommand{\theequation}{\Alph{section}\arabic{equation}}
  \setcounter{equation}{0}

 \section{: Derivation of storage equation (Refer Eqn. \eqref{eqn:storageEquation})}
 \label{appendix:A}
 \paragraph{} We re-write the mass conservation equations for gas and water phase-wise, i.e., for each fluid-phase $\alpha=g,w$,
\begin{align}\label{eqn:MB_GW}
 & \DDt \left( \phi \rho_{\alpha} S_{\alpha} \right) + \nabla \cdot \left( \phi \rho_{\alpha} S_{\alpha} \myvec{v}_{\alpha,t} \right) = \dot q_{\alpha}
\end{align}
where, $\dot q_{\alpha}$ is the volumetric source term for phase $\alpha$ given as $\dot q_{\alpha} = \sum_{\alpha} \left(\chi_{\alpha}^{\kappa} \ \dot g^{\kappa} \right)$.

\paragraph{}
Expanding the partial derivatives in Eqn. \eqref{eqn:MB_GW} and rearranging gives
\begin{align}\label{eqn:MB_GW_expanded}
 & \phi S_{\alpha} \underbrace{\left( \DDt \rho_{\alpha} + \myvec{v}_{\alpha,t} \cdot \nabla \rho_{\alpha} \right)}_{\dfrac{d}{dt}\textstyle{\rho_{\alpha}}} 
	+ \rho_{\alpha} \DDt \left( \phi S_{\alpha} \right) \notag \\
 & + \rho_{\alpha} \nabla \cdot \left( \phi S_{\alpha} \myvec{v}_{\alpha,t} \right) = \dot q_{\alpha} \notag\\
 \implies & \phi S_{\alpha} \frac{d}{dt}\rho_{\alpha} + \rho_{\alpha} \DDt \left( \phi S_{\alpha} \right) + \rho_{\alpha} \nabla \cdot \left( \phi S_{\alpha} \myvec{v}_{\alpha,t} \right) = \dot q_{\alpha}  \ .
\end{align}

\paragraph{}
The rate of change of the fluid density is defined as $\ \dfrac{d}{dt}\rho_{\alpha} = \dfrac{\rho_{\alpha}}{B_{\alpha}}\dfrac{d}{dt}P_{\alpha} \ , \ $ where, $ B_{\alpha}$ is the fluid-phase bulk modulus.
Using this definition in Eqn. \eqref{eqn:MB_GW_expanded} and dividing by $\rho_{\alpha}$, we get
\begin{align}\label{eqn:MB_GW_reduced}
 & \frac{\phi S_{\alpha}}{B_{\alpha}}\frac{d}{dt}P_{\alpha} + \DDt \left( \phi S_{\alpha} \right) + \nabla \cdot \left( \phi S_{\alpha} \myvec{v}_{\alpha,t} \right) = \frac{\dot q_{\alpha}}{\rho_{\alpha}} \ .
\end{align}

\paragraph{}
Since we assume $P_c=0$, 
the phase pressures are equal. So we drop the subscript and assign the symbol $P$ to the phase pressures.

\paragraph{}
Next, we sum Eqn. \eqref{eqn:MB_GW_reduced} over $\alpha=g,w$, which gives
\begin{align}\label{eqn:SUM_MB_GW_reduced}
 & \phi \left(\frac{S_{w}}{B_{w}} + \frac{S_{g}}{B_{g}} \right)\frac{d}{dt}P + \DDt \left[ \phi \left(S_{w} + S_{g}\right) \right] \notag \\
 & + \nabla \cdot \left( \phi S_{w} \myvec{v}_{w,t} \right) 
 + \nabla \cdot \left( \phi S_{g} \myvec{v}_{g,t} \right)
 = \frac{\dot q_{w}}{\rho_{w}} + \frac{\dot q_{g}}{\rho_{g}} \ .
\end{align}

\paragraph{}
We define the effective fluid-phase saturation $S_{\alpha,e}$ as the volume of fluid phase $\alpha=g,w$ in the \textit{effective pore space} which is characterized by the effective porosity $\phi_{e}$.
So, $\quad S_{\alpha,e} = \dfrac{S_{\alpha}}{1-S_h}\ $, and $\quad \phi_e = \phi (1-S_h)\ $.
Also, by the summation relationship, $\quad S_w + S_g = 1-S_h \ $, or $\quad \dfrac{S_w}{1-S_h} + \dfrac{S_g}{1-S_h} = 1\ $.

\paragraph{}
Substituting these definitions in Eqn. \eqref{eqn:SUM_MB_GW_reduced}, we get Eqn. \eqref{eqn:MB_FLUID} which we will call the \textit{fluid mass balance equation}.
\begin{align}\label{eqn:MB_FLUID}
   & \phi_e \left(\frac{S_{w,e}}{B_{w}} + \frac{S_{g,e}}{B_{g}} \right)\frac{d}{dt}P + \DDt \phi_e  
 + \nabla \cdot \left( \phi_e S_{w,e} \myvec{v}_{w,t} \right) \notag\\
 & + \nabla \cdot \left( \phi_e S_{g,e} \myvec{v}_{g,t} \right)
 = \frac{\dot q_{w}}{\rho_{w}} + \frac{\dot q_{g}}{\rho_{g}} \ .
\end{align}


\paragraph{} For the hydrate and the soil phases, the mass conservation equations (Eqn. \eqref{eqn:massBalance_Hydrate} and Eqn. \eqref{eqn:massBalance_Soil}) described in Section \ref{sec:governingEquations} are used. 
Adding Eqn. \eqref{eqn:massBalance_Hydrate} and Eqn. \eqref{eqn:massBalance_Soil}, we get
\begin{align}\label{eqn:MB_SH}
 & \DDt \left[ \phi S_h \rho_h + \left(1-\phi\right) \rho_s\right] 
 + \nabla \cdot \left[ \phi S_h \rho_h + \left(1-\phi\right) \rho_s\right] \myvec{v}_s  = \dot q_h \notag\\
 & \implies \DDt \left(1-\phi_e\right)\rho_{sh} + \nabla \cdot \left[\left(1-\phi_e\right) \myvec{v}_s\right] = \dot q_h \ .
\end{align}

\paragraph{}
Eqn. \eqref{eqn:MB_SH} is the mass balance relationship for the \textit{hydrate-soil-composite} matrix.
Expanding the derivatives in Eqn. \eqref{eqn:MB_SH} and rearranging gives
\begin{align}\label{eqn:MB_SH_expanded}
& \left(1-\phi_e\right) \underbrace{\left[ \DDt \rho_{sh} + \myvec{v}_s \cdot \nabla \rho_{sh} \right]}_{\dfrac{d}{dt}\textstyle{\rho_{sh}}} 
 - \rho_{sh} \DDt \phi_e \notag \\
& + \rho_{sh} \nabla \cdot \left[\left(1-\phi_e\right)\myvec{v}_s\right] = \dot q_h \notag \\
 \implies  & \left(1-\phi_e\right) \frac{d}{dt} \rho_{sh} - \rho_{sh} \DDt \phi_e + \rho_{sh} \nabla \cdot \left[\left(1-\phi_e\right)\myvec{v}_s\right] \notag \\
 & = \dot q_h  \ .
\end{align}

\paragraph{}
Using the expression for rate of change of density of the composite matrix from Eqn. \eqref{eqn:RateOfChangeOfSolidDensity} in Eqn.\eqref{eqn:MB_SH_expanded}, we get
\begin{align}\label{eqn:PorosityEqn}
 \DDt \phi_e = \frac{1}{B_{sh}} \frac{d}{dt}\sigma - \frac{\phi_e}{B_{sh}} \frac{d}{dt}P + \nabla \cdot \left[\left(1-\phi_e\right)\myvec{v}_s\right] - \frac{\dot q_h}{\rho_{sh}}
\end{align}
where, $B_{sh}$ is bulk modulus of the composite solid.

\paragraph{}
Eqn. \eqref{eqn:PorosityEqn} describes the rate of change of the effective porosity due to external stress $\sigma$ and internal fluid pore-pressure $P$.
Finally, substituting Eqn. \eqref{eqn:PorosityEqn} in Eqn. \eqref{eqn:MB_FLUID}, we obtain
\begin{align}\label{eqn:MB_FLUID+SH}
 &\phi_e \left[ \left(\frac{S_{w,e}}{B_{w}}-\frac{S_{w,e}}{B_{sh}}\right) +  \left(\frac{S_{g,e}}{B_{g}}-\frac{S_{g,e}}{B_{sh}}\right) \right]\frac{d}{dt}P 
 + \frac{1}{B_{sh}} \frac{d}{dt}\sigma \notag\\
 &+ \nabla \cdot \left[ \phi_e S_{w,e} \left( \myvec{v}_{w,t} - \myvec{v}_s \right)  \right] 
 + \nabla \cdot \left[ \phi_e S_{g,e} \left( \myvec{v}_{g,t} - \myvec{v}_s \right) \right] \notag\\
 &+ \nabla \cdot \myvec{v}_s = \frac{\dot q_{w}}{\rho_{w}} + \frac{\dot q_{g}}{\rho_{g}} + \frac{\dot q_h}{\rho_{sh}} \ .
\end{align}


\paragraph{}
Eqn. \eqref{eqn:MB_FLUID+SH} expresses the \textit{total mass balance} for the whole porous medium consisting of phases $\gamma=g,w,h,s$.
We now substitute the simplified constitutive relationships listed in Table \ref{table:constitutiveRelationships} in Eqn. \eqref{eqn:MB_FLUID+SH},
\begin{align}\label{eqn:storageEquationExpanded}
 &\phi_e \left[ \left(\frac{S_{w,e}}{B_{w}}-\frac{S_{w,e}}{B_{sh}}\right) +  \left(\frac{S_{g,e}}{B_{g}}-\frac{S_{g,e}}{B_{sh}}\right) \right]\frac{d}{dt}P 
  + \frac{d}{dt}\epsilon  \notag\\
 &+ \frac{1}{B_{sh}} \left( \frac{d}{dt}\sigma' + \alpha \frac{d}{dt}P \right) 
 - \nabla \cdot \left[ K \left( \frac{k_{r,w}}{\mu_{w}} + \frac{k_{r,g}}{\mu_{g}} \right) \nabla P \right] \notag\\
 & = \left(N_{h} \frac{M_{w}}{\rho_{w}} + \frac{M_{g}}{\rho_{g}} - \frac{M_{h}}{\rho_{sh}}\right) k_0 \ A_{s,0} \ S_h \ \left(P_e - P \right) \notag\\
 & \implies 
  \underbrace{\left[\phi_e \left( \frac{S_{w,e}}{B_{w}} + \frac{S_{g,e}}{B_{g}}\right) + \left(\frac{ \alpha - \phi_e }{B_{sh}} \right) \right]}_{\textstyle{\text{Storativity } S}}\frac{d}{dt}P \notag\\
 & + \underbrace{\left( 1 - \frac{B_{m}}{B_{sh}} \right)}_{\textstyle{\alpha}} \frac{d}{dt}\epsilon 
 - \nabla \cdot \underbrace{\frac{K}{2}\left(\frac{1}{\mu_{w}} + \frac{1}{\mu_{g}}\right)}_{\textstyle{\dfrac{K}{\mu_{f}}}} \nabla P  \notag\\
 & = \underbrace{\left(N_{h} \frac{M_{w}}{\rho_{w}} + \frac{M_{g}}{\rho_{g}} - \frac{M_{h}}{\rho_{sh}}\right) k_0 \ A_{s,0}}_{\textstyle{C}} \ S_h \ \left(P_e - P \right) .
\end{align}

\begin{table}[h]
 \centering
    \captionof{table}{Simplified constitutive relationships}\label{table:constitutiveRelationships}
    \begin{tabular}{L{2.7cm} L{4.8cm} }\hline
      Hydrate reaction 	kinetics 	& $\dot q_g = k_0 \ A_{s} \ M_{g} \ \left(P_{e}-P\right)$	\\ [0.3em]
					& $A_{s} = A_{s,0}\times S_h$ 					\\ [0.5em]
					& $\dot q_w = N_h \dfrac{M_w}{M_g} \ \dot q_g$ , $- \dot q_h = \dfrac{M_h}{M_g} \ \dot q_g$ \\ [1em] \hline
      Darcy velocity 			& $\myvec{v}_{\alpha,r} = -K \dfrac{k_{r,\alpha}}{\mu_{\alpha}} \nabla P$ 			\\
					& $\ k_{r,g} = k_{r,w} = 0.5$ 			\\ [0.5em] \hline
      Effective stress principle 	& $\boldmath{\tilde\sigma} = \boldmath{\tilde\sigma'} + \alpha P \boldmath{\tilde I}$				\\ [-0.5em]
					& isotropic stress $\rightarrow$ $\sigma = \sigma' + \alpha P $ \\ [0.5em] \hline
      Linear elastic stress-strain law	& $\tilde \sigma' = - 2 G_{sh}\ \tilde \epsilon - \lambda_{sh} \left( tr \ \tilde \epsilon \right) \tilde I$ \\ [0.3em]
					& $\tilde \epsilon = \frac{1}{2} \left( \nabla \myvec{u} + \nabla^T \myvec{u} \right) $			\\ [0.75em]
					& isotropic strain $\rightarrow$ $\epsilon = \nabla \cdot \myvec{u} = - \left( B_{sh} - \dfrac{4}{3}G_{sh}\right)^{-1} $ \\
					& $\sigma' = - (1/B_{m})\ \sigma'$ \\\hline
    \end{tabular}
\end{table}

\paragraph{}
Eqn. \eqref{eqn:storageEquationExpanded} can be rewritten in a condensed form as
\begin{align}\label{eqn:storageEquationAppendix}
 \alpha \frac{d}{dt}\epsilon + S \frac{d}{dt}P = \nabla \cdot \frac{K}{\mu_{f}} \nabla P + C \ S_h \ \left(P_e - P \right) \ .
\end{align}
This is the \textit{storage equation} describing the pressure response in a poroelastic hydrate soil. 
 


\begin{thebibliography}{500}

  \bibitem{Ahmadi2004} G. Ahmadi, C. Ji, H.S. Duane, Numerical solution for natural gas production from methane hydrate dissociation, \textit{J. Pet. Sci. Eng.}, \textbf{41}, 169-185, 2004.

  \bibitem{Paraview} J. Ahrens, B. Geveci, C. Law, ParaView: An End-User Tool for Large Data Visualization, Visualization Handbook, \textit{Elsevier}, 2005, ISBN-13: 978-0123875822.
    
  \bibitem{AtkinsPhysicalChemistry} P. Atkins, and J. de Paula, Physical Chemistry, \textit{Oxford University Press}, \textbf{8 ed.}, 2006.

  \bibitem{Yuhu2009} Y. Bai, Q. Li, Y. Zhao, X. Li, Y. Du, The experimental and numerical studies on gas production from hydrate reservoir by depressurization, \textit{Transp. Porous Med.}, \textbf{79}, 443-468, 2009.

  \bibitem{DUNEBastian2007} P. Bastian, M. Blatt, A. Dedner, C. Engwer, R. Kl\"ofkorn, R. Kornhuber, M. Ohlberger, O. Sander, 
  A generic grid interface for parallel and adaptive scientific computing. Part II: implementation and tests in DUNE, preprint, no. 404, DFG Research Center MATHEON, 2007.

  \bibitem{Bear1972} J. Bear, Dynamics of fluids in porous media, Elsevier, New York, 1972. 

  \bibitem{BiotWillis1957} 
    M.A. Biot, and D.G. Willis, The elastic coefficients of the theory of consolidation, \textit{J. Appl. Mech.}, \textbf{24}, 594-601, 1957.

%

  \bibitem{KimBishnoi1987} Kim, H.C., Bishnoi, P.R., Heidemann, R.A., and Rizvi, \emph{Kinetics of Methane Hydrate Decomposition}, Chem. Eng. Sci. 
    \textbf{42} (1987), 1645-1653.


  \bibitem{BrooksCorey1964} R.H. Brooks., and A.T. Corey, Hydraulic properties of porous media, \textit{Hydrology Papers}, Colorado State University, Fort Collins, \textbf{3}, 24, 1964.

  \bibitem{Burdine1953} N.T. Burdine, Relative permeability calculations from pore-size distribution data, Technical report, \textit{Petroleum Transanctions, AIME}, 1953.

  \bibitem{Civan2000} F. Civan, Predictability of Porosity and Permeability Alterations by Geochemical and Geomechanical Rock and Fluid Interactions, 
    \textit{SPE 58746, Proceedings of the SPE International Symposium on Formation Damage}, Lafayette, Louisiana, 23-24 February, 2000.


  \bibitem{ClementHookerSkeen1996} T.P. Clement, B.S. Hooker, R.S. Skeen, Macroscopic models for predicting changes in saturated porous media properties caused by microbial growth, \textit{Ground Water}, \textbf{34}, 934-942, 1996.

  \bibitem{DaweThomas2007}
  R.A. Dawe, S. Thomas, A large potential methane source - natural gas hydrates. Energ. Source Part A 29, 217-229, 2007.

  \bibitem{SuperLU} J.W. Demmel, J.R. Gilbert, X.S. Li, 
  SuperLU Users’ Guide. Technical Report LBNL-44289, \textit{Lawrence Berkeley National Laboratory}, September 1999. http://crd-legacy.lbl.gov/~xiaoye/SuperLU/. Last update: September 2007.
  


  \bibitem{Hassanizadeh1986} M. Hassanizadeh, Derivation of basic equations of mass transport in porous media, part 2: generalized Darcy's and Ficks's laws, \textit{Advances in Water Resources}, \textbf{9}, 207-222, 1986.



  \bibitem{Helmig1997} R. Helmig, Multiphase Flow and Transport Processes in the Subsurface: A Contribution to the Modeling of Hydrosystems, \textit{Springer-Verlag}, Berlin, 1997.

  \bibitem{Himmelblau1964} D. M. Himmelblau, Diffusion of Dissolved Gases in Liquids, \textit{Chem. Rev.}, \textbf{64(5)}, 527-550, 1964.

  \bibitem{Hyodo2005} M. Hyodo, Y. Nakata, N. Yoshimoto, T. Ebinuma, Basic research on the mechanical behaviour of methane hydrate sediments mixture, \textit{Soils. Found.}, \textbf{45(1)}, 75-85, 2005.

  \bibitem{Hyodo2014} M. Hyodo, Y. Li, J. Yoneda, Y. Nakata, N. Yoshimoto, A. Nishimura, Effects of dissociation on the shear strength and deformation behavior of methane hydrate-bearing sediments,
  \textit{Marine and Petroleum Geology}, \textbf{51}, 52-62, 2014.

  \bibitem{Itasca2009} Itasca, $\text{FLAC}^{3D}$, Fast lagrangian analysisof continua in 3 dimensions, Version 4.0, Itasca Consulting Group, Minneapolis, Minnesota, 438pp, 2009.


  \bibitem{KamathHolder1987} V.A. Kamath, G.D. Holder, \textit{Dissociation Heat Transfer Characteristics of Methane Hydrates}, AIChE J. \textbf{33} (1987), 347-350.


  \bibitem{Kimoto2010} S. Kimoto,F. Oka, T. Fushita, A chemo-thermo-mechanically coupled analysis of ground deformation induced by gas hydrate dissociation, 
    \textit{International Journal of Mechanical Sciences}, \textbf{52(2)}, 365-376, 2010.

  \bibitem{KlarSogaNg2010} A. Klar, K. Soga, and Y.A. Ng, Coupled deformation-flow analysis for methane hydrate extraction, \textit{G\'eotechnique}, \textbf{60}, 765-776, 2010.

  \bibitem{Kolditz2012} O. Kolditz, U.J. G{\"o}rke, H. Shao, and W. Wang, Thermo-Hydro-Mechanical-Chemical Processes in Porous Media: Benchmarks and Examples,
  \textit{Lecture Notes in Computational Science and Engineering, Springer Berlin Heidelberg}, 2012, ISBN:9783642271779.
  

  \bibitem{Lee2010} J. Y. Lee, F.M. Francisca, J.C. Santamarina, and C. Ruppel, 
  Parametric study of the physical properties of hydrate-bearing sand, silt and clay sediments: 2. small-strain mechanical properties, 
  \textit{J. Geophys. Res.}, \textbf{115}, B11105, doi:10.1029/2009JB006670, 2010. 

  \bibitem{LewisSchrefler2000} R. W. Lewis, and B. A. Schrefler, The finite element method in the static and dynamic deformation and consolidation of porous media, 2nd ed., \textit{John Wiley and Sons}, England, 2000.

  \bibitem{LiuFlemming2007} X. Liu, P.B. Flemings, Dynamic multiphase flow model of hydrate formation in marine sediments, \textit{J. Geophys. Res.}, \textbf{112}, B03101, 2007.


  \bibitem{Masui2006}
    A. Masui, H. Haneda, Y. Ogata, and K. Aoki, Triaxial compression test on submarine sediment containing methane hydrate in deep sea off the coast of Japan (in Japanese), 
    paper presented at the \textit{41st Annual Conference, Jpn. Geotech. Soc.}, Kagoshima, Japan, 12-14 July, (2006).


  \bibitem{Milkov2004}
   A.V. Milkov, Global estimates of hydrate-bound gas in marine sediments: how much is really out there?, Earth-Sci. Rev. 66, 183-197, 2004.

  \bibitem{Miyazaki2011}
    K. Miyazaki, A. Masui, Y. Sakamoto, K. Aoki, N. Tenma, and T. Yamaguchi, Triaxial compressive properties of artificial methane-hydrate-bearing sediment, 
    \textit{J. Geophys. Res.}, \textbf{116}, B06102, doi:10.1029/2010JB008049, 2011.
%

  \bibitem{Moridis2003} G.J. Moridis, Numerical studies of gas productionfrom methane hydrates, \textit{SPE 87330, SPE J.}, \textbf{32}, No.8, 359-370, 2003.

  \bibitem{Moridis2007} G.J. Moridis, M. Kowalsky, K. Pruess, Depressurization-inducedgas production from class-I hydrate deposits, \textit{SPE J. Reservoir Evaluation Engng}, \textbf{10}, No.5, 458-481, 2007.

  \bibitem{Mualem1976} Y. Mualem, A new model for predicting the hydraulic conductivity of unsaturated porous media, \textit{Water Resour. Res.}, (1976), \textbf{12(3)}, 513-522 p.
  
  \bibitem{PengRobinson1970} a) D.Y. Peng, D.B. Robinson, A New Two-Constant Equation of State, \textit{4th Intl. Heat Transfer Conf.}, Paris-Versailles, "Heat Transfer 1970", \textbf {4}, paper B.7.b., 1970\newline
			     b) Y.S. Wei, R.J. Sadus, Equations of State for the Calculation of Fluid-Phase Equilibria, \textit{AIChE Journal}, \textbf{46}, 169-196 p, 2000.


  \bibitem{Rockhold2001} M.L. Rockhold, R.R. Yarwood, M.R. Niemet, P.J. Bottomley, J.S. Selker, Considerations for modelling bacterial-induced changes in hydraulic properties of variably sarurated porous media, 
    \textit{Advances in Water Resources}, \textbf{25}, 477-495, 2002.


  \bibitem{RutqvistMoridis2008} J. Rutqvist, G.J. Moridis, Development of a numerical simulator for analyzing the geomechanical performance of hydrate-bearing sediments, 
    \textit{Proceedings of the 42nd U.S. Rock Mechanics Symposium}, San Francisco, California, USA, June 29-July 2, 2008, American Rock Mechanics Association, ARMA, Paper No. 139. 


  \bibitem{NISTChemHenryConst} R. Sander, Henry's Law Constants, \textit{NIST Chemistry WebBook, 
  NIST Standard Reference Database Number 69, Eds. P.J. Linstrom and W.G. Mallard}, National Institute of Standards and Technology, Gaithersburg MD, 20899, http://webbook.nist.gov, (retrieved February 23, 2015). 

  \bibitem{SantamarinaRuppel2010}
    J.C. Santamarina, and C. Ruppel, The impact of hydrate saturation on the mechanical, electrical, and thermal properties of hydrate-bearing sand, silts, and clay, 
    in \textit{Geophysical Characterization of Gas Hydrates, Geophys. Dev. Ser.}, edited by M. Riedel, C. Willoughby, and S. Chopra, Soc. of Explor. Geophys., Tulsa, Okla, \textbf{14}, 373–384, 2010.

  \bibitem{Sloan1998}
  E.D. Sloan, Gas hydrates: review of physical/chemical properties. \textit{Energ. Fuel.}, \textbf{12}, 191-196, 1998.


  \bibitem{Soga2006}
    K. Soga, S. L. Lee, M. Y. A. Ng, and A. Klar, Characterisation and engineering properties of methane hydrate soils, 
    in \textit{Characterisation and Engineering Properties of Natural Soils}, edited by T. S. Tan et al., A. A. Balkema, Lisse, Netherlands, 2591–2642, 2006.



  \bibitem{StatteryBird1958} J.C. Stattery, R.B. Bird, Calculation of diffusion coefficient of dilute gases and of the self diffusion coefficient of dense gases, \textit{AIChE J.}, \textbf{4(2)}, 137-142, June 1958.

  \bibitem{SultanCochonat2004} a) N.C. Sultan,P.Canals,M.Cattaneo,A.Dennielou,B. Haflidason, H. Laberg, J.S. Long, D. Mienert, J. Trincardi, F. Urgeles, R.Vorren,T.O., C. Wilson,
  Triggering mechanisms of slope instability processes and sediment failures on continental margins: a geotechnical approach, \textit{Marine Geology}, \textbf{213(1-4)}, 291-321, 2004a.
  b) N.C. Sultan, P. Foucher, J.P., J. Mienert, Effectofgashydratesmeltingonseafloorslopeinstability, \textit{Marine Geology}, \textbf{213(1)}, 379-401, 2004b.

  \bibitem{SunMohanty2005} F.X. Sun, K.K. Mohanty, N. Nanchary, 1-D modelling of hydrate depressurization in porous media, \textit{Transp. Porous Med.}, \textbf{58}, 315-338, 2005.

  \bibitem{SunMohanty2006} F.X. Sun, K.K. Mohanty, Kinetic simulation of methane hydrate formation and dissociation in porous media, \textit{Chemical Engineering Science}, \textbf{61}, 3476-3495, 2006.

  \bibitem{Tang2007} L.G. Tang, X.S. Li, Z.P. Feng., et. al., Control mechanisms for gas hydrate production by depressurization in different scale hydrate reservoirs, \textit{Energy Fuels}, \textbf{21(1)}, 227-233, 2007. 

%
%
  \bibitem{Terzaghi1925} K. Terzaghi, Erdbaumechanik auf bodenphysikalischer Grundlage, Deuticke, Wien, 1925.

  \bibitem{Tsypkin1991}
G.G. Tsypkin, Effect of liquid phase mobility on gas hydrate dissociation in reservoirs. Izvestiya Akad. Nauk SSSR. Mekh. Zhidkosti i Gaza. 4: 105-114, 1991.


  \bibitem{Verruijt2013} A. Verruijt, Theory and problems of poroelasticity, 2013.

  \bibitem{Waite2009} W. F. Waite, et al., Physical properties of hydrate-bearing sediments, Rev. Geophys., 47, RG4003, doi:10.1029/2008RG000279, 2009.


  \bibitem{Yousif1991} M.H. Yousif, \textit{Experimental and theoretical investigation of methane-gas-hydrate dissociation in porous media}, SPE Reserv. Eng. (1991), 69-76.

  \bibitem{TOUGH-HYDRATEmanual} G.J. Moridis, M.B. Kowalsky, K. Pruess, TOUGH-Fx/HYDRATE v1.0 User's manual: A code for the simulation of system behaviour in hydrate-bearing geologic media, 
  Earth Sciences Division, Lawrence Berkeley National Laboratory, Berkeley, California, March 2008.

  \bibitem{MH21-HYDRES} M. Kurihara, A. Sato, H. Ouchi, H. Narita, Y. Masuda, T. Saeki, T. Fujii, Prediction of gas productivity from Eastern Nankai Trough methane-hydrate reservoirs,
  \textit{SPE125481, SPE Reservoir Eval. Eng.}, \textbf{12}, 477-499, 2009.

  \bibitem{STOMP-HYD} M.D. White, B.P. McGrail, STOMP-HYD: A new numerical simulator for analysis of methane hydrate production from geologic formations, 
  \textit{Proceedings of 2nd International Symposium on Gas Hydrate Technology}, 1-2 November 2006, KIGAM, Daejeon, Korea.
  
  \bibitem{UMSICHT-HyRes} G. Janicki, S. Schlüter, To. Hennig, H. Lyko, G. Deerberg, Simulation of Methane Recovery from Gas Hydrates Combined with Storing Carbon Dioxide as Hydrates, 
  \textit{Journal of Geological Research}, vol. 2011, Article ID 462156, 15 pages, 2011. doi:10.1155/2011/462156

  
\end{thebibliography}
\end{document}